\documentclass[10pt]{article}
\usepackage[utf8]{inputenc}
\usepackage[a4paper,margin=1.2in]{geometry}

\title{Non existence and strong ill-posedness in $C^k$ and Sobolev spaces for SQG}
\author{Diego C\'ordoba\footnote{dcg@icmat.es}\quad and Luis Mart\'inez-Zoroa\footnote{luis.martinez@icmat.es}\\ \\ \small Instituto de Ciencias Matem\'aticas CSIC-UAM-UCM-UC3M }

\usepackage{graphicx}
\usepackage{amsmath}
\usepackage{ dsfont }
\usepackage{amsfonts}
\usepackage{amsthm}

\newtheorem{theorem}{Theorem}[section]
\newtheorem{corollary}{Corollary}[theorem]
\newtheorem{lemma}[theorem]{Lemma}
\newtheorem{remark}{Remark}
\newtheorem{definition}{Definition}

\begin{document}

\maketitle

\begin{abstract}
    We construct solutions in $\mathds{R}^2$ with finite energy of the surface quasi-geostrophic equations (SQG) that initially are in $C^k$ ($k\geq 2$) but that are not in $C^{k}$ for $t>0$. We prove a similar result also for $H^{s}$ in the range $s\in(\frac32,2)$. Moreover, we prove strong ill-posedness in the critical space $H^{2}$. 
\end{abstract}

\tableofcontents

\newpage

\section{Introduction}

We say a function $\theta(x,t):\mathds{R}^2\times \mathds{R}_+\rightarrow \mathds{R}$ is a solution to the SQG equation with initial conditions $\theta(x,0)=\theta_{0}(x)$ if the equation

\begin{eqnarray}
\frac{\partial \theta}{\partial t} + v_{1}\frac{\partial \theta}{\partial x_{1}} + v_{2}\frac{\partial \theta}{\partial x_{2}}= 0 \label{SQG}
\end{eqnarray}
is fulfilled for every $x\in\mathds{R}^2$ and the derivatives exists for every $x\in \mathds{R}^2$. The velocity field $v=(v_1,v_2)$ is defined by
$$v_{1}=-\frac{\partial}{\partial x_{2}}\Lambda^{-1} \theta=-\mathcal{R}_{2}\theta$$
$$v_{2}=\frac{\partial}{\partial x_{1}}\Lambda^{-1} \theta=\mathcal{R}_{1}\theta$$
where $\mathcal{R}_{i}$ are the Riesz transforms in 2 dimensions, with the integral expression
$$\mathcal{R}_{j}\theta=\frac{\Gamma (3/2)}{\pi^{3/2}}P.V. \int_{\mathds{R}^2} \frac{(x_{j}-y_{j})\theta(y)}{|x-y|^{3}}dy_{1}dy_{2} $$
for $j=1,2.$ We denote $\Lambda^{\alpha} f\equiv (-\Delta)^{\frac{\alpha}{2}} f$ by the Fourier transform $\widehat{\Lambda^{\alpha} f}(\xi) = |\xi|^{\alpha}\widehat{f} (\xi)$. 

This model arises in a geophysical fluid dynamics context (see \cite{Held} and \cite{Ped}) and its mathematical analysis was initially treated by Constantin, Majda and Tabak in \cite{Majda} motivated by the number of traits it shares with 3-D incompressible Euler system, where they already established local existence in $H^{s}$ (see also \cite{const1} for bounded domains) and in the case of  $C^{k,\alpha}$ ($k\geq 1$ and $1>\alpha>0$) see \cite{Wu} by Wu. In the critical Sobolev space $H^2$ Chae and Wu \cite{ChaeWu} proved local existence for a logarithmic inviscid regularization of SQG (see also \cite{Jolly}).  Finite time formation of singularities for smooth initial data with finite energy remains an open problem for both SQG and 3-D incompressible Euler equations. 

Due to incompressibility and the transport structure of SQG the $L^p$ ($1\leq p\leq \infty$) norms of the scalar $\theta$ and the $L^2$ norm of the velocity field $v=(v_1,v_2)$ (kinetic energy) are conserved quantities of the system (\ref{SQG}) for sufficiently regular solutions. Global existence of weak solutions in $L^2$ was proven by Resnick in \cite{Resnic} (see also \cite{const2} in the case of  bounded domains) and extended by Marchand in \cite{Marchand} to the class of initial data in $L^p$ with $p>\frac43$. However non-uniqueness of weak solutions was obtained by Buckmaster, Shkoller and Vicol in \cite{Buckmaster} for solutions in $\Lambda^{-1}\theta\in C_t^{\sigma}C_x^{\beta}$ with $\frac12<\beta<\frac45$ and $\sigma < \frac{\beta}{2-\beta}$. 

One of the main objectives of this paper is to construct solutions  in $\mathds{R}^2$ of SQG  that initially are in $C^k\cap L^2$ ($k\geq 2$) but are not in $C^{k}$ for $t>0$. Note that if we consider a velocity field $v(\theta)=\nabla^{\perp}\Lambda^{-(1+\epsilon)}(\theta)$ with $\epsilon>0$, then we have local existence in $C^{k}$ for (\ref{SQG}). 
We also prove strong ill-posedness  in $H^{s}$ for critical and supercritical spaces in the range $s\in(\frac32,2]$. Moreover we construct solutions that are initially in $H^{s}$ for $s\in(\frac32,2)$ but are not in $H^{s}$ for $t>0$, and that are unique in a certain sense that we will specify later. For the SQG equation, there were no strong ill-posedness results in $H^{s}$ and $C^{k}$ prior to the ones obtained in this paper. There are ill-posedness results for active scalars with more singular velocities obtained by Kukavica, Vicol and Wang in \cite{Kukavica}  and, in the case of SQG, in \cite{Elgindi} Elgindi and Masmoudi a mild ill-posedness result is obtained for perturbations of a stationary solution. This, however, does not imply mild or strong ill-posedness for SQG. For more details about this as well as the specific definitions of mild and strong ill-posedness, see subsection \ref{illposedness}  below. A few days after our result appeared on the arXiv, Jeong and Kim \cite{Injee} posted an article on the arXiv with a similar result to the one we have for the critical space $H^{2}$.

 There are some remarkable results regarding norm growth in the periodic setting for SQG. Kiselev and Nazarov \cite{Nazarov} showed that there exists initial conditions with arbitrarily small norm in $H^s$ ($s\geq 11$)  that become large after a long period of time. Recently, He and Kiselev proved in \cite{SmallscaleSQG} an exponential in time growth for the $C^2$ norm
 $$sup_{t\leq T} |\nabla^2\theta|_{L^{\infty}}\geq \exp{\gamma T} \quad\quad \text{for $\gamma(\theta_0)>0$}.$$
 
 On the other hand numerical simulations suggested the existence of solutions with very fast growth of $|\nabla \theta|$ starting with a smooth profile by a collapsing hyperbolic saddle scenario (see \cite{Majda}, \cite{oy} and \cite{cls}). Such a scenario cannot developed a singularity as shown analytically in \cite{c} and \cite{cf}, where a double exponential bound on $|\nabla \theta|$ is obtained. A different blow-up scenario was proposed in \cite{scott} where the fast growth of $|\nabla \theta|$ is associated to a cascade of filament instabilities.

\subsection{The main theorems}

In this paper we prove the following results:

\begin{theorem}\label{illposckintro} (Strong ill-posedness in $C^k$)
For any $c_{0}>0$, $M>0$, $2\leq k\in \mathds{N}$ and $t_{*}>0$, we can find a function $\theta_{0}(x)\in H^{k+\frac14}\cap C^k$  with $||\theta_{0}(x)||_{C^k}\leq c_{0}$ such that the unique solution $\theta(x,t)\in H^{k+\frac14}$ to the SQG equation (\ref{SQG}) with initial  conditions $\theta_{0}(x)$ satisfies $||\theta(x,t_{*})||_{C^{k}}\geq M  c_{0}$.
\end{theorem}

\begin{theorem}\label{noexistenciackintro} (Non existence in $C^k$)
Given $c_{0}>0$, $t_{*}>0$ and $2\leq k\in \mathds{N}$, there are initial conditions $\theta_{0}\in H^{k+1/8}\cap C^k$ for the SQG equation (\ref{SQG}) such that $||\theta_{0}||_{C^k}\leq c_{0}$ and the unique solution $\theta(x,t)\in H^{k+1/8}$ exists and satisfies that $||\theta(x,t)||_{C^k}=\infty$ for all $t\in(0,t_{*}]$.
\end{theorem}
In fact, for the initial conditions given by theorem \ref{noexistenciackintro} there cannot be a solution $\theta(x,t)\in L^{\infty}_{t} L^{2}_{x}$ to (\ref{SQG}) with those initial conditions 
and $||\theta(x,t)||_{C^{k}}\leq M(t)$, $M(t):\mathds{R}_{+}\rightarrow \mathds{R}_{+}$, even if we allow for $||M(t)||_{L^{\infty}}=\infty$. For more details see remark \ref{remarkunic2} after theorem \ref{nonexc2}.

\begin{theorem}\label{illposednessintro} (Strong ill-posedness in $H^{s}$)
For any $c_{0}>0$, $M>0$, $s\in (\frac32, 2]$ and $t_{*}>0$, we can find a $H^{\beta}$ function $\theta_{0}(x)$  with $||\theta_{0}(x)||_{H^{s}}\leq c_{0}$ such that the only solution $\theta(x,t)\in H^{\beta}$, with $\beta(s)>2$ to the SQG equation (\ref{SQG}) with initial conditions $\theta_{0}(x)$ satisfies $||\theta(x,t_{*})||_{H^{s}}\geq M  c_{0}$.
\end{theorem}

\begin{remark}
The purpose of this paper is not to obtain the optimal range of Sobolev spaces in which strong ill-posedness is achieved. There are refinements to the methods used in theorem \ref{illposednessintro} that would allow us to decrease the lower bound in the interval of ill-posedness. 
\end{remark}

\begin{theorem}\label{nonexsuperintro}
(Non existence in $H^{s}$ in the supercritical case)
For any $t_{*}$, $c_{0}>0$ and $s\in(\frac32,2)$ we can find initial conditions $\theta_{0}(x)$,  with $||\theta_{0}(x)||_{H^{s}}\leq c_{0}$ such that there exists a solution $\theta(x,t)$ to (\ref{SQG}) with $\theta(x,0)=\theta_{0}(x)$ satisfying $||\theta(x,t)||_{H^{s}}=\infty$ for all $t\in(0,t_{*}]$. Furthermore, it is the only  solution with initial conditions $\theta_{0}(x)$ such that $\theta(x,t)\in L^{\infty}_{t}C^{\alpha_{1}}_{x}\cap L^{\infty}_{t}L^{2}_{x}$ ($0<\alpha_{1}<\frac12$) with the property that $||\theta(x,t)||_{H^{\alpha_{2}}}\leq M(t)$ ($1<\alpha_{2}\leq \frac32$) for some function $M(t)$.
\end{theorem}

\begin{theorem}\label{nonuniformex}(Non uniform existence in $H^{2}$)
For any $c_{0}>0$ there exist initial conditions $\theta(x,0)$ with $||\theta(x,0)||_{H^{2}}\leq c_{0}$ such that any solution $\theta(x,t)$ to (\ref{SQG}) satisfies

$$\text{ess-sup}_{t\in[0,\epsilon]}||\theta(x,t)||_{H^{2}}=\infty$$
for any $\epsilon>0$.
\end{theorem}

The proof of theorems \ref{nonexsuperintro} and \ref{nonuniformex} can be adapted to work in the critical spaces $W^{1+\frac2p,p}$, $p\in(1,\infty]$,  but we will not go into detail since that is not the goal of the paper. For more information regarding the necessary changes to adapt the proof for these cases, see remark \ref{remarkfinal} after theorem \ref{teoremacritnon}.

\subsection{The strategy of the proof}



%

Ill posedness in critical spaces for  the incompressible Euler equations was already considered in papers by Bourgain and Li (see \cite{Bourgainsobolev} and \cite{Bourgaincm}) obtaining strong ill-posedness for the velocity in the 2D and 3D Euler equations in $C^{k}$, $k\geq 1$ and for the vorticity in the space $H^{d/2}$ (with $d$ the dimension). In fact, they obtained stronger results, in \cite{Bourgaincm} they obtain a velocity $u$ satisfying that, for $0<t_0\leq 1$
$$ \text{ess-sup}_{0<t<t_0}||u(t,\cdot)||_{C^k}=\infty,$$ 
$$||u(0,\cdot)||_{C^k}\leq c_{0}$$
and in \cite{Bourgainsobolev} 
the vorticity $\omega$ satisfies
$$ \text{ess-sup}_{0<t<t_0}||\omega(t,\cdot)||_{\dot {H}^{\frac{d}{2}}}=\infty,$$
$$||\omega(0,\cdot)||_{\dot {H}^{\frac{d}{2}}}\leq c_{0}.$$
Later, analogous results were obtained by Elgindi and Masmoudi in \cite{Elgindi} and Elgindi and Jeong in \cite{Elgindisobolev} with a different approach. Recently, Kwon proved in \cite{Kwon} that there is still strong ill-posedness in $H^1$ for a  regularized version of the 2D incompressible Euler equations. 

Our strategy in this paper for proving strong ill-posedness for SQG differs from the previous works mentioned above since there is no global existence result for SQG in $H^s$. More precisely for theorems \ref{illposckintro}, \ref{noexistenciackintro}, \ref{illposednessintro} and \ref{nonexsuperintro}, we construct solutions by perturbing radial stationary solutions $\theta=\theta(r)$ and,  in order to obtain precise bounds of the errors, we consider an explicit in time pseudo-solution of SQG. We say that a function $\bar{\theta}$ is a pseudo-solution to the SQG equation if it fulfils the evolution equation with  an appropiate small source term (for a more precise definition see section 2.2 below). Namely to prove strong ill-posedness in $C^k$ we will use the following family of pseudo-solutions in the time interval $t\in[0, T]$
\begin{align*}
    \bar{\theta}_{\lambda,J,N}(r,\alpha,t):=&\lambda f_{1}(r)\\
    &+\lambda f_{2}(N^{1/2}(r-1)+1)\sum_{j=1}^{J}\frac{sin(Nj\alpha-\lambda tNj\frac{v_{\alpha}( f_{1})}{r}-\lambda C_{0}t- \frac{\pi}{2}j)}{N^kj^{k+1}},
\end{align*}  
where $(r, \alpha)$ are the polar coordinates,  $f_i$ are smooth compactly supported radial functions, $v_{\alpha}(f_1)$ is the angular velocity generated by the function $f_1$, the parameters fulfil $\lambda,J,N\in(\mathds{R}_+$, $\mathds{N}$, $\mathds{N})$ and $C_0$ is a constant that arises from the velocity operator. This $\bar{\theta}_{\lambda,J,N}$ fulfills the evolution equation

$$\frac{\partial  \bar{\theta}_{\lambda,J,N}}{\partial t}+\frac{\partial \bar{\theta}_{\lambda,J,N}}{\partial \alpha}\frac{v_{\alpha}(\lambda f_{1})}{r}+\lambda C_{0} H(\bar{\theta}_{\lambda,J,N})=0$$
where $H$ is the Hilbert transform with respect to the $\alpha$ variable. The ill-posedness arises from the unboundedness of the operator $H$ in the $C^{k}\cap L^{2}$ spaces. Note however that the appearance of an unbounded operator in our evolution equation does not imply directly ill-posedness, since for example in the Burger-Hilbert's equation

$$\frac{\partial f}{\partial t}+f\frac{\partial f}{\partial x}+H(f)=0,$$
although the $L^{\infty}$ norm has a fast growth (see \cite{growthc1delta}) as long as the solution is $C^{1,\delta}$, Bressan and Nguyen \cite{Bressan} proved the surprising result of global existence in $L^{2}\cap L^{\infty}$.

We denote $\theta_{\lambda,J,N}(r,\alpha, t)$ to be the unique $H^{k+\frac14}$ solution 
of (\ref{SQG}) satisfying initially 
$$\theta_{\lambda,J,N}(r,\alpha, 0)= \bar{\theta}_{\lambda,J,N}(r,\alpha,0).$$ 
We will prove that that for sufficiently large N we have
$$||\theta_{\lambda,J,N}(r,\alpha, t)- \bar{\theta}_{\lambda,J,N}(r,\alpha,t)||_{H^k}\leq C t N^{-(\frac14 + a(k))}$$
with $a(k)>0$ and the constant $C$ depends on the parameters $(\lambda,J,k,T)$. With this bound and the properties of the pseudo-solution we obtain
$$||\theta_{\lambda,J,N}(r,\alpha, t)||_{C^k}\geq \tilde{C}\lambda^2 ln(J)t$$
where $\tilde{C}$ is a universal constant.

Once we have solutions with arbitrary large growth in norm we prove non existence of solutions in $C^k$ by considering the following initial conditions  
$$\theta(x, 0)=\sum_{n\in\mathds{N}}T_{R_{n}}(\bar{\theta}_{\lambda_n,J_n,N_n}(x,0))$$
with $T_{R}(f(x_{1},x_{2}))=f(x_{1}+R,x_{2})$. By choosing appropriately the parameters $(\lambda_{n})_{n\in\mathds{N}}$, $(K_{n})_{n\in\mathds{N}}$, $(N_{n})_{n\in\mathds{N}}$ and $(R_{n})_{n\in\mathds{N}}$ we can show that the unique solution $\theta(x,t)\in H^{k + \frac18}$ with this initial data will leave $C^k$ instantly. In particular the solution $\theta(x,t)$ is not in $C^k$ for any time $t\in (0,T]$.

In the case of strong ill-posedness in Sobolev spaces, theorem \ref{illposednessintro}, we will use a  similar strategy in the range below the critical exponent $s=2$, although the proofs are more involved since we do not have any existence result for the supercritical Sobolev spaces. However, in the critical case (theorem \ref{nonuniformex}) it is not clear that a suitable  pseudo-solution could be constructed by perturbing a radial solution. In order to overcome this obstacle we need a different strategy. In this case our initial data is similar to the one consider in \cite{Bourgainsobolev} with the following expression
$$\theta_{c,J,b}(x,0)=\sum_{j=1}^{J}c\frac{f(b^{-j}r)b^{j}sin(2\alpha)}{j },\quad \frac12>b>0,$$ where the radial function $0<f\in C^{\infty}$ has $supp(f)\in [\frac12, \frac32]$, $c>0$ and $J\in \mathds{N}$.
 The main difficulty when considering this type of initial conditions is that the usual energy estimates only give existence for a short time interval which does not provide enough growth in $H^2$. To obtain improved time intervals of existence we decompose our solution as a sum of pseudo-solutions with initial conditions
$$c\frac{f(b^{-j}r)b^{j}sin(2\alpha)}{j }$$ for $j=1, ..., J$. To finish the proof we perturb this solution with a small $H^2$ function localized around the origin that will experience  very large norm growth. 

The paper is organized as follows. First in section 2 we prove strong ill-posedness and non existence for the space $C^{k}$. In section 3 we show strong ill-posedness and non existence for Sobolev spaces in the supercritical case. Finally in section 4 we prove strong ill-posedness for the critical $H^{2}$ space.
\subsection{Notation}

In this paper we will consider functions $f(x):\mathds{R}^{2}\rightarrow \mathds{R}$ in $C^{k}$ with $k$ a positive integer and $H^{s}$ with $s$ a positive real number. These spaces allow many different equivalent norms, but we will specifically use

$$||f(x)||_{C^{k}}=\sum_{i=0}^{k}\sum_{j=0}^{i}||\frac{\partial^{i}f(x)}{\partial^{j}x_{1}\partial^{i-j} x_{2}}||_{L^{\infty}}$$
and for $H^{s}$, when $s$ is a positive integer we will use

$$||f(x)||_{H^{s}}=\sum_{i=0}^{s}\sum_{j=0}^{i}||\frac{\partial^{i}f(x)}{\partial^{j}x_{1}\partial^{i-j} x_{2}}||_{L^{2}},$$
where the derivative is understood in the weak sense.




For $s$ non integer, the standard way of defining the norm is by 

$$||f(x)||_{H^{s}}=||\mathcal{F}^{-1}\Big[(1+|\xi|^2)^{\frac{s}{2}}\mathcal{F}f\Big]||_{L^{2}},$$
where $\mathcal{F}$ is the Fourier transform. We will not require to use this definition to compute the norm in these spaces through this paper. For $s$ a positive integer, we will sometimes write

$$||f(x)1_{A}||_{H^{s}},$$
where $1_{A}$ is the characteristic function in the set $A$.
This is slightly an abuse of notation since the function $f(x)1_{A}$ may not be in $H^{s}$, but we will use this as a more compact notation to write

$$\sum_{i=0}^{s}\sum_{j=0}^{i} (\int_{A} (\frac{\partial^{i}f(x)}{\partial^{j}x_{1}\partial^{i-j} x_{2}})^2 dx)^{\frac12}.$$

Analogously, we will use

$$||f(x) 1_{A}||_{C^{k}}:=\sum_{i=0}^{k}\sum_{j=0}^{i}\text{ess-sup}_{x\in A}(\frac{\partial^{i}f(x)}{\partial^{j}x_{1}\partial^{i-j} x_{2}}).$$
We will work both in normal cartesian coordinates and in polar coordinates, using the change of variables $x_{1}=rcos(\alpha)$, $x_{2}=rsin(\alpha)$. We will sometimes define a function in the space $(x_{1},x_{2})$ $f(x)$ and then refer to $f(r,\alpha)$ (or vice versa), and this is an abuse of notation since we should actually write, if $F(r,\alpha)$ is the change of variables that takes us from $(r,\alpha)$ to $(x_{1},x_{2})$, $f(F(r,\alpha))$. Furthermore, given a function in polar coordinates, we define

$$||f(r,\alpha)||_{H^{s}}:=||f(F(r,\alpha))||_{H^{s}},$$

$$||f(r,\alpha)||_{C^{k}}:=||f(F(r,\alpha))||_{C^{k}}.$$

For two sets $A_{1}$, $A_{2}$, we will use $d(A_{1},A_{2})$ to refer to the distance between the sets.

\subsection{Ill-posedness}\label{illposedness}

Since we will be dealing with ill-posedness through this paper, it is important that we clarify exactly what we mean by mild and strong ill-posedness, specially since one can give similar (but not  necessarily equivalent) definitions of these concepts. Through this paper we will use the same definition as in \cite{Elgindi}, that is
\begin{definition}
Given spaces $X$, $Y$ with $Y$ continuously embedded in $X$ and an evolution equation

$$\frac{\partial f(x,t)}{\partial t}=G(f(x,t))$$

$$f(x,0)=f_{0}(x)$$
we say that the evolution equation is mildly ill-posed if we can find $f_{\epsilon}(x)\in X$ such that there exists a unique solution $f_{\epsilon}(x,t)$ in $L^{\infty}([0,\epsilon]; Y)$ to our evolution equation with initial conditions $f_{\epsilon}(x)$ such that $||f_{\epsilon}(x)||_{X}\leq \epsilon$ but there exists a time $t\in(0,\epsilon]$ such that $||f_{\epsilon}(x)||_{X}\geq c$ with $c>0$ some constant independent of $\epsilon$. Furthermore, if we can take $c=\frac{1}{\epsilon}$, then we say that the problem is strongly ill-posed.
\end{definition}

What these notions tell us about the evolution equation is that it is not well behaved in the space $X$. More precisely, mild ill-posedness tells us that the solution map is not continuous with respect to the initial conditions, and strong ill-posedness shows both that and arbitrarily fast norm growth, which could potentially lead to an instantaneous blow up, and therefore, to non-existence of solutions.

Although strong and mild ill-posedness are related, and in fact in some situations they are equivalent (for example, if your evolution equation has appropriate scaling properties), one does not imply the other. In fact, if we consider a radial function $f(r)=r^{k+\gamma}g(r)$ with $g(r)$ a $C^{\infty}$ function such that $g(r)=1$ if $r\in[0,1]$, $g(r)=0$ if $r\geq2$, we have that the evolution equation for perturbations of $f(r)$ for SQG, which is

$$\frac{\partial \theta_{pert}}{\partial t}+u(\theta_{pert})\cdot\nabla \theta_{pert}+u(f(r))\cdot\nabla \theta_{pert}+u(\theta_{pert})\cdot\nabla f(r)=0$$
is mildly ill-posed in $C^{k,\gamma}$ but not strongly ill-posed in $C^{k,\gamma}$. The mild ill-posedness is easy to obtain by noting that, if $v(\theta_{pert(x,0)})(x=0)=(a,b)\neq(0,0)$, then, by using that $\theta_{pert}(x,t)=w(x,t)-f(r)$, with $w(x,t)$ the solution to SQG with initial conditions $\theta_{pert}(x,0)+f(r)$, we get that

$$\lim_{t\rightarrow 0^{+}}\Big(lim_{h\rightarrow 0^{+}}\frac{\frac{\partial^{k}\theta_{pert(x,t)}}{\partial x_{1}^{k}}(x=0)-\frac{\partial^{k}\theta_{pert(x,t)}}{\partial x_{1}^{k}}(x=h)}{h^{\gamma}}\Big)=\prod_{i=1}^{k}(i+\gamma), $$
and thus 
$$\lim_{t\rightarrow 0^{+}}||\theta_{pert}(x,t)||_{C^{k\gamma}}\geq \prod_{i=1}^{k}(i+\gamma)$$
and since this can be obtained independently of the norm of $\theta_{pert}(x,0)$, we obtain mild ill-posedness.
But since we know that solutions of SQG in $C^{k,\gamma}$ fulfil

$$\frac{\partial ||\theta||_{C^{k,\gamma}}}{\partial t}\leq C ||\theta||_{C^{k,\gamma}}^2$$
then
\begin{align*}
    &\frac{\partial ||\theta_{pert}||_{C^{k,\gamma}}}{\partial t}\leq \frac{\partial ||f(r)+\theta_{pert}||_{C^{k,\gamma}}}{\partial t}+\frac{\partial ||f(r)||_{C^{k,\gamma}}}{\partial t}\\
    &\leq C ||f(r)+\theta_{pert}||_{C^{k,\gamma}}^2\leq C (||f(r)||_{C^{k,\gamma}}^2+||\theta_{pert}||_{C^{k,\gamma}}^2)
\end{align*}
which implies that strong ill-posedness is not possible.

\section{Strong ill-posedness and non existence in $C^{k}$}

To prove ill-posedness in $C^k$ we construct fast growth solutions by perturbing in a suitable way a stationary smooth radial solution. In contrast, there are previous results  (\cite{Globalsmooth} and \cite{ccz}) where the perturbation of a radial function led to global $C^{4}$ rotating solutions and enhanced lifespan of solutions respectively.

In this section we will show that, for a specific kind of perturbation we can predict the behaviour of the solution with a very small error. The perturbation will be composed of functions of the form

$$f(N^{1/2}(r-1)+1)sin(Nn\alpha) $$
with $f$ a given smooth function and $N,n$ integers. Below we will obtain the properties that will alow us to work with this kind of functions.

\subsection{Estimates on the velocity field.}

In this section we will use the following expression of the velocity field 
$$v(\theta(.))(x)=\frac{\Gamma (3/2)}{\pi^{3/2}}P.V. \int_{\mathds{R}^2} \frac{(x-y)^{\perp}\theta(y)}{|x-y|^{3}}dy_{1}dy_{2}$$
with $v=(v_{1},v_{2})$ and for a vector $(a,b)$ we define $(a,b)^{\perp}:=(-b,a)$.

We will omit the constant on the outside of the integral from now on, since all the results we will obtain would remain the same if we were to change $\frac{\Gamma (3/2)}{\pi^{3/2}}$ for an arbitrary (non-zero) constant.

\begin{lemma}\label{velocidaderror}
Given natural numbers $n$, $N$ and a $L^{\infty}$ function $g_{N}(r):[0,\infty)\rightarrow \mathds{R}$ with support in $(1-\frac{N^{-1/2}}{2},1+\frac{N^{-1/2}}{2})$  we have that, for $\theta(r,\alpha)=g_{N}(r)sin(Nn\alpha)$, there exists a constant $C$ (depending on $n$) such that, for $N$ big enough and $r\in [1-N^{-\frac12},1+N^{-\frac12}]$

$$|v_{r}(\theta(.,.))(r,\alpha)-cos(Nn\alpha)\int_{\mathds{R}\times[-\pi,\pi]} \frac{r^2 \alpha' g_{N}(r+h)sin(Nn\alpha')}{|h^2+r^2(\alpha')^2|^{3/2}}d\alpha'dh|$$
$$\leq C  ||g_{N}||_{L^{\infty}}  N^{-1/2}.$$

Analogously, for $\theta(r,\alpha)=g_{N}(r)cos(Nn\alpha)$ we have that 

$$|v_{r}(\theta(.,.))(r,\alpha)+sin(Nn\alpha)\int_{\mathds{R}\times[-\pi,\pi]} \frac{r^2 \alpha' g_{N}(r+h)sin(Nn\alpha')}{|h^2+r^2(\alpha')^2|^{3/2}}d\alpha'dh|$$
$$\leq C||g_{N}(r)||_{L^{\infty}}N^{-1/2}.$$
\end{lemma}
Before we get into the proof, a couple of comments need to be made. First, $v_{r}$, refers to the radial component of the velocity at a given point, that is to say, if we call $\hat{x}$ to the unitary vector in the direction of $x$ then
$$v_{r}(\theta(.))(x)=P.V. \int_{\mathds{R}^2} \hat{x} \frac{(x-y)^{\perp}\theta(y)}{|x-y|^{3}}dy_{1}dy_{2}.$$
However, the expression obtained in lemma \ref{velocidaderror} requires us to work in polar coordinates. Therefore, considering a generic function $f(r)sin(k\alpha)$ and making the usual changes of variables $(x_{1},x_{2})=r(cos(\alpha),sin(\alpha))$, $(y_{1},y_{2})=r'(cos(\alpha'),sin(\alpha'))$ we obtain
\begin{align*}
    &v_{r}(\theta(\ .\ ,\ .\ ))(r,\alpha)
\end{align*}
\begin{equation*}
       =P.V.\int_{\mathds{R}\times[-\pi,\pi]}(r')^2\frac{(cos(\alpha)sin(\alpha')-sin(\alpha)cos(\alpha'))f(r')sin(k\alpha')}{|(rcos(\alpha)-r'cos(\alpha'))^2+(rsin(\alpha)-r'sin(\alpha'))^2|^{3/2}}d\alpha'dr'
\end{equation*}
 \begin{align}\label{velocidadradial}
   &=P.V.\int_{\mathds{R}\times[-\pi,\pi]}(r')^2\frac{sin(\alpha'-\alpha)}{|(r-r')^2+2rr'(1-cos(\alpha-\alpha'))|^{3/2}}f(r')sin(k\alpha')d\alpha'dr'\nonumber\\
   &=cos(k\alpha)P.V.\int_{\mathds{R}\times[-\pi,\pi]}(r')^2\frac{sin(\alpha'-\alpha)f(r')sin(k\alpha'-k\alpha)}{|(r-r')^2+2rr'(1-cos(\alpha-\alpha'))|^{3/2}}d\alpha'dr'\nonumber\\
   &=cos(k\alpha)P.V.\int_{\mathds{R}\times[-\pi,\pi]}(r+h)^2\frac{sin(\alpha')f(r+h)sin(k\alpha')}{|h^2+2(r+h)r(1-cos(\alpha'))|^{3/2}}d\alpha'dh,
\end{align}
where we have used trigonometric identities and eliminated the terms that are odd with respect to $\alpha'-\alpha$. Note that in the last line we have relabeled $\alpha'-\alpha$ as $\alpha'$ for a more compact notation. 
Analogously if $\theta(r,\alpha)=f(r)cos(k\alpha)$ we obtain
\begin{align*}
    &v_{r}(\theta(\ .\ ,\ .\ ))(r,\alpha)=-sin(k\alpha)\int_{\mathds{R}\times[-\pi,\pi]}(r+h)^2\frac{sin(\alpha')f(r+h)sin(k\alpha')}{|h^2+2(r+h)r(1-cos(\alpha'))|^{3/2}}d\alpha'dh.\\
\end{align*}

With this, we are now ready to start the proof of lemma \ref{velocidaderror}.

\begin{proof}

We need to find bounds for

$$\int_{\mathds{R}\times[-\pi,\pi]} \frac{r^2 \alpha' g_{N}(r+h)sin(Nn\alpha')}{|h^2+r^2(\alpha')^2|^{3/2}}d\alpha'dh$$
$$-\int_{\mathds{R}\times[-\pi,\pi]}(r+h)^2\frac{sin(\alpha')g_{N}(r+h)sin(Nn\alpha')}{|h^2+2(r+h)r(1-cos(\alpha'))|^{3/2}}d\alpha'dh$$
with $g_{N}(r)$ satisfying our hypothesis. We will first focus on 

\begin{equation}\label{integral A}
    I_{A}:=|\int_{A} \frac{r^2 \alpha' g_{N}(r+h)sin(Nn\alpha')}{|h^2+r^2(\alpha')^2|^{3/2}}d\alpha'dh
\end{equation}
$$-\int_{A}(r+h)^2\frac{sin(\alpha')g_{N}(r+h)sin(Nn\alpha')}{|h^2+2(r+h)r(1-cos(\alpha'))|^{3/2}}d\alpha'dh|$$
with $A:=[-2N^{-1/2},2N^{-1/2}]\times[-2N^{-1/2},2N^{-1/2}]$. This is accomplished in several steps. It should be noted that the constant $C$  may depend on $n$ and it may change through the proof, as it is the name we use for a generic constant that is independent of $N$ and $g$.

Step 1:
\begin{align*}
    &|\int_{A}(r+h)^2\frac{(sin(\alpha')-\alpha')g_{N}(r+h)sin(Nn\alpha')}{|h^2+2(r+h)r(1-cos(\alpha'))|^{3/2}}d\alpha'dh|\\
    &\leq C\int_{A}(r+h)^2\frac{|\alpha'|^3 |g_{N}(r+h)|}{|h^2+2(r+h)r(1-cos(\alpha'))|^{3/2}}d\alpha'dh\\
    &\leq C\int_{A}|g_{N}(r+h)|d\alpha'dh\\
    &\leq C N^{-1}||g_{N}||_{L^{\infty}}\\
\end{align*}

Step 2:
Defining 
$$F(r,h,\alpha'):=\frac{1}{|h^2+2(r+h)r(1-cos(\alpha'))|^{3/2}}-\frac{1}{|h^2+(r+h)r(\alpha')^2)|^{3/2}}$$
we estimate the following integral by
\begin{align*}
    &|\int_{A}(r+h)^2\alpha'g_{N}(r+h)sin(Nn\alpha')F(r,h,\alpha')d\alpha'dh|\\
    &\leq C\int_{A}|\alpha'||g_{N}(r+h)|\frac{(\alpha')^4}{|h^2+2(r+h)r(1-cos(\alpha'))|^{5/2}}d\alpha'dh\\
    &\leq C\int_{A}|g_{N}(r+h)|d\alpha'dh\\
    &\leq C N^{-1}||g_{N}||_{L^{\infty}}.\\
\end{align*}
Step 3:
\begin{align*}
    &|\int_{A}((r+h)^2-r^2)\frac{\alpha'g_{N}(r+h)sin(Nn\alpha')}{|h^2+(r+h)r(\alpha')^2)|^{3/2}}d\alpha'dh|\\
    &\leq C\int_{A}|h|\frac{|\alpha'||g_{N}(r+h)|}{|h^2+(r+h)r(\alpha')^2)|^{3/2}}d\alpha'dh\\
    &\leq C\int_{A}\frac{|g_{N}(r+h)|}{|h^2+(r+h)r(\alpha')^2|^{1/2}}d\alpha'dh\\
    &\leq C N^{-1/2}||g_{N}||_{L^{\infty}}\\
\end{align*}

Combining all these three steps we conclude 

$$|\int_{A} \frac{r^2 \alpha' g_{N}(r+h)sin(Nn\alpha')}{|h^2+r(r+h)(\alpha')^2|^{3/2}}d\alpha'dh-\int_{A}\frac{(r+h)^2sin(\alpha')g_{N}(r+h)sin(Nn\alpha')}{|h^2+2(r+h)r(1-cos(\alpha'))|^{3/2}}d\alpha'dh|$$
$$\leq C||g_{N}||_{L^{\infty}}N^{-1/2},$$

and to bound the contribution of the integral in $A$ we also need

\begin{align*}
    &|\int_{A}\frac{r^2\alpha'g_{N}(r+h)sin(Nn\alpha')}{|h^2+(r+h)r(\alpha')^2)|^{3/2}}-\frac{r^2\alpha'g_{N}(r+h)sin(Nn\alpha')}{|h^2+r^2(\alpha')^2)|^{3/2}}d\alpha'dh|\\
    &\leq C\int_{A}|\alpha'||g_{N}(r+h)|\frac{(\alpha')^2|h|}{|h^2+\frac{r^2}{2}(\alpha')^2)|^{5/2}}d\alpha'dh\\
    &\leq C\int_{A}\frac{|g_{N}(r+h)|}{|h^2+\frac{r^2}{2}(\alpha')^2|^{1/2}}d\alpha'dh\\
    &\leq C N^{-1/2}||g_{N}||_{L^{\infty}}.\\
\end{align*}
Therefore adding and subtracting 
$$\int_{A}\frac{r^2\alpha'g_{N}(r+h)sin(Nn\alpha')}{|h^2+(r+h)r(\alpha')^2)|^{3/2}}$$
to (\ref{integral A}) we obtain that

$$I_{A}\leq C||g_{N}||_{L^{\infty}}N^{-1/2}.$$
Finally, we need to deal with the integral outside of $A$. First we bound the following integral

\begin{align*}
    &\int_{\mathds{R}\times [-\pi,\pi]\setminus A} \frac{r^2 \alpha' g_{N}(r+h)sin(Nn\alpha')}{|h^2+r^2(\alpha')^2|^{3/2}}d\alpha'dh\\
    & =2\int_{[-2N^{1/2},2N^{-1/2}]}\int_{[2N^{-1/2},\pi]}  \frac{r^2 \alpha' g_{N}(r+h)sin(Nn\alpha')}{|h^2+r^2(\alpha')^2|^{3/2}}d\alpha'dh.\\
\end{align*}
To do this we compute, fixed arbitrary $h$ and $r$, the integral over an interval of the form $\alpha\in[k\frac{2\pi}{Nn}-\frac{\pi}{2Nn},(k+1)\frac{2\pi}{Nn}-\frac{\pi}{2Nn}]$ (which we will denote $[\alpha_{k}, \alpha_{k+1}]$). Note that it has the length of the period of $sin(Nn\alpha)$ and that $sin(Nn\alpha)$ is an even function around the point $k\frac{2\pi}{Nn}+\frac{\pi}{2Nn}$. 




If we define 
$$H(\alpha',h,r):=\frac{ \alpha' }{|h^2+r^2(\alpha')^2|^{3/2}}$$
we have that

\begin{align*}
    &\int_{[\alpha_{k},\alpha_{k+1}]} sin(Nn\alpha') \frac{ \alpha' }{|h^2+r^2(\alpha')^2|^{3/2}}d\alpha'\\
    &=\int_{[\alpha_{k},\alpha_{k+1}]} sin(Nn\alpha') \Big(H(\frac{\alpha_{k}+\alpha_{k+1}}{2},h,r)\\
    &+\frac{\partial H(\frac{\alpha_{k}+\alpha_{k+1}}{2},h,r)}{\partial \alpha'}(\alpha'-\frac{\alpha_{k}+\alpha_{k+1}}{2})+\frac{\partial^2 H(c(\alpha'),h,r)}{\partial \alpha'^2}\frac12(\alpha'-\frac{\alpha_{k}+\alpha_{k+1}}{2})^2d\alpha'\Big)\\
    &=\int_{[\alpha_{k},\alpha_{k+1}]} sin(Nn\alpha')\frac{\partial^2 H(c(\alpha'),h,r)}{\partial \alpha'^2}\frac12(\alpha'-\frac{\alpha_{k}+\alpha_{k+1}}{2})^2d\alpha'\\
    &\leq C\int_{[\alpha_{k},\alpha_{k+1}]}| sin(Nn\alpha')|(\alpha'-\frac{\alpha_{k}+\alpha_{k+1}}{2})^2 \frac{ 1 }{|h^2+r^2\alpha_{k}^2|^{2}}d\alpha'\\
    &\leq C\Big(\frac{2\pi}{Nn}\Big)^{3}\frac{1}{|h^2+r^2\alpha_{k}^2|^{2}},\\
\end{align*}
where we have used a second degree taylor expansion around $\frac{\alpha_{k}+\alpha_{k+1}}{2}$ for $H$, and $c(\alpha')$ is where we need to evaluate the second derivative to actually obtain an equality.
Now, adding over all the intervals $[\alpha_{k},\alpha_{k+1}]$ with $\pi-\frac{2\pi}{Nn}\geq \alpha_{k}\geq 2N^{-1/2}$, we get the upper bound

\begin{align*}
    &\sum^{\pi-\frac{2\pi}{Nn}}_{\alpha_{k}\geq 2N^{-1/2}} C\Big(\frac{2\pi}{Nn}\Big)^{3}\frac{1}{|h^2+r^2\alpha_{k}^2|^{2}}\leq \sum^{\infty}_{k\geq \frac{N^{1/2}n}{\pi}} C\Big(\frac{2\pi}{Nn}\Big)^{3}\frac{1}{|h^2+r^2\alpha_{k}^2|^{2}}\\
    &\leq C\Big(\frac{2\pi}{Nn}\Big)^{3} \int^{\infty}_{\frac{N^{1/2}n}{\pi}-1} \frac{1}{|h^2+r^2( x \frac{2\pi}{Nn}-\frac{\pi}{2Nn})^2|^{2}}dx \\
    &\leq C\Big(\frac{2\pi}{Nn}\Big)^{3} \int^{\infty}_{\frac{N^{1/2}n}{\pi}-2} \frac{1}{|h^2+(r x \frac{2\pi}{Nn})^2|^{2}}dx \\
    &\leq C\Big(\frac{2\pi}{Nn}\Big)^{3} \Big(\frac{2\pi}{Nn}\Big)^{-4}(\frac{2\pi}{N^{1/2}n})^{3}\leq CN^{-1/2}, \\
\end{align*}
where we took $N$ big to pass from the third to the fourth line.
The only contribution missing now from the integral in the $\alpha'$ variable, if we call $\alpha_{k_{0}}$ the smallest $\alpha_{k}$ such that $\alpha_{k}\geq 2N^{-1/2}$ and $\alpha_{\infty}$ the biggest one with $\pi\geq \alpha_{\infty}$, is

$$\int_{[2N^{-1/2},\alpha_{k_{0}}]\cup [\alpha_{\infty},\pi]} sin(Nn\alpha') \frac{ \alpha' }{|h^2+r^2(\alpha')^2|^{3/2}}d\alpha',$$
but

$$|\int^{\alpha_{k_{0}}}_{2N^{-1/2}} sin(Nn\alpha') \frac{ \alpha' }{|h^2+r^2(\alpha')^2|^{3/2}}d\alpha'|\leq C,$$
$$|\int_{\alpha_{k_{\infty}}}^{\pi} sin(Nn\alpha') \frac{ \alpha' }{|h^2+r^2(\alpha')^2|^{3/2}}d\alpha'|\leq \frac{C}{N}.$$
Combining all three contributions and integrating with respect to $h$ we get

\begin{align*}
    & |2\int_{[-2N^{1/2},2N^{-1/2}]}\int_{[2N^{-1/2},\pi]}  \frac{r^2 \alpha' g_{N}(r+h)sin(Nn\alpha')}{|h^2+r^2(\alpha')^2|^{3/2}}d\alpha'dh|\\
    &\leq \int_{[-2N^{1/2},2N^{-1/2}]}C|g_{N}(r+h)| dh\leq C||g_{N}||_{L^{\infty}}N^{-1/2}.\\
\end{align*}

The term

$$|\int_{\mathds{R}\times[-\pi,\pi]\setminus A}(r+h)^2\frac{sin(\alpha')g_{N}(r+h)sin(Nn\alpha')}{|h^2+2(r+h)r(1-cos(\alpha'))|^{3/2}}d\alpha'dh|$$
is bounded in a similar fashion, integrating first with respect to $\alpha'$ in intervals of the form $[\alpha_{k},\alpha_{k+1}]$ and then bounding by brute force the parts that are not covered exactly by said intervals, and with that we would be done.

\end{proof}

Now that we have a manageable expression for the radial velocity  we are ready to compute it explicitly (with some error) for some special kind of functions.

\begin{lemma}\label{velocidadcartesianas}
Given natural numbers $n,$ $N$ and a $C^2$ function $g_{N}(.):\mathds{R}\rightarrow\mathds{R}$ with support in the interval $(1-\frac{N^{-1/2}}{2},1+\frac{N^{1/2}}{2})$ satisfying $||g_{N}||_{C^{i}}\leq MN^{i/2}$ for $i=0,1,2$ there exists a constant $C_{0}\neq 0$ (independent of $N$, $n$ and $g_{N}$) such that for $\tilde{x}\in [1-N^{-\frac12},1+N^{-\frac12}]$

\begin{equation*}
    |C_{0}g_{N}(\tilde{x})-\int_{\mathds{R}\times[-\pi,\pi]} g_{N}(\tilde{x}+h_{1})\frac{sin(Nnh_{2})h_{2}}{(h_{1}^2+h_{2}^2)^{3/2}}dh_{1}dh_{2}|
\end{equation*}
\begin{equation}\label{ecuacionlemmatres}
        \leq C M N^{-1/2},
\end{equation}
with $C$ depending on $n$.
\end{lemma}

\begin{proof}

The strategy of this proof is to first show that

\begin{align*}
    &|\int_{\mathds{R}\times[-\pi,\pi]} (g_{N}(\tilde{x}+h_{1})-g_{N}(\tilde{x}))\frac{sin(Nnh_{2})h_{2}}{(h_{1}^2+h_{2}^2)^{3/2}}dh_{1}dh_{2}|
\end{align*}
\begin{equation}\label{gconstante}
        \leq C M N^{-1/2},
\end{equation}
and then prove that 
\begin{equation}\label{INn}
    I_{N,n}:=\int_{\mathds{R}\times[-\pi,\pi]} sin(Nnh_{2})\frac{h_{2}}{(h_{1}^2+h_{2}^2)^{3/2}}dh_{1}dh_{2}
\end{equation}
is a Cauchy series with respect to $N$, satisfying
\begin{equation}\label{cauchy}
    |I_{N_{1},n}-I_{N_{2},n}|\leq C sup(N_{1},N_{2})^{-1/2}
\end{equation}
with $C$ depending on $n$.

Combining both of these results and taking 

$$C_{0}=\text{lim}_{N\rightarrow\infty } I_{N,n}$$
we  obtain  (\ref{ecuacionlemmatres}), and we only need to check that $C_{0}$ is different from zero and independent of $n$.

We start by obtaining bound (\ref{gconstante}), noting that, by parity
\begin{align*}
    &|\int_{[-2N^{-1/2},2N^{-1/2}]\times[-\pi,\pi]} (g_{N}(\tilde{x}+h_{1})-g_{N}(\tilde{x}))sin(Nnh_{2})\frac{h_{2}}{(h_{1}^2+h_{2}^2)^{3/2}}dh_{1}dh_{2}|\\
    &=|\int_{[0,2N^{-1/2}]\times[-\pi,\pi]} \frac{(g_{N}(\tilde{x}+h_{1})+g_{N}(\tilde{x}-h_{1})-2g_{N}(\tilde{x}))sin(Nnh_{2})h_{2}}{(h_{1}^2+h_{2}^2)^{3/2}}dh_{1}dh_{2}|.\\
\end{align*}

We start by fixing  some $h_{1}$ and obtaining bounds for the integral with respect to $h_{2}$. This is done as in lemma \ref{velocidaderror}, dividing in periods of length $\frac{2\pi}{Nn}$ starting at $\frac{\pi}{2Nn}$, and approximating $\frac{h_{2}}{h_{1}^2+h_{2}^2}$ by its second order Taylor expansion, since the first two orders will cancel. That way, for the interval with $h_{2}\in[k\frac{2\pi}{Nn}+\frac{\pi}{2Nn},(k+1)\frac{2\pi}{Nn}+\frac{\pi}{2Nn}]$ we obtain the bound
\begin{align}
    &|\int_{k\frac{2\pi}{Nn}+\frac{\pi}{2Nn}}^{(k+1)\frac{2\pi}{Nn}+\frac{\pi}{2Nn}}sin(Nnh_{2})\frac{h_{2}}{(h_{1}^2+h_{2}^2)^{\frac32}}dh_{2}|\leq C\Big(\frac{2\pi}{Nn}\Big)^{3}\frac{1}{(h_{1}^2+(\frac{k2\pi}{Nn})^2)^{2}}.\label{periodo}
\end{align}
We can add periods contained in the interval $[0,2N^{-1/2}]$ and, if we denote by $k_{\infty}=k_{\infty}(N,n)$ the biggest integer $k$ such that $(k+1)\frac{2\pi}{Nn}+\frac{\pi}{2Nn}\leq 2N^{-1/2}$, we get that
\begin{align*}
    &|\int_{\frac{5\pi}{2Nn}}^{(k_{\infty}+1)\frac{2\pi}{Nn}+\frac{\pi}{2Nn}}sin(Nnh_{2})\frac{h_{2}}{(h_{1}^2+h_{2}^2)^{\frac32}}dh_{2}|\\
    &\leq \sum_{k=1}^{k_{\infty}}C\Big(\frac{2\pi}{Nn}\Big)^{3}\frac{1}{(h_{1}^2+(\frac{k2\pi}{Nn})^2)^{2}}\leq C\Big(\frac{2\pi}{Nn}\Big)^{3} \int_{0}^{k_{\infty}}\frac{1}{(h_{1}^2+(\frac{x2\pi}{Nn})^2)^{2}}dx\\
    &\leq C\Big(\frac{2\pi}{Nn}\Big)^{3} \int_{0}^{k_{\infty}}\frac{1}{(h_{1}+\frac{x2\pi}{Nn})^{4}}dx= C\Big(\frac{2\pi}{Nn}\Big)^{3} \int_{\frac{h_{1}Nn}{2\pi}}^{k_{\infty}+\frac{h_{1}Nn}{2\pi}} \frac{1}{(\frac{x2\pi}{Nn})^{4}}dx\\
    &\leq C \Big(\frac{2\pi}{Nn}\Big)^{2} \frac{1}{h_{1}^{3}}.\\
\end{align*}
This allows us to bound the contribution when $h_{1}\geq \frac{2\pi}{Nn}$ by dividing it in three parts:

1) If $h_{2}\leq \frac{5\pi}{2Nn}$:
\begin{align*}
    &|\int_{\frac{2\pi}{Nn}}^{2N^{-1/2}}\int_{0}^{\frac{5\pi}{2Nn}}(g_{N}(\tilde{x}+h_{1})+g_{N}(\tilde{x}-h_{1})-2g_{N}(\tilde{x}))sin(Nnh_{2})\frac{h_{2}}{(h_{1}^2+h_{2}^2)^{3/2}}dh_{2}dh_{1}|\\
    &\leq |C\int_{\frac{2\pi}{Nn}}^{2N^{-1/2}}\frac{5\pi}{2Nn}M h_{1}^2N\frac{1}{h_{1}^2}dh_{1}|\leq C M N^{-1/2}.\\
\end{align*}

2) If $\frac{5\pi}{2Nn}\leq h_{2}\leq (k_{\infty}+1)\frac{2\pi}{Nn}+\frac{\pi}{2Nn}$:
\begin{align*}
    &|\int_{\frac{2\pi}{Nn}}^{2N^{-1/2}}\int_{\frac{5\pi}{2Nn}}^{(k_{\infty}+1)\frac{2\pi}{Nn}+\frac{\pi}{2Nn}}(g_{N}(\tilde{x}+h_{1})+g_{N}(\tilde{x}-h_{1})-2g_{N}(\tilde{x}))\frac{sin(Nnh_{2})h_{2}}{(h_{1}^2+h_{2}^2)^{3/2}}dh_{2}dh_{1}|\\
    &\leq \int_{\frac{2\pi}{Nn}}^{2N^{-1/2}}|g_{N}(\tilde{x}+h_{1})+g_{N}(\tilde{x}-h_{1})-2g_{N}(\tilde{x})|\ |\int_{\frac{5\pi}{2Nn}}^{(k_{\infty}+1)\frac{2\pi}{Nn}+\frac{\pi}{2Nn}}\frac{sin(Nnh_{2})h_{2}}{(h_{1}^2+h_{2}^2)^{3/2}}dh_{2}| dh_{1}\\
    &\leq C\int_{\frac{2\pi}{Nn}}^{2N^{-1/2}}Mh_{1}^2N\Big(\frac{2\pi}{Nn}\Big)^{2} \frac{1}{h_{1}^{3}}dh_{1}\leq C MN^{-1}log(N).\\
\end{align*}

3) If $(k_{\infty}+1)\frac{2\pi}{Nn}+\frac{\pi}{2Nn}\leq h_{2}\leq 2N^{-1/2}$:

\begin{align*}
    &|\int_{\frac{2\pi}{Nn}}^{2N^{-1/2}}\int_{(k_{\infty}+1)\frac{2\pi}{Nn}+\frac{\pi}{2Nn}}^{2N^{-1/2}}(g_{N}(\tilde{x}+h_{1})+g_{N}(\tilde{x}-h_{1})-2g_{N}(\tilde{x}))\frac{sin(Nnh_{2})h_{2}}{(h_{1}^2+h_{2}^2)^{3/2}}dh_{2}dh_{1}|\\
    &\leq C\int_{\frac{2\pi}{Nn}}^{2N^{-1/2}}Mdh_{1}\leq C MN^{-\frac{1}{2}}.\\
\end{align*}
Finally, we bound the error when $h_{1}\leq \frac{2\pi}{Nn}$:

1) If $|h_{2}|\leq 2N^{-1/2}$
\begin{align*}
    &|\int_{0}^{\frac{2\pi}{Nn}}\int_{0}^{2N^{-1/2}}(g_{N}(\tilde{x}+h_{1})+g_{N}(\tilde{x}-h_{1})-2g_{N}(\tilde{x}))sin(Nnh_{2})\frac{h_{2}}{(h_{1}^2+h_{2}^2)^{3/2}}dh_{2}dh_{1}|\\
    &\leq \int_{0}^{\frac{2\pi}{Nn}}\int_{0}^{2N^{-1/2}}Mh_{1}^2N\frac{1}{(h_{1}^2+h_{2}^2)}dh_{2}dh_{1}\leq C MN^{-1/2}.\\
\end{align*}

2) If $|h_{2}|\geq 2N^{-1/2}$
\begin{align*}
    &|\int_{0}^{\frac{2\pi}{Nn}}\int_{2N^{-1/2}}^{\pi}(g_{N}(\tilde{x}+h_{1})+g_{N}(\tilde{x}-h_{1})-2g_{N}(\tilde{x}))\frac{sin(Nnh_{2})h_{2}}{(h_{1}^2+h_{2}^2)^{3/2}}dh_{2}dh_{1}|\\
    &\leq \int_{0}^{\frac{2\pi}{Nn}}\int_{2N^{-1/2}}^{\pi}\frac{M}{2}h_{1}^2N\frac{1}{(h_{1}^2+h_{2}^2)}dh_{2}dh_{1}\leq C MN^{-1}.\\
\end{align*}

Combining all these bounds we obtain (\ref{gconstante}).
Therefore, we have that it is enough to prove that

\begin{align*}
    &|C_{0}g_{N}(\tilde{x})-\int_{\mathds{R}\times[-\pi,\pi]} g_{N}(\tilde{x})sin(Nnh_{2})\frac{h_{2}}{(h_{1}^2+h_{2}^2)^{3/2}}dh_{1}dh_{2}|\\
    &\leq C MN^{-1/2},\\
\end{align*}
which is equivalent to study the behaviour of $I_{N,n}$, defined as in (\ref{INn}).

To obtain the properties of $I_{N,n}$, we start by transforming the integral with a change of variables $\bar{h}_{1}=Nnh_{1}$, $\bar{h}_{2}=Nnh_{2}$, although we will relabel $\bar{h}_{1},\bar{h}_{2}$ as $h_{1}$, $h_{2}$ to simplify the notation.

\begin{align*}
    &g_{N}(\tilde{x})\int_{[-2N^{-1/2},2N^{-1/2}]\times[-\pi,\pi]} sin(Nnh_{2})\frac{h_{2}}{(h_{1}^2+h_{2}^2)^{3/2}}dh_{1}dh_{2}\\
    &=g_{N}(\tilde{x})\int_{[-2nN^{1/2},2nN^{1/2}]\times[-Nn\pi,Nn\pi]} sin(h_{2})\frac{h_{2}}{(h_{1}^2+h_{2}^2)^{3/2}}dh_{1}dh_{2}.\\
\end{align*}

If we compare the integral for different values of $N$, $N_{1}\geq N_{2}$ we get

\begin{align*}
    &I_{N_{1},n}-I_{N_{2},n}=\int_{A\cup B} sin(h_{2})\frac{h_{2}}{(h_{1}^2+h_{2}^2)^{3/2}}dh_{1}dh_{2}\\
\end{align*}
with 
$$A=[-2nN_{2}^{1/2},2nN_{2}^{1/2}]\times[N_{2}n\pi,N_{1}n\pi]\cup[-2nN_{2}^{1/2},2nN_{2}^{1/2}]\times[-N_{1}n\pi,-N_{2}n\pi],$$
$$B=[2nN^{1/2}_{2},2nN^{1/2}_{1}]\times[-nN_{1}\pi,nN_{1}\pi]\cup[-2nN^{1/2}_{1},-2nN^{1/2}_{2}]\times[-nN_{1}\pi,nN_{1}\pi].$$

To get an estimate for the integral on $A$ we use symmetry to focus on $h_{2}>0$ and  we separate the integral into three parts, $h_{2}\in[2\pi k_{0}+\frac{\pi}{2}, 2\pi (k_{\infty}+1)+\frac{\pi}{2} ]$ (with $k_{0}=k_{0}(N_{2},n)$ the smallest integer with $2\pi k_{0}+\frac{\pi}{2}\geq N_{2}n\pi$ and $k_{\infty}=k_{\infty}(N_{1},n)$ the biggest one such that $(k_{\infty}+1)2\pi+\frac{\pi}{2}\leq N_{1}n\pi$), $h_{2}\in [N_{2}n\pi,2\pi k_{0}+\frac{\pi}{2}]$ and $h_{2}\in[(k_{\infty}+1)2\pi+\frac{\pi}{2},N_{1}n\pi]$, and we estimate each part separately:

1) If  $h_{2}\in [2\pi k_{0}+\frac{\pi}{2}, 2\pi (k_{\infty}+1)+\frac{\pi}{2} ]$
\begin{align*}
    &|\int_{-2nN_{2}^{1/2}}^{2nN_{2}^{1/2}}\int_{2\pi k_{0}+\frac{\pi}{2}}^{ 2\pi (k_{\infty}+1)+\frac{\pi}{2} } sin(h_{2})\frac{h_{2}}{(h_{1}^2+h_{2}^2)^{3/2}}dh_{2}dh_{1}|\\
    &\leq\int_{-2nN_{2}^{1/2}}^{2nN_{2}^{1/2}} C\sum_{k=k_{0}}^{k_{\infty}}\frac{1}{(h_{1}^2+(k2\pi)^2)^{2}}dh_{1}\leq C\int_{-2nN_{2}^{1/2}}^{2nN_{2}^{1/2}} \sum_{k=k_{0}}^{k_{\infty}}\frac{1}{(h_{1}+k2\pi)^{4}}dh_{1}\\
    &\leq C\int_{-2nN_{2}^{1/2}}^{2nN_{2}^{1/2}} \int_{\frac{N_{2}n}{2}-\frac{5}{4}}^{\frac{N_{1}n}{2}}\frac{1}{(h_{1}+x2\pi)^{4}}dxdh_{1}
    \leq C\int_{-2nN_{2}^{1/2}}^{2nN_{2}^{1/2}} \frac{1}{(h_{1}+N_{2}n)^{3}}dh_{1}\\
    &\leq \frac{C}{N_{2}^{5/2}n^2}, \\
\end{align*}

2) If $h_{2}\in[N_{2}n\pi,2\pi k_{0}+\frac{\pi}{2}]$
\begin{align*}
    &|\int_{-2nN_{2}^{1/2}}^{2nN_{2}^{1/2}}\int_{N_{2}n\pi}^{2\pi k_{0}+\frac{\pi}{2}} sin(h_{2})\frac{h_{2}}{(h_{1}^2+h_{2}^2)^{3/2}}dh_{2}dh_{1}|\leq \frac{C}{N_{2}^{3/2}n},\\
\end{align*}

3) If $h_{2}\in[(k_{\infty}+1)2\pi+\frac{\pi}{2},N_{1}n\pi]$
\begin{align*}
    &|\int_{-2nN_{2}^{1/2}}^{2nN_{2}^{1/2}}\int_{2\pi (k_{\infty}+1)+\frac{\pi}{2}}^{N_{1}n\pi} sin(h_{2})\frac{h_{2}}{(h_{1}^2+h_{2}^2)^{3/2}}dh_{2}dh_{1}|\leq \frac{C}{N_{2}^{3/2}n}.\\
\end{align*}

For the integration in $B$ we use a similar trick, using parity to consider only $h_ {2}\geq 0$ and separating in the parts $h_{2}\leq \frac{5\pi}{2}$, $\frac{5\pi}{2}\leq h_{2}\leq 2\pi (k_{\infty}+1)+\frac{\pi}{2}$ and $ 2\pi (k_{\infty}+1)+\frac{\pi}{2}\leq h_{2}\leq N_{1}n\pi$, with $k_{\infty}=k_{\infty}(N_{1},n)$  the biggest integer such that $(k_{\infty}+1)2\pi+\frac{\pi}{2}\leq N_{1}n\pi$:
1) If $\frac{5\pi}{2}\leq h_{2}\leq 2\pi (k_{\infty}+1)+\frac{\pi}{2}$
\begin{align*}
    &|\int_{2nN_{2}^{1/2}}^{2nN_{1}^{1/2}}\int_{\frac{5\pi}{2}}^{ 2\pi (k_{\infty}+1)+\frac{\pi}{2}} sin(h_{2})\frac{h_{2}}{(h_{1}^2+h_{2}^2)^{3/2}}dh_{2}dh_{1}|\\
    &\leq\int_{2nN_{2}^{1/2}}^{2nN_{1}^{1/2}} C\sum_{k=1}^{k_{\infty}}\frac{1}{(h_{1}^2+(k2\pi)^2)^{2}}dh_{1}\leq C\int_{2nN_{2}^{1/2}}^{2nN_{1}^{1/2}} \sum_{k=1}^{k_{\infty}}\frac{1}{(h_{1}+k2\pi)^{4}}dh_{1}\\
    & \leq C\int_{2nN_{2}^{1/2}}^{2nN_{1}^{1/2}} \int_{0}^{\frac{N_{1}n}{2}}\frac{1}{(h_{1}+x2\pi)^{4}}dxdh_{1}
    \leq C\int_{2nN_{2}^{1/2}}^{2nN_{1}^{1/2}} \frac{1}{h_{1}^{3}}dh_{1} \\
    &\leq \frac{C}{N_{2}n^2}\\
\end{align*}
2) If $h_{2}\leq \frac{5\pi}{2}$
\begin{align*}
    &|\int_{2nN_{2}^{1/2}}^{2nN_{1}^{1/2}}\int_{0}^{\frac{5\pi}{2}} sin(h_{2})\frac{h_{2}}{(h_{1}^2+h_{2}^2)^{3/2}}dh_{2}dh_{1}|\leq \frac{C}{N_{2}n^2},\\
\end{align*}
3) If $ 2\pi (k_{\infty}+1)+\frac{\pi}{2}\leq h_{2}\leq N_{1}n\pi$
\begin{align*}
    &|\int_{2nN_{2}^{1/2}}^{2nN_{1}^{1/2}}\int_{2\pi (k_{\infty}+1)+\frac{\pi}{2}}^{N_{1}n\pi} sin(h_{2})\frac{h_{2}}{(h_{1}^2+h_{2}^2)^{3/2}}dh_{2}dh_{1}|\leq \frac{C}{N_{2}^{3/2}n}.\\
\end{align*}

Putting together the estimates in the regions $A$ and $B$ we have that $lim_{N\rightarrow \infty}I_{N,n}=C_{0}(n)$, and that $|C_{0}-I_{N,n}|\leq C N^{-1/2}$. The only thing left to do is to prove that $C_{0}$ is indeed different from $0$ and independent of $n$.

To prove that $C_{0}(n)$ is actually independent of $n$, it is enough to prove that, for two arbitrary integers $n_{1},$ $n_{2}$,

$$\text{lim}_{N\rightarrow\infty} I_{N,n_{1}}-I_{N,n_{2}}=0.$$
The proof is equivalent to that of (\ref{cauchy}), so we will omit it.

To prove that $C_{0}\neq 0$,  we start by focusing on the integral with respect to $h_{2}$ for any fixed $h_{1}$ on an interval of the form $[-K\pi,K\pi,]$ with $K\in\mathds{N}$

\begin{align*}
    &\int_{[-K\pi,K\pi]} sin(h_{2})\frac{h_{2}}{(h_{1}^2+h_{2}^2)^{3/2}}dh_{2}\\
    &=\int_{[-K\pi,K\pi]} cos(h_{2})\frac{1}{(h_{1}^2+h_{2}^2)^{1/2}}dh_{2}-\Big[sin(h_{2})\frac{1}{(h_{1}^2+h_{2}^2)^{1/2}}\Big]_{h_{2}=-K\pi}^{K\pi}\\
    &=\int_{[-K\pi,K\pi]} cos(h_{2})\frac{1}{(h_{1}^2+h_{2}^2)^{1/2}}dh_{2}=2\int_{[0,K\pi]} cos(h_{2})\frac{1}{(h_{1}^2+h_{2}^2)^{1/2}}dh_{2},\\
\end{align*}
and we can use this property to compute the integral in $[-2nN^{1/2},2nN^{1/2}]\times[-K\pi,K\pi]$ as

\begin{align*}
    &\int_{-2nN^{1/2}}^{2nN^{1/2}}\int_{0}^{K\pi} cos(h_{2})\frac{1}{(h_{1}^2+h_{2}^2)^{1/2}}dh_{2}dh_{1}=\int_{0}^{K\pi}\int_{-2nN^{1/2}}^{2nN^{1/2}} cos(h_{2})\frac{1}{(h_{1}^2+h_{2}^2)^{1/2}}dh_{1}dh_{2}\\
    &=\int_{0}^{K\pi}cos(h_{2})\int_{-\frac{2nN^{1/2}}{h_{2}}}^{\frac{2nN^{1/2}}{h_{2}}} \frac{1}{(x^2+1)^{1/2}}dxdh_{2}=2\int_{0}^{K\pi}cos(h_{2})log(\frac{2nN^{1/2}}{h_{2}}+(1+\frac{4n^2N}{h_{2}^2})^{1/2})dh_{2}\\
    &=2\int_{0}^{K\pi}cos(h_{2})(log(\frac{2nN^{1/2}}{h_{2}}+(1+\frac{4n^2N}{h_{2}^2})^{1/2})-log(\frac{4nN^{1/2}}{h_{2}}))dh_{2}\\
    &+2\int_{0}^{K\pi}cos(h_{2})log(\frac{4nN^{1/2}}{h_{2}})dh_{2}.\\
\end{align*}

And we can evaluate the last line  by checking the two integrals separately
\begin{align*}
    &\int_{0}^{K\pi}cos(h_{2})log(\frac{4nN^{1/2}}{h_{2}})dh_{2}=-\int_{0}^{K\pi}cos(h_{2})log(h_{2})dh_{2}\\
    &=-\Big[log(x)sin(x)-Si(x)\Big]^{K\pi}_{0}=Si(K\pi)>0,\\
\end{align*}
where $Si(x)\equiv \int_{0}^{x}\frac{sin(t)}{t}dt$ is the sine integral function, and
\begin{align*}
    &|\int_{0}^{K\pi}cos(h_{2})(log(\frac{2nN^{1/2}}{h_{2}}+(1+\frac{4n^2N}{h_{2}^2})^{1/2})-log(\frac{4nN^{1/2}}{h_{2}}))dh_{2}|\\
    &\leq \int_{0}^{K\pi}\frac{h_{2}}{4nN^{1/2}}((1+\frac{4n^2N}{h_{2}^2})^{1/2}-\frac{2nN^{1/2}}{h_{2}})dh_{2}\leq \frac{CK^3}{N}.\\
\end{align*}
Furthermore, we can bound the integral outside of the interval $h_{2}\in [-K\pi,K,\pi]$. The particular way we divide the integral depends on the parity of $K$ and $Nn$. Here we will obtain the bounds in the case $K$ even and $Nn$ odd, the other cases being analogous:
\begin{align*}
    & |\int_{-2nN^{1/2}}^{2nN^{1/2}}\int_{K\pi}^{Nn\pi} cos(h_{2})\frac{h_{2}}{(h_{1}^2+h_{2}^2)^{3/2}}dh_{2}dh_{1}| \\
    &\leq \int_{-2nN^{1/2}}^{2nN^{1/2}}\sum_{k=\frac{K}{2}}^{\frac{Nn-1}{2}-1} \frac{1}{(h_{1}^2+(2\pi k)^2)^{2}}dh_{1}+\int_{-2nN^{1/2}}^{2nN^{1/2}}\int_{(Nn-1)\pi}^{Nn\pi} \frac{1}{h_{1}^2+h_{2}^2}dh_{2}dh_{1}\\
    &\leq C\int_{0}^{2nN^{1/2}}\sum_{k=\frac{K}{2}}^{\frac{Nn-1}{2}-1} \frac{1}{(h_{1}+2\pi k)^{4}}dh_{1}+\frac{C}{N^{\frac32}}\\
    &\leq C\int_{0}^{2nN^{1/2}} \frac{1}{(h_{1}+2\pi (\frac{K}{2}-1))^{3}}dh_{1}+\frac{C}{N^{\frac32}}\leq (\frac{C}{K-2})^2+\frac{C}{N^{\frac32}}.\\
\end{align*}
Combining all these together we get that, for any $K\leq nN$

\begin{align*}
    &\int_{-2nN^{1/2}}^{2nN^{1/2}}\int_{0}^{nN\pi} cos(h_{2})\frac{1}{(h_{1}^2+h_{2}^2)^{1/2}}dh_{2}dh_{1}\\
    &\geq Si(K\pi)- (\frac{C}{K-2})^2-\frac{CK^3}{N}-\frac{C}{N^{\frac32}}\\
\end{align*}
and by taking $K$ big enough so that $\frac{Si(K\pi)}{2}- (\frac{C}{K-2})^2>0$ and then $N$ big enough so that $\frac{Si(K\pi)}{2}-\frac{CK^3}{N}-\frac{C}{N^{\frac12}}>0$ we are done.
\end{proof}

We can now combine both lemmas to obtain

\begin{lemma}\label{aproxfinal}
Given natural numbers $n$, $N$ and a $C^2$ function $g_{N}(.):\mathds{R}\rightarrow\mathds{R}$ with support in the interval $(1-\frac{N^{-1/2}}{2},1+\frac{N^{1/2}}{2})$ and $||g_{N}||_{C^{i}}\leq M N^{i/2}$ for $i=0,1,2$ , we have that there exists a constant $C_{0}\neq 0$ such that, for $r\in (1-N^{-1/2},1+N^{-1/2})$,

\begin{equation}
    |v_{r}(g_{N}(r)sin(Nn\alpha))-C_{0}cos(Nn\alpha)g_{N}(r)|\leq CMN^{-1/2}
\end{equation}
with $C$ depending on $n$ but not on $N$ or $g$.

Analogously, we have that
\begin{equation}\label{coslemma}
    |v_{r}(g_{N}(r)cos(Nn\alpha))+C_{0}sin(Nn\alpha)g_{N}(r)|\leq CMN^{-1/2}
\end{equation}
with $C$ depending only on $n$.

\end{lemma}
\begin{proof}

We already know by lemma \ref{velocidaderror} that
\begin{align}\label{cosa1}
    &|v_{r}(g_{N}(r)sin(Nn\alpha))-cos(Nn\alpha)\int_{\mathds{R}\times[-\pi,\pi]} \frac{r^2 \alpha' g_{N}(r+h)sin(Nn\alpha')}{|h^2+r^2(\alpha')^2|^{3/2}}d\alpha'dh|\\
    &\leq C||g_{N}(r)||_{L^{\infty}}N^{-1/2}\nonumber\\\nonumber
\end{align}
and, by a change of variables we have that

\begin{align*}
   &\int_{\mathds{R}\times[-\pi,\pi]} \frac{r^2 \alpha' g_{N}(r+h)sin(Nn\alpha')}{|h^2+r^2(\alpha')^2|^{3/2}}dhd\alpha'=\int_{\mathds{R}}\int_{[-\pi,\pi]} \frac{ \alpha' g_{N}(r+hr)sin(Nn\alpha')}{|h^2+\alpha'^2|^{3/2}}d\alpha'dh.\\
\end{align*}

But, for any fixed $r\in [1/2,3/2]$, we have $||g_{N}(r+rh)||_{C^{i}}\leq 2^{i} ||g_{N}(r+h)||_{C^{i}}$ and thus applying lemma \ref{velocidadcartesianas} we get

\begin{align}\label{cosa2}
   &|C_{0}g_{N}(r)-\int_{\mathds{R}}\int_{[-\pi,\pi]} \frac{ \alpha' g_{N}(r+hr)sin(Nn\alpha')}{|h^2+\alpha'^2|^{3/2}}d\alpha'dh|\leq 2CMN^{-1/2},\\\nonumber
\end{align}

and combining (\ref{cosa1}) and (\ref{cosa2}) we get the desired result.
We omit the proof of (\ref{coslemma}) since it is completely analogous to the previous result.
\end{proof}

All these results will allow us to compute locally the radial velocity with a small error, but we would like to also have decay as we go far away from $r=1$. For that we have the following lemma.

\begin{lemma}\label{decaimiento}
Given a $L^{\infty}$ function $g_{N}(.):\mathds{R}\rightarrow \mathds{R}$ with support in the interval $(1-\frac{N^{-1/2}}{2},1+\frac{N^{-1/2}}{2})$, and let $\theta$ be defined as

$$\theta(r,\alpha):=sin(Nn\alpha)g_{N}(r)$$
with $N,n$ natural numbers.

Then there is a constant $C$ such that, if $N$ is big enough and  $1/2>|r-1|\geq N^{-1/2}$ or $r\geq 3/2$, we have

$$|v_{r}(\theta)(r,\alpha)|\leq \frac{C||g_{N}||_{L^{\infty}}}{N^{3/2}|r-1|^2}.$$
\end{lemma}

\begin{proof}
To estimate $|v_{r}(\theta)(r,\alpha)|$ we will use expression (\ref{velocidadradial}) and therefore we need to find upper bounds for

$$|\int_{\mathds{R}\times[-\pi,\pi]}(r+h)^2\frac{sin(\alpha')g_{N}(r+h)sin(Nn\alpha')}{|h^2+2(r+h)r(1-cos(\alpha'))|^{3/2}}d\alpha'dh|.$$

Let us fix $h$ such that $r+h\in (1-\frac{N^{-1/2}}{2},1+\frac{N^{-1/2}}{2})$ and with $r\geq 1/2$. Using that $\int_{i\frac{2\pi}{Nn}}^{(i+1)\frac{2\pi}{Nn}}sin(Nn\alpha)d\alpha=0$ and a degree one Taylor expansion around $\alpha'=k\frac{2\pi}{Nn}+\frac{\pi}{Nn}$ for $\frac{sin(\alpha')}{|h^2+2(r+h)r(1-cos(\alpha'))|^{3/2}}$  we can bound the integral over a single period

\begin{align*}
    &|\int_{[k\frac{2\pi}{Nn},(k+1)\frac{2\pi}{Nn}]}\frac{sin(\alpha')sin(Nn\alpha')}{|h^2+2(r+h)r(1-cos(\alpha'))|^{3/2}}d\alpha'|\\
    &\leq \int_{[k\frac{2\pi}{Nn},(k+1)\frac{2\pi}{Nn}]}\frac{C}{Nn}\frac{1}{|h^2+2(r+h)r(1-cos(\alpha'))|^{3/2}}d\alpha'\\
    &\leq \frac{C}{(Nn)^2}\frac{1}{|h+ck\frac{2\pi}{Nn}|^{3}},\\
\end{align*}
with $c$ small and $C$ big, where we used that $r+h,r\geq 1/2$ and that there exists $c>0$ such that $\frac{1}{c}(1-cos(\alpha'))\geq (\alpha')^2$ if $\alpha'\in[-\pi,\pi]$. Adding over all the relevant periods we obtain

\begin{align*}
    &\sum_{k=0}^{nN}\frac{C}{(Nn)^2}\frac{1}{|h+ck\frac{2\pi}{Nn}|^{3}}\leq \int_{-1}^{Nn}\frac{C}{(Nn)^2}\frac{1}{|h+cx\frac{2\pi}{Nn}|^{3}}dx\\
    &\int_{-1}^{Nn}CNn\frac{1}{|h\frac{Nn}{2\pi c}+x|^{3}}dx\leq CNn \frac{1}{|h\frac{Nn}{2\pi c}-1|^{2}} \\
    &=\frac{C}{Nn}\frac{1}{|h-\frac{2\pi c}{Nn}|^{2}}\leq \frac{C}{Nn}\frac{1}{h^{2}}.\\
\end{align*}

Furthermore, since the support of $g_{N}(r)$ lies in $(1-\frac{N^{-1/2}}{2},1+\frac{N^{-1/2}}{2})$ and $|r-1|>N^{-1/2}$ we have that $|h|\geq \frac{|r-1|}{2}$, so, by integrating in $h$ we get

\begin{align*}
    &\int_{\mathds{R}}\frac{C}{Nn}(r+h)^2\frac{|g_{N}(r+h)|}{h^2}dh\\
    &\leq \int_{r+h-1\in(-\frac{N^{-1/2}}{2},\frac{N^{-1/2}}{2})}\frac{C}{Nn|r-1|^2}||g||_{L^{\infty}}dh\leq\frac{C}{N^{3/2}n|r-1|^2}||g||_{L^{\infty}}. \\
\end{align*}

\end{proof}

 \subsection{The pseudo-solution  method for ill-posedness in $C^{k}$}
 We say that a function $\bar{\theta}$ is a pseudo-solution to the SQG equation if it fulfils that

$$\frac{\partial \bar{\theta}}{\partial t} + v_{1}(\bar{\theta})\frac{\partial \bar{\theta}}{\partial x_{1}}+ v_{2}(\bar{\theta})\frac{\partial \bar{\theta}}{\partial x_{2}}+F(x,t)=0$$
$$v_{1}(\bar{\theta})=-\frac{\partial}{\partial x_{2}}(-\Delta)^{1/2} \bar{\theta}=-\mathcal{R}_{2}\bar{\theta}$$
$$v_{2}(\bar{\theta})=\frac{\partial}{\partial x_{1}}(-\Delta)^{1/2} \bar{\theta}=\mathcal{R}_{1}\bar{\theta}$$
$$\bar{\theta}(x,0)=\theta_{0}(x),$$
for some $F(x,t)$.
Obviously, this definition is not restrictive at all, since you can get essentially anything by choosing the right $F(x,t)$. We will, however, try and use the term pseudo-solution only for functions where $F(x,t)$ is small in a suitable norm.






With this in mind we are ready to discuss the initial conditions we will be considering. Namely, in polar coordinates we will work with initial conditions of the form

$$\lambda (f_{1}(r)+f_{2}(N^{1/2}(r-1)+1)\sum_{k=1}^{K}\frac{sin(Nk\alpha)}{N^2k^3})$$
with $N$ and $K$ natural numbers, $\lambda>0$ and where $f_{1}$ and $f_{2}$ satisfy the following conditions:
\begin{itemize}
    \item Both $f_{1}(r)$ and $f_{2}(r)$ are $C^{\infty}$ functions.
    \item $f_{2}(r)$ has its support contained in the interval $(1/2, 3/2)$ and $f_{1}$ has its support in $(1/2,3/2)\cup (M_{1},M_{2})$ with some $M_{1}$, $M_{2}$ big.
    \item $\frac{\partial f_{1}(r)}{\partial r}=1$ in $(3/4,5/4)$.
    \item  $f_{2}(r)=1$ in $(3/4,5/4)$.
    \item $\frac{\partial^{k}\frac{ v_{\alpha}(f_{1})(r)}{r}}{\partial r^{k}}$ is $0$ when $r=1$, $k=1,2$, where $v_{\alpha}(f_{1})$ is the velocity produced by $f_{1}$ in the angular direction.
\end{itemize}

We will use these pseudo-solutions to prove ill-posedness in $C^2$, and at the end of this section we will explain how to extend the proof to $C^k$, $k>2$.

It is not obvious that the properties we require for $f_{1}$ can be obtained, so we need the following lemma.

\begin{lemma}\label{seleccionvelocidadnagular}
There exists  a $C^{\infty}$  compactly supported function  $g(.):[0,\infty)\rightarrow\mathds{R}$ with support in $ (2,\infty)$ such that $\frac{\partial^{i}\frac{ v_{\alpha}(g(.))(r)}{r}}{\partial r^{i}}(r=1)=a_{i}$ with $i=1,2$ and $a_{i}$ arbitrary.
\end{lemma}

\begin{proof}
We start by considering a $C^{\infty}$  function $h(x):\mathds{R}\rightarrow\mathds{R}$ which is positive, with support in $(-1/2,1/2)$ and $\int hdx=1$. We define the family of functions

$$f_{n_{1},n_{2}}(r):=n_{1}h(n_{1}(r-n_{2})),$$
with  $n_{2}\geq n_{1}\geq 2$, $n_{1},n_{2}\in\mathds{N}$. These functions are $C^{\infty}$ for any $n_{1},n_{2}$ , and have their support in the interval $(n_{2}-\frac{1}{2n_{1}},n_{2}+\frac{1}{2n_{1}})$. Now let's consider the associated family of vectors

$$V=\cup V_{n_{1},n_{2}},$$
with

$$V_{n_{1},n_{2}}:=(\frac{\partial v_{\alpha}(f_{n_{1},n_{2}})}{\partial r}(r=1),\frac{\partial^{2}v_{\alpha}(f_{n_{1},n_{2}})}{\partial r^{2}}(r=1)).$$

Note that to prove our lemma it is sufficient to show that this family is in fact a base of the space $\mathds{R}^{2}$. Before we can actually prove that this is the case, we need to find expressions for $V_{n_{1},n_{2}}$. For our purposes it is enough to compute $\lambda_{n_{1},n_{2}}V_{n_{1},n_{1}}$ since this vectors will span the same space as long as $\lambda_{n_{1},n_{2}}\neq 0$.

To begin with, we need to start deducing the expression for $v_{\alpha}$. Proceeding in a similar way as for $v_{r}$, and for simplicity only considering the case when $\theta(r,\alpha)=f(r)$ we get
\begin{align}\label{velocidadangular}
    &v_{\alpha}(\theta(\ .\ ,\ .\ ))\nonumber\\
    &=P.V. \int_{\mathds{R}^2} \hat{x}^{\perp} \frac{(x-y)^{\perp}\theta(y)}{|x-y|^{3}}dy_{1}dy_{2}=P.V. \int_{\mathds{R}^2} \hat{x}^{\perp} \frac{(x-y)^{\perp}(\theta(y)-\theta(x))}{|x-y|^{3}}dy_{1}dy_{2}\nonumber\\
   & =P.V.\int_{\mathds{R}_{+}\times[-\pi,\pi]}r'\frac{(f(r')-f(r))(r-r'(cos(\alpha)cos(\alpha')+sin(\alpha)sin(\alpha')))}{|(rcos(\alpha)-r'cos(\alpha'))^2+(rsin(\alpha)-r'sin(\alpha'))^2|^{3/2}}d\alpha'dr'\nonumber\\
   &=P.V.\int_{\mathds{R}_{+}\times[-\pi,\pi]}r'\frac{r-r'cos(\alpha'-\alpha)}{|r^2+(r')^2-2rr'cos(\alpha-\alpha'))|^{3/2}}(f(r')-f(r))d\alpha'dr'.\\\nonumber
\end{align}

And, since we will be considering functions with support in $(2,\infty)$,  after relabeling $\alpha-\alpha'$ as $\alpha'$ we end up with the expression
$$P.V.\int_{2}^{\infty}\int_{-\pi}^{\pi}r'\frac{r-r'cos(\alpha')}{|r^2+(r')^2-2rr'cos(\alpha'))|^{3/2}}(f(r')-f(r))dr'd\alpha'.$$

Furthermore, if we write 
$$F(r,r',\alpha'):=r'\frac{r-r'cos(\alpha')}{|r^2+(r')^2-2rr'cos(\alpha'))|^{3/2}}$$
for $r=1$ we can use differentiation under the integral sign and obtain
\begin{align*}
    &\frac{\partial^{j} v_{\alpha}(f(.))}{\partial r^{j}}(r=1)=\int_{(2,\infty)\times[-\pi,\pi]}\frac{\partial^{j} F}{\partial r^{j}}(r,r',\alpha')(r=1)f(r')dr'd\alpha'.\\
\end{align*}

But for $f=f_{n_{1},n_{2}}$ we have that

\begin{align*}
    |\int_{(2,\infty)\times[-\pi,\pi]}\frac{\partial^{j} F}{\partial r^{j}}&(r,r',\alpha')f_{n_{1},n_{2}}(r')d\alpha'dr'-\int_{[-\pi,\pi]}\frac{\partial^{j} F}{\partial r^{j}}(r,n_{2},\alpha')d\alpha'|\leq \frac{C}{n_{1}},\\
\end{align*}
with $C$ depending on $r$ and, in particular, since $span(V)$ is a closed set, by taking $lim_{n_{1}\rightarrow\infty}V_{n_{1},n_{2}}$ we get that

$$(\int_{[-\pi,\pi]}\frac{\partial F}{\partial r}(1,n_{2},\alpha')d\alpha',\int_{[-\pi,\pi]}\frac{\partial^{2} F}{\partial r^{2}}(1,n_{2},\alpha')d\alpha')\in span(V)$$
Furthermore, we have that
\begin{align*}
    \frac{\partial F}{\partial r}(r,r'&,\alpha')(r=1)\\
    =r'\Big(\frac{1}{|r^2+(r')^2-2rr'cos(\alpha')|^{3/2}}&-\frac{3(r-r'cos(\alpha'))^2}{|r^2+(r')^2-2rr'cos(\alpha')|^{5/2}}\Big)(r=1)\\
\end{align*}
and, integrating with respect to $\alpha'$ we get
$$\int_{[-\pi,\pi]}\frac{\partial F}{\partial r}(r,n_{2},\alpha')(r=1)d\alpha'=-\frac{\pi}{(r')^2}(1+O(\frac{1}{r'})).$$

With the second derivative we obtain

\begin{align*}
    &\frac{\partial^{2} F}{\partial r^{2}}(r,r',\alpha')(r=1)\\
    &=-r'\Big(\frac{9(r-r'cos(\alpha'))}{|r^2+(r')^2-2rr'cos(\alpha')|^{5/2}}-\frac{15(r-r'cos(\alpha'))^3}{|r^2+(r')^2-2rr'cos(\alpha')|^{7/2}}\Big)(r=1).\\
\end{align*}

Now before we get into more details regarding this value, note that 

$$r'\Big(|\frac{9(r-r'cos(\alpha'))}{|r^2+(r')^2-2rr'cos(\alpha')|^{5/2}}-\frac{15(r-r'cos(\alpha'))^3}{|r^2+(r')^2-2rr'cos(\alpha')|^{7/2}}|\Big)(r=1)\leq \frac{1}{(r')^3}.$$

Therefore, we have that

$$(\frac{1}{(n_{2})^2}+O(\frac{1}{(n_{2})^{3}}),O(\frac{1}{(n_{2})^{3}}))\in span(V),$$
and again since $span(V)$ is a closed set, the vector $(1,0)$ belongs to $span(V)$. Now we only need to prove that there exists a point $r'$ such that

$$\int_{[-\pi,\pi]}-r'\Big(\frac{9(1-r'cos(\alpha'))}{|1+(r')^2-2r'cos(\alpha')|^{5/2}}-\frac{15(1-r'cos(\alpha'))^3}{|1+(r')^2-2r'cos(\alpha')|^{7/2}|}\Big)d\alpha'\neq 0,$$
so that we can find a vector $V_{n_{1},n_{2}}$ of the form $(a,b)$ with $b\neq 0$.
But, for example, using that, for $\delta>0$ and $r'$ big 
$$\frac{1}{(1+(r')^2-2r'cos(\alpha'))^{\delta}}+\frac{1}{(1+(r')^2+2r'cos(\alpha'))^{\delta}}-\frac{2}{(1+(r')^2)^{\delta}}\leq \frac{C}{(r')^{2(\delta+1)}}$$
one can check that
\begin{align*}
    &\int_{[-\pi,\pi]}-r'\Big(\frac{9(1-r'cos(\alpha'))}{|1+(r')^2-2r'cos(\alpha')|^{5/2}}-\frac{15(1-r'cos(\alpha'))^3}{|1+(r')^2-2r'cos(\alpha')|^{7/2}|}\Big)d\alpha'\\
    &=\frac{C}{(r')^4}+O(\frac{1}{(r')^5})\\
\end{align*}
with $C\neq 0$, and taking $r'$ big enough we are done.

\end{proof}

Therefore, to obtain $f_{1}$ with the desired properties, we first consider a radial $C^{\infty}$ function with support in $(\frac12,\frac32)$ and derivative $1$ in $(\frac34,\frac54)$ and then add a $C^{\infty}$ function with support in $[2,M]$ that cancels the derivatives of the velocity at $r=1$, and such a function exists thanks to lemma \ref{seleccionvelocidadnagular}.

Once we choose specific $f_{1}$ and $f_{2}$, this family of initial conditions has some useful properties that we will use later. First, for any fixed  $K$ and $\lambda$ our initial conditions are bounded in $H^{2+1/4}$ independently of the choice of $N$. Furthermore, the $C^2$ norm is bounded for any fixed $\lambda$  independently of both $N$ and $K$, and can be taken as small as we want by taking $\lambda$ small.

For any such initial conditions, we consider the associated pseudo-solution
\begin{equation}\label{pseudosolution}
    \bar{\theta}_{\lambda,K,N}(r,\alpha,t):=\lambda (f_{1}(r)+f_{2}(N^{1/2}(r-1)+1)\sum_{k=1}^{K}\frac{sin(Nk\alpha-\lambda tNk\frac{v_{\alpha}( f_{1})}{r}-\lambda C_{0}t)}{N^2k^3}),
\end{equation}

with $C_{0}$ the constant from lemma \ref{velocidadcartesianas} and \ref{aproxfinal}. We don't add subindexes for $f_{1}$ and $f_{2}$ since we consider them fixed from now on. Furthermore, the constants appearing in most results will also depend on $f_{1}$ and $f_{2}$, but since we consider them fixed we will not mention this.

This function for $N\geq 4$ satisfies 
\begin{equation}\label{pseudoevolucion}
    \frac{\partial \bar{\theta}_{\lambda,K,N}(r,\alpha,t)}{\partial t}+\frac{\partial \bar{\theta}_{\lambda,K,N}}{\partial \alpha} \frac{v_{\alpha}(\lambda f_{1})}{r}+\frac{\partial \lambda f_{1}}{\partial r} \bar{v}_{r}(\bar{\theta}_{\lambda,K,N})=0
\end{equation}
with 
$$\bar{v}_{r}(f(r)sin(k\alpha+g(r)))=C_{0}f(r)sin(k\alpha+g(r)+\frac{\pi}{2})$$
if $k\neq 0$, and $\bar{v}_{r}(f(r))=0$. Note that, for arbitrary fixed $T$, these functions satisfy that $||\bar{\theta}_{\lambda,K,N}||_{H^{2+1/4}}\leq C\lambda K$, with $C$ depending only on $T$.

Furthermore, we can rewrite (\ref{pseudoevolucion}) as

\begin{align*}
    &\frac{\partial \bar{\theta}_{\lambda,K,N}(r,\alpha,t)}{\partial t}+  \frac{\partial \bar{\theta}_{\lambda,K,N}}{\partial \alpha} \frac{v_{\alpha}(\bar{\theta}_{\lambda,K,N})}{r}+\frac{\partial \bar{\theta}_{\lambda,K,N}}{\partial r} v_{r}(\bar{\theta}_{\lambda,K,N})\\
    &+\frac{\partial \bar{\theta}_{\lambda,K,N}}{\partial \alpha} \frac{v_{\alpha}(\lambda f_{1}-\bar{\theta}_{\lambda,K,N})}{r}+\frac{\partial (\lambda f_{1}-\bar{\theta}_{\lambda,K,N})}{\partial r} v_{r}(\bar{\theta}_{\lambda,K,N})\\
    &+\frac{\partial \lambda f_{1}}{\partial r} (\bar{v}_{r}(\bar{\theta}_{\lambda,K,N})-v_{r}(\bar{\theta}_{\lambda,K,N}))=0\\
\end{align*}

Therefore $\bar{\theta}$ is a pseudo-solution with source term 
\begin{align*}
   &F_{\lambda,K,N}(x,t)=\\
   &\frac{\partial \bar{\theta}_{\lambda,K,N}}{\partial \alpha} \frac{v_{\alpha}(\lambda f_{1}-\bar{\theta}_{\lambda,K,N})}{r}+\frac{\partial (\lambda f_{1}-\bar{\theta}_{\lambda,K,N})}{\partial r}v_{r}(\bar{\theta}_{\lambda,K,N})+\frac{\partial \lambda f_{1}}{\partial r} (\bar{v}_{r}(\bar{\theta}_{\lambda,K,N})-v_{r}(\bar{\theta}_{\lambda,K,N})).\\
\end{align*}

Next we would like to prove that this source term is, indeed, small enough to obtain the desired results. We start by proving bounds on $L^{2}$ and in $H^{3}$ for $F_{\lambda,K,N}(x,t)$.

\begin{lemma}
For $t\in[0,T]$ and a pseudo-solution $\bar{\theta}_{\lambda,K,N}$ as in (\ref{pseudosolution}) the source term $F_{\lambda,K,N}(x,t)$ satisfies

$$||F_{\lambda,K,N}(x,t)||_{L^{2}}\leq CN^{-(2+3/4)}$$
with $C$ depending on  $K$, $\lambda$ and $T$.
\end{lemma}

\begin{proof}
We start bounding the term 
$\frac{\partial \lambda f_{1}}{\partial r} (\bar{v}_{r}(\bar{\theta}_{\lambda,K,N})-v_{r}(\bar{\theta}_{\lambda,K,N})).$ First we decompose  each function

\begin{align*}
    &\frac{sin(Nk\alpha-\lambda tN\frac{v_{\alpha}( f_{1})}{r}-\lambda C_{0}t)}{N^2k^3}\\
    &=\frac{sin(Nk\alpha)cos(\lambda tN\frac{v_{\alpha}( f_{1})}{r}+\lambda C_{0}t)-cos(Nk\alpha)sin(\lambda t\frac{v_{\alpha}( f_{1})}{r}+\lambda C_{0}t)}{N^2k^3}\\
\end{align*}
and using that $\frac{\partial^{k}\frac{v_{\alpha}( f_{1})}{r}}{\partial r^{k}}(r=1)=0 $ for $k=1,2$, then for $r\in(1-2N^{-1/2},1+2N^{-1/2})$ we have that

$$|\frac{\partial\frac{v_{\alpha}( f_{1})}{r}}{\partial r}|\leq\frac{C}{N}$$
and thus

$$||\frac{\partial cos(\lambda tNk\frac{v_{\alpha}( f_{1})}{r}+\lambda C_{0}t)}{\partial r}||_{L^{\infty}}\leq C$$
$$||\frac{\partial sin(\lambda tNk\frac{v_{\alpha}( f_{1})}{r}+\lambda C_{0}t)}{\partial r}||_{L^{\infty}}\leq C.$$
 Therefore, we can directly apply lemma \ref{aproxfinal} to obtain

\begin{align*}
    &|v_{r}(f_{2}(N^{1/2}(r-1)+1)\frac{sin(Nk\alpha)cos(\lambda t\frac{v_{\alpha}( f_{1})}{r}+\lambda C_{0}t)}{N^2k^3})\\
    &-\bar{v}_{r}(f_{2}(N^{1/2}(r-1)+1)\frac{sin(Nk\alpha)cos(\lambda t\frac{v_{\alpha}( f_{1})}{r}+\lambda C_{0}t)}{N^2k^3})|\leq \frac{C}{N^{5/2}k^{3}},\\
    &|v_{r}(f_{2}(N^{1/2}(r-1)+1)\frac{cos(Nk\alpha)sin(\lambda t\frac{v_{\alpha}( f_{1})}{r}+\lambda C_{0}t)}{N^2k^3})\\
    &-\bar{v}_{r}(f_{2}(N^{1/2}(r-1)+1)\frac{cos(Nk\alpha)sin(\lambda t\frac{v_{\alpha}( f_{1})}{r}+\lambda C_{0}t)}{N^2k^3})|\leq \frac{C}{N^{5/2}k^{3}}.\\
\end{align*}
With this we can estimate
\begin{align*}
    &\int_{1-N^{-1/2}}^{1+N^{-1/2}}\int_{-\pi}^{\pi}(\frac{\partial \lambda f_{1}}{\partial r} (v_{r}(\bar{\theta}_{\lambda,K,N})-\bar{v}_{r}(\bar{\theta}_{\lambda,K,N})))^2d\alpha dr\leq (||\frac{\partial f_{1}}{\partial r} ||_{L^{\infty}})^2\frac{C}{N^{5+1/2}}.\\
\end{align*}

For $r\in (1/2, 1-N^{-1/2})\cup(1+N^{-1/2},\infty)$, we use that $\bar{v}$ is zero in those points and lemma \ref{decaimiento} to obtain

\begin{align*}
    &\int_{1/2}^{1-N^{-1/2}}\int_{-\pi}^{\pi} (\frac{\partial \lambda f_{1}}{\partial r} (v_{r}(\bar{\theta}_{\lambda,K,N})-\bar{v}_{r}(\bar{\theta}_{\lambda,K,N})))^2d\alpha dr\\
    &\leq C (||\frac{\partial f_{1}}{\partial r} ||_{L^{\infty}})^2\int_{1/2}^{1-N^{-1/2}}\int_{-\pi}^{\pi} (\frac{||f_{2}||_{L^{\infty}}}{N^{3/2}|r-1|^2} )^2d\alpha dr \\
    &\leq \frac{C}{N^{5+1/2}}(||\frac{\partial f_{1}}{\partial r} ||_{L^{\infty}})^2 (||f_{2}||_{L^{\infty}} )^2\\
\end{align*}
and similarly

\begin{align*}
    &\int_{1+N^{-1/2}}^{\infty}\int_{-\pi}^{\pi} (\frac{\partial \lambda f_{1}}{\partial r} (v_{r}(\bar{\theta})-\bar{v}_{r}(\bar{\theta})))^2d\alpha dr\leq \frac{C}{N^{5+1/2}}(||\frac{\partial f_{1}}{\partial r} ||_{L^{\infty}})^2 (||f_{2}||_{L^{\infty}} )^2.\\
\end{align*}

Combining all of these inequalities we get

$$||\frac{\partial \lambda  f_{1}}{\partial r} (v_{r}(\bar{\theta})-\bar{v}_{r}(\bar{\theta}))||_{L^{2}}\leq \frac{C}{N^{2+3/4}}$$

with $C$ depending on $\lambda$, $K$ and $T$.

For the term $\frac{\partial \bar{\theta}_{\lambda,K,N}}{\partial \alpha} \frac{v_{\alpha}(\lambda f_{1}-\bar{\theta}_{\lambda,K,N})}{r}$ we simply use $||\frac{\partial \bar{\theta}_{\lambda,K,N}}{\partial \alpha}||_{L^{\infty}}\leq \frac{C}{N}$ and $$||\frac{v_{\alpha}(\lambda f_{1}-\bar{\theta}_{\lambda,K,N})}{r}1_{\text{supp}(\bar{\theta}_{\lambda,K,N})}||_{L^2}\leq \frac{C}{N^{2+1/4}}$$ so

$$||\frac{\partial \bar{\theta}_{\lambda,K,N}}{\partial \alpha} \frac{v_{\alpha}(\lambda  f_{1}-\bar{\theta}_{\lambda,K,N})}{r}||_{L^{2}}\leq \frac{C}{N^{3+1/4}}.$$

Similarly for $\frac{\partial (\lambda f_{1}-\bar{\theta}_{\lambda,K,N})}{\partial r}v_{r}(\bar{\theta}_{\lambda,K,N})$ we have that $$||v_{r}(\bar{\theta}_{\lambda,K,N})||_{L^{2}}\leq \frac{C}{N^{2+1/4}} \quad \text{and}\quad  ||\frac{\partial (\bar{\theta}_{\lambda,K,N}-\lambda f_{1})}{\partial r}||_{L^{\infty}}\leq \frac{C}{N}$$ so

$$||\frac{\partial (\bar{\theta}_{\lambda,K,N}-\lambda f_{1})}{\partial r}v_{r}(\bar{\theta}_{\lambda,K,N})||_{L^{2}}\leq \frac{C}{N^{3+1/4}}$$
and we are done.
\end{proof}

\begin{lemma}
For $t\in [0,T]$, given a  pseudo-solution $\bar{\theta}_{\lambda,K,N}$ as in (\ref{pseudosolution}) the source term $F_{\lambda,K,N}(x,t)$ satisfies

$$||F_{\lambda,K,N}(x,t)||_{H^{3}}\leq CN^{3/4}$$
with $C$ depending on $K$, $\lambda$, and $T$.
\end{lemma}
\begin{proof}
To prove this we will use that, given the product of two functions, we have

$$||fg||_{H^{3}}\leq C (||f||_{L^{\infty}}||g||_{H^{3}}+||f||_{C^{1}}||g||_{H^2}+||f||_{C^2}||g||_{H^{1}}+||f||_{C^3}||g||_{L^2}).$$

Furthermore, for the pseudo-solutions considered, we have that $||\bar{\theta}_{\lambda,K,N}-\lambda f_{1}||_{C^{k}}\leq CN^{k-2}$, $||\bar{\theta}_{\lambda,K,N}-\lambda f_{1}||_{H^{k}}\leq CN^{k-2-1/4}$, $||\lambda f_{1}||_{C^{k}}\leq C$ with the constants $C$ depending on $k$, $\lambda$ and $K$.

Therefore we have that, using the bounds for the support of $\bar{\theta}_{\lambda,K,N}$

\begin{align*}
    &||\frac{\partial \bar{\theta}_{\lambda,K,N}}{\partial \alpha} \frac{v_{\alpha}(\lambda f_{1}-\bar{\theta}_{\lambda,K,N})}{r}||_{H^{3}}\\
    &\leq C (||\frac{\partial \bar{\theta}_{\lambda,K,N}}{\partial \alpha}||_{L^{\infty}}||v_{\alpha}(\lambda f_{1}-\bar{\theta}_{\lambda,K,N})||_{H^{3}}+||\frac{\partial \bar{\theta}_{\lambda,K,N}}{\partial \alpha}||_{C^{1}}||v_{\alpha}(\lambda f_{1}-\bar{\theta}_{\lambda,K,N})||_{H^2}\\
    &+||\frac{\partial \bar{\theta}_{\lambda,K,N}}{\partial \alpha}||_{C^2}||v_{\alpha}(\lambda f_{1}-\bar{\theta}_{\lambda,K,N})||_{H^{1}}+||\frac{\partial \bar{\theta}_{\lambda,K,N}}{\partial \alpha}||_{C^3}||v_{\alpha}(\lambda f_{1}-\bar{\theta}_{\lambda,K,N})||_{L^2})\\
    &\leq C  N^{-1/4},\\
\end{align*}
analogously

\begin{align*}
    &||\frac{\partial (\bar{\theta}_{\lambda,K,N}-\lambda f_{1})}{\partial r}v_{r}(\bar{\theta}_{\lambda,K,N})||_{H^{3}}\\
    &\leq C (||\frac{\partial (\bar{\theta}_{\lambda,K,N}-\lambda f_{1})}{\partial r}||_{L^{\infty}}||v_{r}(\bar{\theta}_{\lambda,K,N})||_{H^{3}}+||\frac{\partial (\bar{\theta}_{\lambda,K,N}-\lambda f_{1})}{\partial r}||_{C^{1}}||v_{r}(\bar{\theta}_{\lambda,K,N})||_{H^2}\\
    &+||\frac{\partial (\bar{\theta}_{\lambda,K,N}-\lambda f_{1})}{\partial r}||_{C^2}||v_{r}(\bar{\theta}_{\lambda,K,N})||_{H^{1}}+||\frac{\partial (\bar{\theta}_{\lambda,K,N}-\lambda f_{1})}{\partial r}||_{C^3}||v_{r}(\bar{\theta}_{\lambda,K,N})||_{L^2})\\
    &\leq C N^{-1/4},\\
\end{align*}
and finally

\begin{align*}
    &||\frac{\partial \lambda f_{1}}{\partial r} (v_{r}(\bar{\theta}_{\lambda,K,N})-\bar{v}_{r}(\bar{\theta}_{\lambda,K,N}))||_{H^{3}}\\
    & \leq C (||\frac{\partial \lambda f_{1}}{\partial r} ||_{L^{\infty}}||(v_{r}(\bar{\theta}_{\lambda,K,N})-\bar{v}_{r}(\bar{\theta}_{\lambda,K,N}))||_{H^{3}}+||\frac{\partial \lambda f_{1}}{\partial r} ||_{C^{1}}||(v_{r}(\bar{\theta}_{\lambda,K,N})-\bar{v}_{r}(\bar{\theta}_{\lambda,K,N}))||_{H^2}\\
    &+||\frac{\partial \lambda f_{1}}{\partial r} ||_{C^2}||(v_{r}(\bar{\theta}_{\lambda,K,N})-\bar{v}_{r}(\bar{\theta}_{\lambda,K,N}))||_{H^{1}}+||\frac{\partial \lambda f_{1}}{\partial r} ||_{C^3}||(v_{r}(\bar{\theta}_{\lambda,K,N})-\bar{v}_{r}(\bar{\theta}_{\lambda,K,N}))||_{L^2})\\
    &\leq C N^{3/4}
\end{align*}

\end{proof}

We can combine these two lemmas and use the interpolation inequality for Sobolev spaces to obtain that


$$||F||_{H^{2+1/4}}\leq C (N^{-(2+3/4)})^{1/4}(N^{3/4})^{3/4}\leq CN^{-1/8}.$$

With this, we are ready to study how the real solution behaves. If we define $$\Theta_{\lambda,K,N}=\theta_{\lambda,K,N}-\bar{\theta}_{\lambda,K,N},$$ with $\theta_{\lambda,K,N}$ the only $H^{2+\frac14}$ solution to the SQG equation with the same initial conditions as $\bar{\theta}_{\lambda,K,N}$, we have that

\begin{align}\label{evoluciondiferencia}
    &\frac{\partial \Theta_{\lambda,K,N}}{\partial t}+ v_{1}(\Theta_{\lambda,K,N})\frac{\partial \Theta_{\lambda,K,N}}{\partial x_{1}}+ v_{2}(\Theta_{\lambda,K,N})\frac{\partial \Theta_{\lambda,K,N}}{\partial x_{2}}\nonumber\\
    &+v_{1}(\Theta_{\lambda,K,N})\frac{\partial \bar{\theta}_{\lambda,K,N}}{\partial x_{1}}+ v_{2}(\Theta_{\lambda,K,N})\frac{\partial \bar{\theta}_{\lambda,K,N}}{\partial x_{2}}\\
    &+ v_{1}(\bar{\theta}_{\lambda,K,N})\frac{\partial \Theta_{\lambda,K,N}}{\partial x_{1}}+ v_{2}(\bar{\theta}_{\lambda,K,N})\frac{\partial \Theta_{\lambda,K,N}}{\partial x_{2}}-F_{\lambda,K,N}(x,t)=0,\nonumber\\\nonumber
\end{align}
and we have the following results regarding the evolution of $\Theta_{\lambda,K,N}$.

\begin{lemma}\label{errorl2}
Let $\Theta_{\lambda,K,N}$ defined as in (\ref{evoluciondiferencia}), then if $\theta_{\lambda,K,N}$ exists for $t\in [0,T]$, we have that

$$||\Theta_{\lambda,K,N}(x,t)||_{L^2}\leq \frac{Ct}{N^{(2+3/4)}}$$
with $C$ depending on $\lambda$, $K$ and $T$.
\end{lemma}

\begin{proof}
We start by noting that
\begin{align*}
    &\frac{\partial }{\partial t} \frac{||\Theta_{\lambda,K,N}||_{L^2}^2}{2}=-\int_{\mathds{R}^2}\Theta_{\lambda,K,N}\\
    &\Big((v_{1}(\Theta_{\lambda,K,N})+v_{1}(\bar{\theta}_{\lambda,K,N}))\frac{\partial \Theta_{\lambda,K,N}}{\partial x_{1}}+ (v_{2}(\Theta_{\lambda,K,N})+v_{2}(\bar{\theta}_{\lambda,K,N}))\frac{\partial \Theta_{\lambda,K,N}}{\partial x_{2}}\\
    &+v_{1}(\Theta_{\lambda,K,N})\frac{\partial \bar{\theta}_{\lambda,K,N}}{\partial x_{1}}+ v_{2}(\Theta_{\lambda,K,N})\frac{\partial \bar{\theta}_{\lambda,K,N}}{\partial x_{2}}-F_{\lambda,K,N}(x,t)\Big)dx,\\
\end{align*}
but, by incompressibility we have that 

\begin{align*}
\int_{\mathds{R}^2}\Theta_{\lambda,K,N}\Big((v_{1}(\Theta_{\lambda,K,N})+v_{1}(\bar{\theta}_{\lambda,K,N}))\frac{\partial \Theta_{\lambda,K,N}}{\partial x_{1}}+ (v_{2}(\Theta_{\lambda,K,N})+v_{2}(\bar{\theta}_{\lambda,K,N}))\frac{\partial \Theta_{\lambda,K,N}}{\partial x_{2}}\Big) dx=0,\\
\end{align*}
 and therefore we get that
 
 \begin{align*}
    &\frac{\partial }{\partial t} \frac{||\Theta_{\lambda,K,N}||_{L^2}^2}{2}\\
    &\leq |\int_{\mathds{R}^2}\Theta_{\lambda,K,N}\Big( v_{1}(\Theta_{\lambda,K,N})\frac{\partial \bar{\theta}_{\lambda,K,N}}{\partial x_{1}}+ v_{2}(\Theta_{\lambda,K,N})\frac{\partial \bar{\theta}_{\lambda,K,N}}{\partial x_{2}}+F_{\lambda,K,N}(x,t)\Big)dx|\\
    &\leq ||\Theta_{\lambda,K,N}||_{L^2} \Big(||\Theta_{\lambda,K,N}||_{L^2}||\bar{\theta}_{\lambda,K,N}||_{C^1}+||F_{\lambda,K,N}(x,t)||_{L^{2}}\Big),\\
\end{align*}
and using that $||F_{\lambda,K,N}||_{L^{2}}\leq \frac{C}{N^{(2+3/4)}}$, $||\bar{\theta}_{\lambda,K,N}||_{C^1}\leq C$ and integrating we get that

$$||\Theta_{\lambda,K,N}||_{L^2}\leq \frac{C (e^{Ct}-1)}{N^{(2+3/4)}}.$$

\end{proof}

\begin{lemma}\label{errorhs}
Let $\Theta_{\lambda,K,N}$ defined as in (\ref{evoluciondiferencia}), then for $N$ big enough, $\theta_{\lambda,K,N}$ exists for $t\in [0,T]$ and

$$||\Theta_{\lambda,K,N}(x,t)||_{H^{2+1/4}}\leq \frac{Ct}{N^{1/8}}$$
with $C$ depending on $\lambda$, $K$ and $T$.
\end{lemma}

\begin{proof}

It is enough to prove that 

$$||D^{2+1/4}\Theta_{\lambda,K,N}||_{L^2}\leq \frac{Ct}{N^{1/8}}$$
since 
$$||f||_{H^{s}}\leq C(||D^{s}f||_{L^2}+||f||_{L^2})$$
with $D^{s}=(-\Delta)^{s/2}$ and we already have the result

$$||\Theta_{\lambda,K,N}||_{L^2}\leq \frac{Ct}{N^{(2+3/4)}}.$$

We will use the following result found in \cite{Dongli}.

\begin{lemma}\label{katoponce}
Let $s > 0$. Then for any $s_{1}, s_{2} \geq 0$ with $s_{1} + s_{2} = s$, and any $f$, $g \in \mathcal{S}(\mathds{R}^2)$, the following holds:
\begin{equation}\label{Dsconmutador}
    ||D^{s}(fg)-\sum_{|\mathbf{k}|\leq s_{1}}\frac{1}{\mathbf{k}!}\partial^{\mathbf{k}}f D^{s,\mathbf{k}}g-\sum_{|\mathbf{j}|\leq s_{2}}\frac{1}{\mathbf{j}!}\partial^{\mathbf{j}}g D^{s,\mathbf{j}}f||_{L^{2}}\leq C ||D^{s_{1}}f||_{L^{2}} ||D^{s_{2}}g||_{BMO}
\end{equation}
where $\mathbf{j}$ and $\mathbf{k}$ are multi-indexes, $\partial^{\mathbf{j}}=\frac{\partial}{\partial x^{j_{1}}_{1}\partial x^{j_{2}}_{2}}$, $\partial_{\xi}^{\mathbf{j}}=\frac{\partial}{\partial \xi^{j_{1}}_{1}\partial \xi^{j_{2}}_{2}}$ and $ D^{s,\mathbf{j}}$ is defined using

$$\widehat{D^{s,\mathbf{j}}f}(\xi)=\widehat{D^{s,\mathbf{j}}}(\xi)\hat{f}(\xi)$$
$$\widehat{D^{s,\mathbf{j}}}(\xi)=i^{-|\mathbf{j}|}\partial^{\mathbf{j}}_{\xi}(|\xi|^s).$$
\end{lemma}

Although this result is for functions in the  Schwartz space $\mathcal{S}$, since we only consider compactly supported functions we can apply it to functions in $H^s$. We will consider $s=2+1/4$, although we will just write $s$ for compactness of notation.

Then

\begin{align*}
    &\frac{\partial}{\partial t} \frac{||D^{s}\Theta_{\lambda,K,N}||_{L^{2}}^2}{2}=- \int_{\mathds{R}^2} D^{s}\Theta_{\lambda,K,N} \\
        &D^{s}\Big((v_{1}(\Theta_{\lambda,K,N})+v_{1}(\bar{\theta}_{\lambda,K,N}))\frac{\partial \Theta_{\lambda,K,N}}{\partial x_{1}}+ (v_{2}(\Theta_{\lambda,K,N})+v_{2}(\bar{\theta}_{\lambda,K,N}))\frac{\partial \Theta_{\lambda,K,N}}{\partial x_{2}}\\
    &+v_{1}(\Theta_{\lambda,K,N})\frac{\partial \bar{\theta}_{\lambda,K,N}}{\partial x_{1}}+ v_{2}(\Theta_{\lambda,K,N})\frac{\partial \bar{\theta}_{\lambda,K,N}}{\partial x_{2}}+F_{\lambda,K,N})(x,t)\Big)dx.
\end{align*}

We will focus for now on
 
 $$\int_{\mathds{R}^2} D^{s}\Theta_{\lambda,K,N}
D^{s}\Big(v_{1}(\bar{\theta}_{\lambda,K,N}))\frac{\partial \Theta_{\lambda,K,N}}{\partial x_{1}}+ v_{2}(\bar{\theta}_{\lambda,K,N})\frac{\partial \Theta_{\lambda,K,N}}{\partial x_{2}}\Big)dx.$$

 Applying (\ref{Dsconmutador}) with $s_{2}=1$, $g=v_{i}(\bar{\theta}_{\lambda,K,N}))$, $f=\frac{\partial \Theta_{\lambda,K,N}}{\partial x_{i}}$, $i=1,2$ we would get that

\begin{align*}
    &(D^{s}\Theta_{\lambda,K,N},D^{s}(fg)-\sum_{|\mathbf{k}|\leq s_{1}}\frac{1}{\mathbf{j}!}\partial^{\mathbf{j}}f D^{s,\mathbf{j}}g-\sum_{|\mathbf{k}|\leq s_{2}}\frac{1}{\mathbf{k}!}\partial^{\mathbf{k}}g D^{s,\mathbf{k}}f)_{L^{2}}\\
    &\leq C ||D^{s}\Theta_{\lambda,K,N}||_{L^{2}}||D^{s_{1}}f||_{L^{2}} ||D^{s_{2}}g||_{BMO}\\
    &\leq C||D^{s}\Theta_{\lambda,K,N}||_{L^{2}} ||\bar{\theta}_{\lambda,K,N}||_{H^{s}}||\Theta_{\lambda,K,N}||_{H^s}.\\
\end{align*}

Furthermore we have that

\begin{align*}
    &(D^{s}\Theta_{\lambda,K,N}, D^{s}(\frac{\partial \Theta_{\lambda,K,N}}{\partial x_{1}})v_{1}(\bar{\theta}_{\lambda,K,N})+ D^{s}(\frac{\partial \Theta_{\lambda,K,N}}{\partial x_{2}})v_{2}(\bar{\theta}_{\lambda,K,N}))_{L^{2}}\\
    &=\frac{1}{2}\int_{\mathds{R}^2} \frac{\partial}{\partial x_{1}}(D^{s}\Theta_{\lambda,K,N})^2
v_{1}(\bar{\theta}_{\lambda,K,N}))+ \frac{\partial}{\partial x_{2}}(D^{s}\Theta_{\lambda,K,N})^2v_{2}(\bar{\theta}_{\lambda,K,N})dx=0\\
\end{align*}
and, for $i=1,2$, using that the operators $D^{s,c}$ are continuous from $H^{a}$ to $H^{a-s+c}$, we have the following three estimates

1)
\begin{align*}
    &|(D^{s}\Theta_{\lambda,K,N},\sum_{|\mathbf{k}|=1}\frac{1}{\mathbf{k}!}\partial^{\mathbf{k}}v_{i}(\bar{\theta}_{\lambda,K,N})) D^{s,\mathbf{k}}\frac{\partial \Theta_{\lambda,K,N}}{\partial x_{i}})_{L^{2}}|\\
    &\leq C ||D^{s}\Theta_{\lambda,K,N}||_{L^{2}}||v_{i}(\bar{\theta}_{\lambda,K,N})||_{H^{2+\epsilon}}||\Theta_{\lambda,K,N}||_{H^s}\\
    &\leq C ||D^{s}\Theta_{\lambda,K,N}||_{L^{2}}||\bar{\theta}_{\lambda,K,N}||_{H^{s}}||\Theta_{\lambda,K,N}||_{H^s}\\
\end{align*}

2)
\begin{align*}
    &|(D^{s}\Theta_{\lambda,K,N},\sum_{|\mathbf{j}|=1}\frac{1}{\mathbf{j}!}\partial^{\mathbf{j}}\frac{\partial \Theta_{\lambda,K,N}}{\partial x_{i}} D^{s,\mathbf{j}}v_{i}(\bar{\theta}_{\lambda,K,N}))_{L^{2}}|\\
    &\leq C\sum_{|\mathbf{j}|=1} ||D^{s}\Theta_{\lambda,K,N}||_{L^{2}}||\frac{1}{\mathbf{j}!}\partial^{\mathbf{j}}\frac{\partial \Theta_{\lambda,K,N}}{\partial x_{i}}||_{L^{2/(3-s)}}||D^{s,\mathbf{j}}v_{i}(\bar{\theta}_{\lambda,K,N}))||_{L^{2/(s-2)}}\\
    &\leq C ||D^{s}\Theta_{\lambda,K,N}||_{L^{2}}||\Theta_{\lambda,K,N}||_{H^s}||\bar{\theta}_{\lambda,K,N}||_{H^{s}},\\
\end{align*}


3)
\begin{align*}
    &|(D^{s}\Theta_{\lambda,K,N}, \frac{\partial \Theta_{\lambda,K,N}}{\partial x_{i}} D^{s}v_{i}(\bar{\theta}_{\lambda,K,N}))_{L^{2}}|\\
    &\leq C ||D^{s}\Theta_{\lambda,K,N}||_{L^{2}}||\Theta_{\lambda,K,N}||_{H^s}||\bar{\theta}_{\lambda,K,N}||_{H^{s}}.\\
\end{align*}

Most of the other terms are bounded in a similar way without any complication, although a comment needs to be made about bounding the terms

\begin{align*}
    & \int_{\mathds{R}^2} D^s(\Theta_{\lambda,K,N}) \Big(v_{1}(\Theta_{\lambda,K,N})\frac{\partial D^{s}\bar{\theta}_{\lambda,K,N}}{\partial x_{1}}+ v_{2}(\Theta_{\lambda,K,N})\frac{\partial D^{s}\bar{\theta}_{\lambda,K,N}}{\partial x_{2}}\Big)dx.
\end{align*}

At first glance one could think that, since we are considering $\bar{\theta}_{\lambda,K,N}$ bounded in $H^{2+1/4}$ but not in higher order spaces, we could have a problem bounding this integral. However, we actually have that 
$$||\frac{\partial D^{s}\bar{\theta}_{\lambda,K,N}}{\partial x_{i}}||_{L^{\infty}}\leq C ||D^{s}\bar{\theta}_{\lambda,K,N}||_{H^{2+\epsilon}}\leq C||\bar{\theta}_{\lambda,K,N}||_{H^{4+1/4+\epsilon}}\leq C N^{2+\epsilon}$$

$$||v_{i}(\Theta_{\lambda,K,N})||_{L^{2}}\leq CTN^{-(2+3/4)}$$
and thus

\begin{align*}
    & |\int_{\mathds{R}^2} D^s(\Theta_{\lambda,K,N}) \Big(v_{1}(\Theta_{\lambda,K,N})\frac{\partial D^{s}\bar{\theta}_{\lambda,K,N}}{\partial x_{1}}+ v_{2}(\Theta_{\lambda,K,N})\frac{\partial D^{s}\bar{\theta}_{\lambda,K,N}}{\partial x_{2}}\Big)dx|\\
    &\leq CT ||D^s(\Theta_{\lambda,K,N})||_{L^{2}} N^{-3/4+\epsilon}\leq  CT ||D^s(\Theta_{\lambda,K,N})||_{L^{2}} N^{-1/8}, \\
\end{align*}
and combining all of this together plus similar bounds for the other terms, and using 
$$||\bar{\theta}_{\lambda,K,N}||_{H^{s}}\leq C,\quad ||F_{\lambda,K,N}||_{H^{s}}\leq C N^{-1/8}$$ with $C$ depending on $\lambda,$ $K$ and $T$, we get

$$\frac{\partial}{\partial t}||D^{s}\Theta_{\lambda,K,N}||^{2}_{L^{2}}\leq ||D^{s}\Theta_{\lambda,K,N}||_{L^{2}}(CN^{-1/8}+C||\Theta_{\lambda,K,N}||_{H^{s}}+C||\Theta_{\lambda,K,N}||_{H^{s}}^2)$$
which gives us, using 
$$||\Theta_{\lambda,K,N}||_{H^{s}}\leq C(||\Theta_{\lambda,K,N}||_{L^{2}}+||D^{s}\Theta_{\lambda,K,N}||_{L^{2}})\leq C(||D^{s}\Theta_{\lambda,K,N}||_{L^{2}}+N^{-(2+3/4)})$$
that

$$\frac{\partial}{\partial t}||D^{s}\Theta_{\lambda,K,N}||_{L^{2}}\leq (CN^{-1/8}+C||D^{s}\Theta_{\lambda,K,N}||_{L^{2}}+C||D^{s}\Theta_{\lambda,K,N}||_{L^{2}}^2).$$

Now, we restrict ourselves to $[0,T_{*}]$, with $T_{*}$ the smallest time such that $||D^{s}\Theta_{\lambda,K,N}||_{L^{2}}\leq 1$ (or $T$ if $T_{*}$ is bigger than $T$ or it does not exist). Integrating  for those times we get

$$||D^{s}\Theta_{\lambda,K,N}||_{L^2}\leq \frac{C (e^{Ct}-1)}{N^{1/8}},$$
and since for $N$ big enough we have that $T\leq T_{*}$ we are done.

\end{proof}


Now we are finally prepared to prove strong ill-posedness in $C^2$ for the SQG equation.

\begin{theorem}\label{strongillpos}
For any $c_{0}>0$, $M>0$ and $t_{*}>0$, we can find a $C^2\cap H^{2+1/4}$ function $\theta_{0}(x)$  with $||\theta_{0}(x)||_{C^2}\leq c_{0}$ such that the only solution $\theta(x,t)\in H^{2+\frac14}$ to the SQG problem (\ref{SQG}) with initial  conditions $\theta_{0}(x)$ that satisfy $||\theta(x,t_{*})||_{C^{2}}\geq M  c_{0}$.
\end{theorem}

\begin{proof}
We will prove this by constructing a solution with the desired properties.
We fix arbitrary $c_{0}>0$, $M>0$ and $t_{*}$, and consider the pseudo-solutions $\bar{\theta}_{\lambda,K,N}$. First, note that, for any $N$, $K$ natural numbers, for $\lambda>0$  small enough our family of pseudo-solutions has a small initial norm in $C^2$, so we consider $\lambda=\lambda_{0}$ small so that $||\theta_{\lambda_{0},K,N}(x,0)||_{C^2}\leq c_{0}$ for all $K$, $N$ natural and such that $|\lambda_{0}C_{0}t_{*}|\leq \frac{\pi}{2}$.

These pseudo-solutions fulfill that, at time $t$, for $\alpha= \lambda_{0} t \frac{v_{\alpha}( f_{1})}{r}$
\begin{align*}
    |\frac{\partial^2 \bar{\theta}_{\lambda_{0,K,N}}(x,t)}{\partial \alpha^2}&|=|\lambda_{0} f_{2}(N^{1/2}(r-1)+1) \sum_{k=1}^{K}\frac{sin(Nk\alpha-\lambda_{0} tNn\frac{v_{\alpha}( f_{1})}{r}-\lambda_{0} C_{0}t)}{k})|\\
    &=|\lambda_{0} f_{2}(N^{1/2}(r-1)+1) \sum_{k=1}^{K}\frac{sin(-\lambda_{0} C_{0}t)}{k})|\\
    &\geq \lambda_{0}|f_{2}(N^{1/2}(r-1)+1)|ln(K)|sin(-\lambda_{0} C_{0}t)|.\\
\end{align*}

Furthermore, we can find $c>0$ small such that , for $\alpha\in[\lambda_{0} t \frac{v_{\alpha}( f_{1})}{r}-c\frac{2\pi}{NK},\lambda_{0} t \frac{v_{\alpha}( f_{1})}{r}+c\frac{2\pi}{NK}]$ we have
$$|\frac{\partial^2 \bar{\theta}_{\lambda_{0},K,N}(x,t)}{\partial \alpha^2}|\geq \lambda_{0}\frac{|f_{2}(N^{1/2}(r-1)+1)|ln(K)|sin(-\lambda_{0} C_{0}t)|}{2}.$$

Therefore by using that $f(r)=1$ if $r\in(3/4,5/4)$ and defining 
$$B=\cup_{j\in \mathds{N}} \Big[j\frac{2\pi}{N}+\lambda_{0} t \frac{v_{\alpha}( f_{1})}{r}-c\frac{2\pi}{NK},j\frac{2\pi}{N}+\lambda_{0} t \frac{v_{\alpha}( f_{1})}{r}+c\frac{2\pi}{NK}\Big]$$
and $A=[1-\frac{N^{-1/2}}{4},1+\frac{N^{-1/2}}{4}]$, then

\begin{equation}\label{ABderivadagrande}
    \int_{A}\int_{B}\frac{1}{r^2}\Big(\frac{\partial^2\bar{\theta}_{\lambda_{0},K,N}}{\partial \alpha^2}\Big)^{2}d\alpha dr\geq \lambda^{2}_{0} \frac{ln(K)^2}{4(1+N^{-\frac12})^2}|A||B|sin(-\lambda_{0} C_{0}t)|^2,
\end{equation}
with $|A|,|B|$ the length of $A$ and $B$ respectively. We now consider $K$ big enough such that 
$$\lambda_{0}ln(K)|sin(-\lambda_{0} C_{0}t_{*})|\geq 16Mc_{0},$$
and thus, for $N$ big 

\begin{equation}
    \int_{A}\int_{B}\frac{1}{r^2}\Big(\frac{\partial^2\bar{\theta}_{\lambda_{0},K,N}}{\partial \alpha^2}\Big)^{2}d\alpha dr\geq 16M^2c^2_{0}|A||B|.
\end{equation}
Now, we can use lemmas \ref{errorhs} and \ref{errorl2} plus the interpolation inequality for sobolev spaces to obtain that, for $N$ big enough,  we have

$$||\Theta_{\lambda_{0},K,N}||_{H^2}\leq CN^{-a-1/4}$$
for some $a>0$ which can be computed explicitly but whose particular value is not relevant for this proof. With this we have that the solution $\theta_{\lambda_{0},K,N}$ satisfies that, at $t=t_{*}$

\begin{align*}
    &\Big(\int_{A}\int_{B}\frac{1}{r^4}\Big(\frac{\partial^2\theta_{\lambda_{0},K,N}}{\partial \alpha^2}\Big)^2d\alpha dr\Big)^{1/2}\\
    &=||\frac{1}{r^2}\frac{\partial^2\theta_{\lambda_{0},K,N}}{\partial \alpha^2}1_{A\times B}||_{L^2}\\
    &\geq ||\frac{1}{r^2}\frac{\partial^2\bar{\theta}_{\lambda_{0},K,N}}{\partial \alpha^2}1_{A\times B}||_{L^2}-||\frac{1}{r^2}\frac{\partial^2\Theta_{\lambda_{0},K,N}}{\partial \alpha^2}1_{A\times B}||_{L^2}\\
    &\geq 4Mc_{0}|A|^{1/2}|B|^{1/2}- Ct_{*}N^{-a-1/4}\\
\end{align*}

where we used that there is a constant $C$ such that if $S\subset \{\frac12\leq|x|\leq \frac32\}$ then
\begin{equation}\label{h2equivalpha}
    ||\frac{1}{r^2}\frac{\partial^2 g}{\partial \alpha^2}1_{S}||_{L^2}\leq C||g1_{S}||_{H^{2}}.
\end{equation}
But $|A||B|\geq CN^{-1/2}$, so, taking $N$ big enough we get

\begin{align*}
    &\Big(\int_{A}\int_{B}\frac{1}{r^4}\Big(\frac{\partial^2\theta_{\lambda_{0},K,N}}{\partial \alpha^2}\Big)^2d\alpha dr\Big)^{1/2}\\
    &\geq 3Mc_{0}|A|^{1/2}|B|^{1/2}.\\
\end{align*}

But, if $S\subset \{\frac12\leq|x|\leq \frac32\}$ then
\begin{equation}\label{c2equivalpha}
    \text{sup}_{x\in S}|\frac{1}{r^2}\frac{\partial^2 g}{\partial \alpha^2}|\leq 2||g1_{S}||_{C^{2}}\leq 2||g||_{C^{2}},
\end{equation}
so
\begin{align*}
    &\Big(\int_{A}\int_{B}\frac{1}{r^4}\Big(\frac{\partial^2\theta_{\lambda_{0},K,N}}{\partial \alpha^2}\Big)^2d\alpha dr\Big)^{1/2}\\
    &\leq 2|A|^{1/2}|B|^{1/2}||\theta_{\lambda_{0},K,N}||_{C^2},\\
\end{align*}
and thus

$$||\theta_{\lambda_{0},K,N}||_{C^2}\geq \frac{3Mc_{0}}{2}.$$

\end{proof}

\subsection{Non existence in $C^{k}$}
Now we can prove the last result of this section.

\begin{theorem}\label{nonexc2}
Given $c_{0}>0$, there are initial conditions $\theta_{0}\in H^{2+1/8}\cap C^2$ for the SQG equation (\ref{SQG}) such that $||\theta_{0}||_{C^2}\leq c_{0}$ and the only solution $\theta(x,t)\in H^{2+\frac18}$ with $\theta(x,0)=\theta_{0}(x)$ satisfies  that there exists a $t_{*}>0$ with $||\theta(x,t)||_{C^2}=\infty$ for all $t$ in the interval $(0,t_{*})$.
\end{theorem}

\begin{remark}\label{remarkunic2}
We can actually prove that, for the initial conditions $\theta_{0}(x)$ obtained in theorem \ref{nonexc2}, there is no solution in $L^{\infty}_{t}L^{2}_{x}$ such that $\theta(x,t)\in C^{2}$ for $t$ in some small time interval (even if we allow $\text{ess-sup}_{t\in[0,\epsilon]}||\theta(x,t)||_{C^{2}}=\infty$), since, if we call $\theta_{1}(x,t)$ the solution found in theorem \ref{nonexc2} and $\theta_{2}(x,t)$ the new solution belonging pointwise in time to $C^{2}$ for a small time interval, we can obtain the bound

$$\frac{ \partial ||\theta_{2}(x,t)-\theta_{1}(x,t)||_{L^{2}}}{\partial t}\leq C ||\theta_{2}(x,t)-\theta_{1}(x,t)||_{L^{2}}$$
which implies that $||\theta_{2}(x,t)-\theta_{1}(x,t)||_{L^{2}}=0$. 
\end{remark}

\begin{remark}
The value of $t_{*}$ can be made arbitrarily big if wanted with very small adjustments on the proof, but for simplicity we provide the proof without worrying about the specific value of $t_{*}$.
\end{remark}

\begin{proof}

We consider a family of pseudo-solutions to the SQG equation 
$$\bar{\theta}_{n}(x,t)=\bar{\theta}_{\lambda_{n},K_{n},N_{n}}(x,t)$$ for $n\in \mathds{N}$, with $\bar{\theta}_{\lambda_{n},K_{n},N_{n}}$ defined as in (\ref{pseudosolution}). Although $\bar{\theta}_{n}$ depends on the choice of $\lambda_{n}$, $K_{n}$ and $N_{n}$, we do not write the dependence explicitly to get a more compact notation. We start by fixing $\lambda_{n}$ satisfying 
 
$$\lambda_{n}\leq 2^{-n},$$
and such that $||\bar{\theta}_{n}(x,0)||_{C^2}\leq c_{0}$ independently of the choice of $K_{n}$ and $N_{n}.$

Note that this already tells us that for any fixed arbitrary $T$, if $0\leq t\leq T$ then  
$$||\bar{\theta}_{n}(x,t)||_{H^{2+1/8}}\leq C2^{-n}(\frac{K_{n}}{N_{n}^{1/8}}+1)$$
with $C$ depending on $T$. We will only consider $N_{n}^{1/8}\geq K_{n}$, so that $||\bar{\theta}_{n}(x,t)||_{H^{2+1/8}}\leq C2^{-n}$.
We fix now $K_{n}$ so that $\lambda^2_{n}ln(K_{n})\geq 16n$.  Note that then,  as seen in the proof of theorem \ref{strongillpos}, we have that there is a set $S_{n}=S_{\lambda_{n},K_{n},N_{n},t,}$ (see (\ref{ABderivadagrande})) with measure $|S_{n}|\geq \frac{c}{K_{n}N_{n}^{1/2}}>0$ such that the function $\bar{\theta}_{n}(x,t)$ fulfils that
\begin{equation}\label{zonaderivadagrande}
    ||\frac{1}{r^2}\frac{\partial^2 \bar{\theta}_{n}(x,t)}{\partial \alpha^2}1_{S_{n}}||_{L^2}\geq 4n\frac{|S_{n}|^{1/2}|sin(\lambda_{n}C_{0}t)|}{\lambda_{n}}.
\end{equation}

Let us consider now the initial conditions

$$\theta((\lambda_{n})_{n\in\mathds{N}},(K_{n})_{n\in\mathds{N}},(N_{n})_{n\in\mathds{N}},(R_{n})_{n\in\mathds{N}})=\sum_{n\in\mathds{N}}T_{R_{n}}(\bar{\theta}_{n}(x,0))$$
with $T_{R}(f(x_{1},x_{2}))=f(x_{1}+R,x_{2})$, with $R_{n}$ yet to be fixed. We will refer to these initial conditions simply as $\theta(x,0)$ and to the unique $H^{2+\frac18}$ solution to the SQG equation (\ref{SQG}) with initial conditions $\theta(x,0)$, as $\theta(x,t)$ for a more compact notation, keeping in mind that the function depends on multiple parameters. Since $||\bar{\theta}_{n}(x,0)||_{H^{2+1/8}}\leq C2^{-n}$ we have that  $|| \theta(x,0)||_{H^{2+1/8}}\leq C$, and thus we can use the a priori bounds to assure the existence of $\theta(x,t)$ for some time interval $[0,t_{ex}]$ and also $||\theta(x,t)||_{H^{2+1/8}}\leq C$ for some big $C$ for $t\in[0,\frac{t_{ex}}{2}]$. This also tells us that, in particular, $||v_{j}(\theta)||_{L^{\infty}}\leq v_{max}$ for some big constant $v_{max}$ for $t\in[0,\frac{t_{ex}}{2}]$ and $j=1,2$.

We restrict ourselves now to study the interval $t\in[0,t_{crit}]$ with
$$t_{crit}=min(\frac{t_{ex}}{2},\frac{\pi}{sup_{n}(\lambda_{n})C_{0}2}).$$

By construction, $\bar{\theta}_{n}(x,0)$ is contained in a ball of a certain radius $D$. Then, if we consider $R_{n}=R_{n-1}+2D+4v_{max}t_{crit}+D_{n}+D_{n-1}$ with $D_{n},D_{n-1}>0$, we have that 
$$d(\text{supp}(1_{B_{D+2v_{max}t_{crit}}(-R_{n},0)}\theta(x,t)),\text{supp}(\theta(x,t)-1_{B_{D+2v_{max}t_{crit}}(-R_{n},0)}\theta(x,t)))> D_{n}$$
and 

$$\tilde{\theta}_{n}(x,t):=\theta(x,t) 1_{B_{D+2v_{max}t_{crit}}(-R_{n},0)}$$
 is a pseudo-solution fulfilling

$$\frac{\partial \tilde{\theta}_{n}}{\partial t}+v_{1}(\tilde{\theta}_{n})\frac{\partial \tilde{\theta}_{n}}{\partial x_{1}}+ v_{2}(\tilde{\theta}_{n})\frac{\partial \tilde{\theta}_{n}}{\partial x_{2}}+\tilde{F}_{n}=0$$
$$v_{1}(\tilde{\theta}_{n})=-\frac{\partial}{\partial x_{2}}\Lambda^{-1} \tilde{\theta}_{n}=-\mathcal{R}_{2}\theta$$
$$v_{2}(\tilde{\theta}_{n})=\frac{\partial}{\partial x_{1}}\Lambda^{-1} \tilde{\theta}_{n}=\mathcal{R}_{1}\theta$$
$$\tilde{F}_{n}=v_{1}(\theta-\tilde{\theta}_{n})\frac{\partial \tilde{\theta}_{n}}{\partial x_{1}}+ v_{2}(\theta-\tilde{\theta}_{n})\frac{\partial \tilde{\theta}_{n}}{\partial x_{2}}$$
$$\tilde{\theta}_{n}(x,0)=\theta(x,0) 1_{B_{D+2v_{max}t_{crit}}(-R_{n},0)}.$$

If we now define $\Theta_{n}:=\tilde{\theta}_{n}-T_{R_{n}}(\bar{\theta}_{n}) $ we get 

\begin{align}
    &\frac{\partial \Theta_{n}}{\partial t} v_{1}(\Theta_{n})\frac{\partial \Theta_{n}}{\partial x_{1}}+ v_{2}(\Theta_{n})\frac{\partial \Theta_{n}}{\partial x_{2}}\nonumber\\
    &+v_{1}(\Theta_{n})\frac{\partial T_{R_{n}}( \bar{\theta}_{n})}{\partial x_{1}}+ v_{2}(\Theta_{\lambda,K,N})\frac{\partial T_{R_{n}}(\bar{\theta}_{n})}{\partial x_{2}}\\
    &+ v_{1}(\bar{\theta}_{n})\frac{\partial \Theta_{n}}{\partial x_{1}}+ v_{2}(T_{R_{n}}(\bar{\theta}_{n}))\frac{\partial \Theta_{n}}{\partial x_{2}}-T_{R_{n}}(F_{\lambda_{n},K_{n},N_{n}}(x,t))+\tilde{F}_{n}=0,\nonumber\\\nonumber
\end{align}
with $F_{\lambda_{n},K_{n},N_{n}}$ the source term of our pseudo-solution $\bar{\theta}_{n}=\bar{\theta}_{\lambda_{n},K_{n},N_{n}}$ and therefore satisfying the bounds given by lemmas \ref{errorl2} and \ref{errorhs},  
$$||F_{\lambda_{n},K_{n},N_{n}}||_{L^{2}}\leq \frac{C}{N_{n}^{2+3/4}}$$
and 
$$||F_{\lambda_{n},K_{n},N_{n}}||_{H^{2+1/4}}\leq \frac{C}{N_{n}^{1/8}}.$$

It is easy to prove that
$$||v_{i}(\theta-\tilde{\theta}_{n})1_{\text{supp}(\tilde{\theta}_{n})}||_{L^{\infty}}\leq \frac{C}{(D_{n})^2}$$ 
and in fact
\begin{equation}\label{cotasckv}
    ||v_{i}(\theta-\tilde{\theta}_{n})1_{\text{supp}(\tilde{\theta}_{n})}||_{C^{k}}\leq\frac{C}{(D_{n})^2}
\end{equation}
since  
$$d(supp(\tilde{\theta}_{n}),supp(\tilde{\theta}-\tilde{\theta}_{n}))\geq D_{n}.$$

Taking, for example, $D_{n}=N_{n}^{\frac{2+3/4}{2}}$ to obtain that $||\tilde{F}_{n}||_{L^{2}}\leq \frac{C}{N_{n}^{2+\frac34}}$
 we can argue as in lemma \ref{errorl2}  to get that 
$$||\Theta_{n}||_{L^{2}}\leq \frac{Ct}{N_{n}^{2+3/4}}$$
 for all $t\in[0,t_{crit}]$.
We can also estimate  $||\Theta_{n}||_{H^{2+1/8}}$ as in lemma \ref{errorhs}, being the only difference that now we have the extra term $\tilde{F}_ {n}$. Therefore, it is enough to obtain bounds for

$$\int_{\mathds{R}^2}D^{s}(\Theta_{n})(D^{s}(v_{1}(\theta-\tilde{\theta}_{n})\frac{\partial \tilde{\theta}_{n}}{\partial x_{1}}+ v_{2}(\theta-\tilde{\theta}_{n})\frac{\partial \tilde{\theta}_{n}}{\partial x_{2}}))dx_{1}dx_{2}$$
with $s=2+1/8$.

Using lemma \ref{katoponce} in the same way as we did in lemma \ref{errorhs}, we can decompose this integral in several terms that are easy to bound using (\ref{cotasckv}) plus the term

$$\int_{\mathds{R}^2}D^{s}(\Theta_{n})(v_{1}(\theta-\tilde{\theta}_{n})D^{s}\frac{\partial \tilde{\theta}_{n}}{\partial x_{1}}+ v_{2}(\theta-\tilde{\theta}_{n})D^{s}\frac{\partial \tilde{\theta}_{n}}{\partial x_{2}})dx_{1}dx_{2}$$
which is, in principle, too irregular to be bounded. However, using incompressibility and $\Theta_{n}=\tilde{\theta}_{n}-T_{R_{n}}(\bar{\theta}_{n})$ we get

\begin{align*}
    &|\int_{\mathds{R}^2}D^{s}(\Theta_{n})(v_{1}(\theta-\tilde{\theta}_{n})D^{s}\frac{\partial \tilde{\theta}_{n}}{\partial x_{1}}+ v_{2}(\theta-\tilde{\theta}_{n})D^{s}\frac{\partial \tilde{\theta}_{n}}{\partial x_{2}})dx_{1}dx_{2}|\\
    &=|\int_{\mathds{R}^2}D^{s}(T_{R_{n}}(\bar{\theta}_{n}))(v_{1}(\theta-\tilde{\theta}_{n})D^{s}\frac{\partial \tilde{\theta}_{n}}{\partial x_{1}}+ v_{2}(\theta-\tilde{\theta}_{n})D^{s}\frac{\partial \tilde{\theta}_{n}}{\partial x_{2}})dx_{1}dx_{2}|\\
    &=|\int_{\mathds{R}^2}D^{s}(\tilde{\theta}_{n})(v_{1}(\theta-\tilde{\theta}_{n})D^{s}\frac{\partial T_{R_{n}}(\bar{\theta}_{n})}{\partial x_{1}}+ v_{2}(\theta-\tilde{\theta}_{n})D^{s}\frac{\partial T_{R_{n}}(\bar{\theta}_{n})}{\partial x_{2}})dx_{1}dx_{2}|\\
    &\leq ||D^{s}\tilde{\theta}_{n}||_{L^{2}}\frac{C}{N_{n}^{2+\frac34}}N_{n}^2\leq \frac{C}{N_{n}^{\frac34}}.\\
    \end{align*}
Therefore, as in lemma \ref{errorhs}, we get

$$||\Theta_{n}||_{H^{2+1/8}}\leq \frac{Ct}{N_{n}^{1/8}}.$$

This combined with the $L^2$ norm and using the interpolation inequality for Sobolev spaces gives us

$$||\Theta_{n}||_{H^2}\leq \frac{Ct}{N_{n}^{11/34}}=\frac{Ct}{N_{n}^{1/4+a}}.$$
with $a>0$, for all $t\in[0,t_{crit}]$.

But, this means that, if we consider the polar coordinates around the point $(-R_{n},0)$, which we will call $(r_{R_{n}},\alpha_{R_{n}})$, and using (\ref{h2equivalpha})
\begin{align*}
    &||\frac{1}{r_{R_{n}}^2}\frac{\partial^2 \theta(x,t)}{\partial \alpha_{R_{n}}^2}T_{R_{n}}(1_{S_{n}})||_{L^2}\\
    &\geq ||T_{R_{n}}(\frac{1}{r^2}\frac{\partial^2 \bar{\theta}(x,t)}{\partial \alpha^2}1_{S_{n}})||_{L^2}-||\tilde{\theta}_{n}-T_{R_{n}}(\bar{\theta}_{n}(x,t))||_{H^{2}}\\
    &\geq 4n|S_{n}|^{1/2}\frac{|sin(\lambda_{n}C_{0}t)|}{\lambda_{n}}-\frac{Ct}{N_{n}^{1/4+a}}\\
\end{align*}
but, using $C_{0}t\lambda_{n}\leq\frac{\pi}{2} $, $|S_{n}|\geq CK^{-1}_{n}N_{n}^{-1/2}$ and taking $N_{n}$ big enough we get

$$||T_{R_{n}}(1_{S_{n}})\frac{1}{r_{R_{n}}^2}\frac{\partial^2 \theta(x,t)}{\partial \alpha^2}||_{L^2}\geq cnt|S_{n}|^{1/2} $$
for some small constant $c$.

But then

\begin{align*}
    &||T_{R_{n}}(1_{S_{n}})\frac{1}{r_{R_{n}}^2}\frac{\partial^2 \tilde{\theta}_{n}(x,t)}{\partial \alpha_{R_{n}}^2}||_{L^2}\\
    &\leq || T_{R_{n}}\tilde{\theta}_{n}(x,t)||_{C^{2}} |S_{n}|^{1/2}\\
\end{align*}
and thus $||T_{R_{n}}(1_{S_{n}}) \theta(x,t)||_{C^2}\geq cnt$ and we are done since we can do this for every $n$.

\end{proof}

Both results in this section can be obtained in $C^{m}$ for $m\geq 2$, using exactly the same method. To do it we consider pseudo-solutions of the form

$$\lambda (f_{1}(r)+f_{2}(N^{1/2}(r-1)+1)\sum_{k=1}^{K}\frac{sin(Nk\alpha)}{N^{m}k^{m+1}}).$$

The proof follows exactly the same method, only this time we have that the associated source terms $F_{\lambda,K,N}$ of these pseudo-solutions fulfil $||F_{\lambda,K,N}||_{L^{2}}\leq \frac{C}{N^{m+3/4}}$, $||F_{\lambda,K,N}||_{H^{k}}\leq C N^{k-m-1/4}$, which gives us, by taking $k$ big an using the interpolation inequality that
$$||F_{\lambda,K,N}||_{H^{m+\frac14}}\leq C N^{-\frac{1}{2}+\delta}$$
for $\delta>0$ arbitrary.

Note also that analogous expressions as (\ref{h2equivalpha}) and (\ref{c2equivalpha}) exists for higher order derivatives in $\alpha$, albeit with different constants.

\section{Strong ill-posedness and non existence in supercritical Sobolev spaces}

\subsection{Pseudo-solutions for $H^{s}$}

The proof for ill-posedness in supercritical Sobolev spaces follows a very similar strategy as before. We find an appropriate pseudo-solution with the desired properties, we find bounds for the source term and then  we obtain bounds for the difference between the real solution and the pseudo-solution. This time, we will consider pseudo-solutions of the form

$$\bar{\theta}(r,\alpha,t)= f_{1}(r)+ f_{2}(r) \frac{sin(N\alpha-Nt\frac{v_{\alpha}(f_{1}(r)}{r}))r_{0}^{\beta}}{N^{\beta}}$$

with $f_{1}$, $f_{2}$ compactly supported $C^{\infty}$ functions,  $r_{0}>0$ and $v_{\alpha}(f_{1}(r))$ is the angular velocity generated by the function $f_{1}(r)$. 

The choice of $f_{1}$, $f_{2}$ and $r_{0}$ will depend on the specific behaviour we want our pseudo-solutions to have. Before we start to specify how we choose them and how we will label the pseudo-solutions, we need the following technical lemma.

\begin{lemma}\label{f1sobolev}
For any $\beta\in (\frac32,2)$ and $K,c>0$, there exists a $C^{\infty}$ radial function $f_{1}(r):\mathds{R}_{+}\times[0,2\pi]\rightarrow \mathds{R}$, with support in some $[a_{1},a_{2}]\times[0,2\pi]$, $0<a_{1}<a_{2}$  depending on $K$, $c$ and $\beta$ such that $||f_{1}(r)||_{H^{\beta}}\leq c$, and  $|\frac{\partial \frac{v_{\alpha}(f_{1}(.))(r)}{r}}{\partial r}(r=\frac{a_{1}}{2})|\geq \frac{2K}{a_{1}}$.
\end{lemma}
\begin{proof}
By lemma \ref{seleccionvelocidadnagular}, we can find a $C^{\infty}$ function $g(r):\mathds{R}_{+}\times[0,2\pi]\rightarrow \mathds{R}$  with support in $r\in[2,M]$ such that $\frac{\partial\frac{ v_{\alpha}(g(.))(r)}{r}}{\partial r}(r=1)=1$. If we consider now the functions

$$g_{\lambda_{1},\lambda_{2}}(r):=\frac{g(\lambda_{1} r)}{\lambda_{2}\lambda_{1}^{\beta-1}},\quad \lambda_{1},\lambda_{2} > 1$$

we have (for example using the interpolation inequalities for Sobolev spaces) that
\begin{equation}\label{clambda2}
   ||g_{\lambda_{1},\lambda_{2}}(r)||_{H^{\beta}}\leq \frac{C}{\lambda_{2}} 
\end{equation}

with $C$ depending on $||g(r)||_{H^{2}}.$

Furthermore, $v_{\alpha}(f(\lambda \cdot ))(\frac{r}{\lambda})=v_{\alpha}(f(\cdot ))(r)$, $ \frac{\partial v_{\alpha}(f(\lambda \cdot ))}{\partial r}(\frac{r}{\lambda})=\lambda\frac{\partial v_{\alpha}(f(\cdot ))}{\partial r}(r)$, so

$$\frac{\partial\frac{ v_{\alpha}(g_{\lambda_{1},\lambda_{2}}(.))(r)}{r}}{\partial r}(r=\frac{1}{\lambda_{1}})=\frac{\lambda_{1}^{3-\beta}}{\lambda_{2}}=\frac{\lambda_{1}^{2-\beta}}{\frac{1}{\lambda_{1}}\lambda_{2}}.$$

Therefore it is enough to take $g_{\lambda_{1},\lambda_{2}}$ with  $\lambda_{2}$ big enough so that $\frac{C}{\lambda_{2}}\leq c$ ($C$ the constant in (\ref{clambda2})) and then $\lambda_{1}$ big enough so that $\frac{\lambda_{1}^{2-\beta}}{\lambda_{2}}\geq K$ and $g_{\lambda_{1},\lambda_{2}}$ with $a_{1}=\frac{2}{\lambda_{1}}$, $a_{2}=\frac{M}{\lambda_{1}}$ will have all the properties desired.
\end{proof}

From now on we consider $\beta$ a fixed value in the interval $ (\frac{3}{2},2)$. The family of pseudo-solutions we consider to obtain ill-posedness in $H^{\beta}$ is, for $N\in\mathds{N}$ 

\begin{equation}\label{pseudoNcK}
   \bar{\theta}_{N,c,K}(r,\alpha,t)= f_{1,c,K}(r)+ f_{2,c,K}(r)r_{c,K}^{\beta} \frac{sin(N\alpha-Nt\frac{v_{\alpha}(f_{1}(r))}{r})}{N^{\beta}} 
\end{equation}

with $f_{1,c,K}$ the function given by lemma \ref{f1sobolev} for the specific values of $c$ and $K$ considered and  $r_{c,K}=\frac{a_{1}}{2}$ given by the lemma. By continuity, we have that there exists an interval $[r_{c,K}-\epsilon,r_{c,K}+\epsilon]$ such that if $\bar{r}\in[r_{c,K}-\epsilon,r_{c,K}+\epsilon]$ then

\begin{equation}\label{velocidadrapida}
    \frac{\partial \frac{v_{\alpha}(f_{1,c,K}(.))(r)}{r}}{\partial r}(r=\bar{r})\geq \frac{K}{2\bar{r}}.
\end{equation}

We take $f_{2,c,K}$ to be a $C^{\infty}$ function with support in $[r_{c,K}-\epsilon,r_{c,K}+\epsilon]\cap [\frac{r_{c,K}}{2},\frac{3r_{c,K}}{2}]$ and fulfilling $||f_{2,c,K}||_{L^{2}}=c.$ 

These pseudo-solutions fulfil the evolution equation

$$\frac{\partial \bar{\theta}_{N,c,K}}{\partial t}+ \frac{v_{\alpha}(f_{1,c,K}(\cdot))}{r} \frac{\partial \bar{\theta}_{N,c,K}}{\partial \alpha}=0$$

and therefore they are pseudo-solutions with source term
\begin{align}\label{sourceNcK}
    &F_{N,c,K}\\
    &=-(\frac{v_{\alpha}(\bar{\theta}_{N,c,K}(\cdot)-f_{1,c,K}(\cdot))}{r} \frac{\partial \bar{\theta}_{N,c,K}}{\partial \alpha}+v_{r}(\bar{\theta}_{N,c,K}(\cdot)) \frac{\partial \bar{\theta}_{N,c,K}}{\partial r})\nonumber\\
    &=-(\frac{v_{\alpha}(\bar{\theta}_{N,c,K}(\cdot)-f_{1,c,K}(\cdot))}{r} \frac{\partial \bar{\theta}_{N,c,K}}{\partial \alpha}+v_{r}(\bar{\theta}_{N,c,K}(\cdot)-f_{1,c,K}(\cdot)) \frac{\partial \bar{\theta}_{N,c,K}}{\partial r}).\nonumber\\\nonumber
\end{align}

Next we need to obtain bounds for our source term. To do this, we start with a lemma analogous to lemma \ref{decaimiento}:

\begin{lemma}\label{decaimiento2}
Given a $L^{\infty}$ function $g(.):\mathds{R}\rightarrow \mathds{R}$ with support in the interval $(a,b)$ then if we define $g_{N}$ as

$$g_{N}(r,\alpha):=sin(N\alpha+\alpha_{0})\tilde{g}_{N}(r)$$

with $N$ a natural number, then there is a constant $C$  depending on $(a,b)$  such that if $r>b $, then

$$|v_{r}(g_{N})|(r,\alpha)\leq \frac{C||g_{N}||_{L^{\infty}}}{N|r-b|^2}.$$

Furthermore, we have that if $||\tilde{g}_{N}||_{C^{i}}\leq M N^{i}$ for $i=0,1,...,m$, then

$$|\frac{\partial^{m}v_{r}(g_{N})}{\partial x_{1}^{m-i}\partial x_{2}^{i}}|(r,\alpha)\leq \frac{CMN^{m-1}}{|r-b|^{2}},$$

with $C$  depending on $(a,b)$ and $m$.

\end{lemma}

\begin{proof}

The proof for the decay of the velocity it is analogous to that of lemma \ref{decaimiento}. As for the higher derivatives, using that

$$v_{r}(w)=cos(\alpha(x)) v_{1}(w)+sin(\alpha(x))v_{2}(w),$$

one can obtain that

\begin{align*}
    &|\frac{\partial^{m}v_{r}(g_{N})}{\partial x_{1}^{m-i}\partial x_{2}^{i}}(r,\alpha)|\\
    &\leq |v_{r}(\frac{\partial^{m}g_{N}}{\partial x_{1}^{m-i}\partial x_{2}^{i}}(r,\alpha))(r,\alpha)|\\
    &+C\sum_{i=0}^{m-1}\sum_{j=0}^{i}(\frac{\partial^{i}v_{1}(g_{N})}{\partial^{j}x_{1}\partial^{i-j} x_{2}})(r,\alpha)\\
    &+C\sum_{i=0}^{m-1}\sum_{j=0}^{i}(\frac{\partial^{i}v_{1}(g_{N})}{\partial^{j}x_{1}\partial^{i-j} x_{2}})(r,\alpha)\\
\end{align*}

with $C$ depending on $m,a$ and $b$, and using the decay for $v_{r},$ and
$$|v_{1}(w)(x)|\leq C\frac{||w||_{L^{1}}}{|d(x,supp(w))|^2}$$
$$|v_{2}(w)(x)|\leq C\frac{||w||_{L^{1}}}{|d(x,supp(w))|^2}$$
we are done.

\end{proof}

With this, we are now ready to obtain the bounds for our source term.

\begin{lemma}
For $t\in[0,T]$ and a pseudo-solution $\bar{\theta}_{N,c,K}$ as in (\ref{pseudoNcK}) then the source term $F_{N,c,K}(x,t)$ as in (\ref{sourceNcK}) satisfies

$$||F_{N,c,K}(x,t)||_{L^{2}}\leq CN^{-(2\beta-1)}$$
with $C$ depending on  $c$, $K$ and $T$.
\end{lemma}

\begin{proof}
In order to obtain the desired estimate we divide the source term in several parts. First we have
\begin{align*}
    &||\frac{v_{\alpha}(\bar{\theta}_{N,c,K}(\cdot)-f_{1,c,K}(\cdot))}{r} \frac{\partial \bar{\theta}_{N,c,K}}{\partial \alpha}||_{L^{2}}\\
    &\leq C||v_{\alpha}(\bar{\theta}_{N,c,K}(\cdot)-f_{1,c,K}(\cdot))||_{L^{2}}||\frac{\partial \bar{\theta}_{N,c,K}}{\partial \alpha}||_{L^{\infty}}\leq \frac{C}{N^{2\beta-1}}\\
\end{align*}
and analogously
\begin{align*}
    &||v_{r}(\bar{\theta}_{N,c,K}(\cdot)-f_{1,c,K}(\cdot)) \frac{\partial (\bar{\theta}_{N,c,K}-f_{1,c,K}(r))}{\partial r}||_{L^{2}}\\
    &\leq  ||v_{r}(\bar{\theta}_{N,c,K}(\cdot)-f_{1,c,K}(\cdot))||_{L^{2}}||\frac{\partial (\bar{\theta}_{N,c,K}-f_{1,c,K}(r))}{\partial r}||_{L^{\infty}}\leq \frac{C}{N^{2\beta-1}}.\\
\end{align*}

Finally, by using that $supp(f_{1,c,K})\in[2r_{c,k},a_{2})$ (see lemma \ref{f1sobolev} and the definition of the pseudo-solution ), $supp(f_{2,c,K})\in[\frac{r_{c,k}}{2},\frac{3r_{c,k}}{2}]$ and together with  lemma \ref{decaimiento2} we have

\begin{align}\label{decaerck}
    &||v_{r}(\bar{\theta}_{N,c,K}(\cdot)-f_{1,c,K}(\cdot)) \frac{\partial f_{1,c,K}(r)}{\partial r}||_{L^{2}}\nonumber\\
    &\leq \Big(\int_{2r_{c,k}}^{a_{2}}\frac{C}{N^{2+2\beta}(r-\frac{3r_{c,K}}{2})^4}rdr\Big)^{1/2}\leq \frac{C}{N^{1+\beta}}.\\\nonumber
\end{align}

Combining all three bounds we obtain the desired result.
\end{proof}

\begin{lemma}
For $t\in[0,T]$ and a pseudo-solution $\bar{\theta}_{N,c,K}$ as in (\ref{pseudoNcK}) then the source term $F_{N,c,K}(x,t)$ as in (\ref{sourceNcK}) satisfies, for $k\in\mathds{N}$

$$||F_{N,c,K}(x,t)||_{H^{k}}\leq \frac{C}{N^{2\beta-1-k}}$$
with $C$ depending on $k$, $c$, $K$ and $T$.
\end{lemma}

\begin{proof}
We separate the source term in three different parts:

1) Using the properties of the support of $\bar{\theta}_{N,c,K}$
\begin{align*}
    &||\frac{\partial \bar{\theta}_{N,c,K}}{\partial \alpha} \frac{v_{\alpha}(f_{1,c,K}-\bar{\theta}_{N,c,K})}{r}||_{H^{k}}\\
    &\leq C \sum_{i=0}^{k}||\frac{\partial \bar{\theta}_{N,c,K}}{\partial \alpha}||_{C^{i}}||v_{\alpha}(f_{1,c,K}-\bar{\theta}_{N,c,K})||_{H^{k-i}}\\
    &\leq \frac{C}{N^{2\beta-1-k}},\\
\end{align*}
2)
\begin{align*}
    &||\frac{\partial (\bar{\theta}_{N,c,K}-f_{1,c,K})}{\partial r}v_{r}(\bar{\theta}_{N,c,K})||_{H^{k}}\\
    &\leq C  \sum_{i=0}^{k}||\frac{\partial (\bar{\theta}_{N,c,K}-f_{1,c,K})}{\partial r}||_{C^{i}}||v_{r}(\bar{\theta}_{N,c,K})||_{H^{k-i}}\\
    &\leq \frac{C}{N^{2\beta-1-k}},\\
\end{align*}
3)
To bound $||\frac{\partial f_{1,c,K}}{\partial r}v_{r}(\bar{\theta}_{N,c,K})||_{H^{3}},$ we just apply lemma \ref{decaimiento2} as in (\ref{decaerck}) to obtain

\begin{align*}
    &||\frac{\partial f_{1,c,K}}{\partial r}v_{r}(\bar{\theta}_{N,c,K})||_{H^{k}}\\
    &\leq \sum_{i=0}^{s}\sum_{j=0}^{i}||\frac{\partial^{i} (\frac{f_{1,c,K}}{\partial r}v_{r}(\bar{\theta}_{N,c,K}))}{\partial^{j}x_{1}\partial^{i-j} x_{2}}||_{L^{2}}\leq \frac{C}{N^{\beta+1-k}}\leq  \frac{C}{N^{2\beta-1-k}}.\\
\end{align*}
\end{proof}
And applying the interpolation inequality for Sobolev spaces (with $L^{2}$ and for example $H^{3}$) we obtain the following corollary:
\begin{corollary}
For $t\in[0,T]$ and a pseudo-solution $\bar{\theta}_{N,c,K}$ as in (\ref{pseudoNcK}) then the source term $F_{N,c,K}(x,t)$ as in (\ref{sourceNcK}) satisfies

$$||F_{N,c,K}(x,t)||_{H^{\beta+\frac12}}\leq CN^{-(\beta-\frac32) }$$
with $C$ depending on  $c$, $K$, $T$.
\end{corollary}

Now, as in last section, we define $\theta_{N,c,K}(x,t)$ to be the unique $H^{\beta+\frac12}$ solution to (\ref{SQG}) with initial conditions $\theta_{N,c,K}(x,0)=\bar{\theta}_{N,c,K}(x,0)$, and we denote
\begin{equation}\label{evoluciondiferencia2}
    \Theta_{N,c,K}:=\theta_{N,c,K}-\bar{\theta}_{N,c,k}.
\end{equation}

The next step now is to find bounds for $\Theta_{N,c,K}$.

\begin{lemma}\label{errorl2sobolev}
Let $\Theta_{\lambda,K,N}$ defined as in (\ref{evoluciondiferencia2}), then if $\theta_{\lambda,K,N}$ exists for $t\in [0,T]$, we have that

$$||\Theta_{N,c,K}(x,t)||_{L^2}\leq \frac{Ct}{N^{(2\beta-1)}}$$

with $C$ depending on $\lambda$, $K$ and $T$.
\end{lemma}
\begin{proof}
As in the proof of lemma \ref{errorl2}, we obtain the equation

 \begin{align*}
    &\frac{\partial }{\partial t} \frac{||\Theta_{N,c,K}||_{L^2}^2}{2}\leq |\int_{\mathds{R}^2}\Theta_{N,c,K}\\
    &\Big( v_{1}(\Theta_{N,c,K})\frac{\partial \bar{\theta}_{N,c,K}}{\partial x_{1}}+ v_{2}(\Theta_{N,c,K})\frac{\partial \bar{\theta}_{N,c,K}}{\partial x_{2}}+F_{N,c,K}(x,t)\Big)dx|\\
    &\leq ||\Theta_{N,c,K}||_{L^2} \Big(||\Theta_{N,c,K}||_{L^2}||\bar{\theta}_{N,c,K}||_{C^1}+||F_{N,c,K}(x,t)||_{L^{2}}\Big).\\
\end{align*}

By using that $||F_{N,c,K}||_{L^{2}}\leq \frac{C}{N^{(2\beta-1)}}$, $||\bar{\theta}_{\lambda,K,N}||_{C^1}\leq C$ and integrating it follows

$$||\Theta_{N,c,K}||_{L^2}\leq \frac{C (e^{Ct}-1)}{N^{(2\beta-1)}}.$$

\end{proof}

Before obtaining the bounds for the higher order norms of $\Theta_{N,c,K}$ we need a couple of technical lemmas:

\begin{lemma}\label{laplacianoestructura}
Given a $C^{1}$ function $h(x):\mathds{R}^{2}\rightarrow \mathds{R}$ with $||h||_{L^{\infty}}\leq M$, $||h||_{C^{1}}\leq MN$ and $a\in(0,1)$, then there exists a constant $C$ depending on $a$ such that

$$||(-\Delta)^{a/2}(h(x))||_{L^{\infty}}\leq C M N^{a}.$$

\end{lemma}
\begin{proof}

Using the integral expression from the fractional Laplacian 

$$(-\Delta)^{a/2}(h(.))(x)=C\int_{\mathds{R}^2}\frac{(h(x)-h(z))}{|x-z|^{2+a}}dz$$

and dividing the integral in two parts depending on the value of $|x-z|$ we get

$$\int_{|x-z|\geq \frac{1}{N}}\frac{(h(x)-h(z))}{|x-z|^{2+a}}dz\leq CN^{a}||h||_{L^{\infty}}=CMN^{a}$$

$$\int_{|x-z|\leq \frac{1}{N}}\frac{(h(x)-h(z))}{|x-z|^{2+a}}dz\leq C N^{a-1}||h||_{C^{1}}=CMN^{a}$$
and we are done.

\end{proof}

\begin{lemma}\label{velocidadacotada}
Given a $C^{1}$ function $h(x):\mathds{R}^{2}\rightarrow \mathds{R}$ with $||h||_{L^{\infty}}\leq M$, $||h||_{C^{1}}\leq MN$ and with support in the set $[-R,R]^2$ for some $R$, we have that there exists a constant $C$ depending on $R$ such that for $i=1,2$

$$||v_{i}(h(x))||_{L^{\infty}}\leq CM log(N) .$$

Furthermore, if $||h||_{C^{n}}\leq M$,$||h||_{C^{n+1}}\leq MN$ for some natural number $n$ we also have that, for $i=1,2$, $k=0,2,...,n$

$$||\frac{\partial^{n} v_{i}(h(x))}{\partial^{n-k} x_{1} \partial^{k} x_{2}}||_{L^{\infty}}\leq CMlog(N) .$$
\end{lemma}

\begin{proof}
The proof  of the first part is the same as in lemma \ref{laplacianoestructura} but using the kernel for $v_{i}$ instead of the one for $(-\Delta)^{a/2}$. For the second part we just need to use that, for sufficiently regular functions we have that 
$$\frac{\partial v_{i}(h(x))}{\partial x_{j}}=v_{i}\Big(\frac{\partial h(x)}{\partial x_{j}}\Big).$$
\end{proof}

\begin{lemma}\label{errorhssobolev}
Let $\Theta_{N,c,K}$ defined as in (\ref{evoluciondiferencia2}), then we have that, for $N$ large,  $\theta_{N,c,K}$ exists for $t\in [0,T]$ and

$$||\Theta_{N,c,K}(x,t)||_{H^{\beta+\frac12}}\leq \frac{Ct}{N^{\beta-\frac32}}$$

with $C$ depending on $\lambda$, $K$ and $T$.
\end{lemma}

\begin{proof}
The proof is very similar to that of lemma \ref{errorhs}. We will prove the inequality for the time interval $[0,T^{*}]$ with $T^{*}$ the smallest time fulfilling  $\quad$ $||\Theta_{N,c,K}(x,t)||_{H^{\beta+1/2}}= log(N)N^{-(\beta-\frac32)}$, (we can just consider $t\in[0,T]$ directly if $T_{*}>T$ or if it does not exists) but note that this is enough since then we can take $N$ big enough so that $T^{*}\geq T$. Note also that, since we have local existence, obtaining this bound also allows us to ensure that we have existence for the times considered.

First we have that, for $s=\beta+\frac12$
\begin{align*}
    &\frac{\partial}{\partial t} \frac{||D^{s}\Theta_{N,c,K}||_{L^{2}}^2}{2}=- \int_{\mathds{R}^2} D^{s}\Theta_{N,c,K} \\
        &D^{s}\Big((v_{1}(\Theta_{N,c,K})+v_{1}(\bar{\theta}_{N,c,K}))\frac{\partial \Theta_{N,c,K}}{\partial x_{1}}+ (v_{2}(\Theta_{N,c,K})+v_{2}(\bar{\theta}_{N,c,K}))\frac{\partial \Theta_{N,c,K}}{\partial x_{2}}\\
    &+v_{1}(\Theta_{N,c,K})\frac{\partial \bar{\theta}_{N,c,K}}{\partial x_{1}}+ v_{2}(\Theta_{N,c,K})\frac{\partial \bar{\theta}_{N,c,K}}{\partial x_{2}}+F_{N,c,K}(x,t)\Big)dx.
\end{align*}

We start bounding
 
 $$\int_{\mathds{R}^2} D^{s}\Theta_{N,c,K}
D^{s}\Big(v_{1}(\bar{\theta}_{N,c,K}))\frac{\partial \Theta_{N,c,K}}{\partial x_{1}}+ v_{2}(\bar{\theta}_{N,c,K})\frac{\partial \Theta_{N,c,K}}{\partial x_{2}}\Big)dx.$$

 Applying lemma \ref{katoponce} with $s_{1}=s-1$, $s_{2}=1$, $f=v_{i}(\bar{\theta}_{N,c,K}))$, $g=\frac{\partial \Theta_{N,c,K}}{\partial x_{i}}$, $i=1,2$ we would get that

\begin{align*}
    &(D^{s}\Theta_{N,c,K},D^{s}(fg)-\sum_{|\mathbf{k}|\leq s_{1}}\frac{1}{\mathbf{k}!}\partial^{\mathbf{k}}f D^{s,\mathbf{k}}g-\sum_{|\mathbf{j}|\leq s_{2}}\frac{1}{\mathbf{j}!}\partial^{\mathbf{j}}g D^{s,\mathbf{j}}f)_{L^{2}}\\
    &\leq C ||D^{s}\Theta_{N,c,K}||_{L^{2}}||D^{s_{1}}f||_{BMO} ||D^{s_{2}}g||_{L^{2}}\\
    &\leq C||D^{s}\Theta_{N,c,K}||_{L^{2}} ||\Theta_{N,c,K}||_{H^s},\\
\end{align*}
where we used $||D^{s_{1}}v_{i}(\bar{\theta}_{N,c,K})||_{L^{\infty}}\leq C$.
Furthermore we have that

\begin{align*}
    &(D^{s}\Theta_{N,c,K}, D^{s}(\frac{\partial \Theta_{N,c,K}}{\partial x_{1}})v_{1}(\bar{\theta}_{N,c,K})+ D^{s}(\frac{\partial \Theta_{N,c,K}}{\partial x_{2}})v_{2}(\bar{\theta}_{N,c,K}))_{L^{2}}\\
    &=\frac{1}{2}\int_{\mathds{R}^2} \frac{\partial}{\partial x_{1}}(D^{s}\Theta_{N,c,K})^2
v_{1}(\bar{\theta}_{N,c,K})+ \frac{\partial}{\partial x_{2}}(D^{s}\Theta_{N,c,K})^2v_{2}(\bar{\theta}_{N,c,K})dx=0\\
\end{align*}

and, for $i=1,2$, using that the operators $D^{s,\mathbf{k}}$ are continuous from $H^{a}$ to $H^{a-s+\mathbf{k}}$,
\begin{align*}
    &|(D^{s}\Theta_{N,c,K},\sum_{|\mathbf{k}|=1}\frac{1}{\mathbf{k}!}\partial^{\mathbf{k}}v_{i}(\bar{\theta}_{N,c,K}) D^{s,\mathbf{k}}\frac{\partial \Theta_{\lambda,K,N}}{\partial x_{i}})_{L^{2}}|\\
    &\leq C ||D^{s}\Theta_{N,c,K}||_{L^{2}}||v_{i}(\bar{\theta}_{N,c,K})||_{C^1}||\Theta_{N,c,K}||_{H^{s}}\\
    &\leq C ||D^{s}\Theta_{N,c,K}||_{L^{2}}||\Theta_{N,c,K}||_{H^s}\\
\end{align*}

where we used $||v_{i}(\bar{\theta}_{_{N,c,K}}-f_{1,c,K})||_{C^1}\leq Clog(N)N^{\beta-1}$ (consequence of lemma \ref{velocidadacotada}) and $||v_{i}(f_{1,c,K})||_{C^1}\leq C$.

We also have
\begin{align*}
    &|(D^{s}\Theta_{\lambda,K,N},\sum_{|\mathbf{j}|=1}\frac{1}{\mathbf{j}!}\partial^{\mathbf{j}}\frac{\partial \Theta_{\lambda,K,N}}{\partial x_{i}} D^{s,\mathbf{j}}v_{i}(\bar{\theta}_{\lambda,K,N}))_{L^{2}}|\\
    &\leq C\sum_{|\mathbf{j}|=1} ||D^{s}\Theta_{\lambda,K,N}||_{L^{2}}||\frac{1}{\mathbf{j}!}\partial^{\mathbf{j}}\frac{\partial \Theta_{\lambda,K,N}}{\partial x_{i}}||_{L^{2}}||D^{s,\mathbf{j}}v_{i}(\bar{\theta}_{\lambda,K,N})||_{L^{\infty}}\\
    &\leq C \sum_{|\mathbf{j}|=1}||D^{s}\Theta_{\lambda,K,N}||_{L^{2}}|| \Theta_{\lambda,K,N}||_{H^{s}}||D^{s-2}\partial^{\mathbf{j}} v_{i}(\bar{\theta}_{\lambda,K,N})||_{L^{\infty}},\\
    & \leq C ||D^{s}\Theta_{\lambda,K,N}||_{L^{2}}|| \Theta_{\lambda,K,N}||_{H^{s}}(\frac{N N^{s-2}log(N)}{N^{\beta}}+C)\\
    &\leq C ||D^{s}\Theta_{\lambda,K,N}||_{L^{2}}|| \Theta_{\lambda,K,N}||_{H^{s}},\\
\end{align*}

where we used lemmas \ref{laplacianoestructura} and \ref{velocidadacotada}, the expression for $D^{s,\mathbf{j}}$ and the bounds for the derivatives of $\bar{\theta}_{N,c,K}$.

The last part to bound from the term with $v_{i}(\bar{\theta}_{N,c,K})$ is, for $i=1,2$
\begin{align*}
    &|(D^{s}\Theta_{N,c,K}, \frac{\partial \Theta_{N,c,K}}{\partial x_{i}} D^{s}v_{i}(\bar{\theta}_{N,c,K}))_{L^{2}}|\\
    &\leq C ||D^{s}\Theta_{N,c,K}||_{L^{2}}||\Theta_{N,c,K}||_{H^1}||D^{s-2}v_{i}(\Delta \bar{\theta}_{N,c,K})||_{L^{\infty}}.\\
    &\leq C ||D^{s}\Theta_{N,c,K}||_{L^{2}}N^{-(\beta-\frac12)}\frac{CN^{s-2}log(N)N^2}{N^{\beta}}\\
    &\leq C ||D^{s}\Theta_{N,c,K}||_{L^{2}}N^{-\beta+1}log(N)\leq C ||D^{s}\Theta_{N,c,K}||_{L^{2}} N^{-\frac12},\\
\end{align*}
where we used that, for the times considered, using lemma \ref{errorl2sobolev} and the interpolation inequality we have $||\Theta_{N,c,K}||_{H^1}\leq CN^{-(\beta-\frac12)}$ (the bound is actually better, but this is enough).






The rest of the terms not depending on $F_{N,c,K}$ are bounded in a similar fashion, and using $||\bar{\theta}_{\lambda,K,N}||_{H^{s}}\leq C$, $||F_{N,c,K}||_{H^{s}}\leq C N^{-(\beta-\frac32)}$ with $C$ depending on $c,$ $K$ and $T$, we get

$$\frac{\partial}{\partial t}||D^{s}\Theta_{N,c,K}||^{2}_{L^{2}}\leq ||D^{s}\Theta_{N,c,K}||_{L^{2}}(CN^{-(\beta-\frac32)}+C||\Theta_{N,c,K}||_{H^{s}}+C||\Theta_{N,c,K}||_{H^{s}}^2)$$

which gives us, using 
$$||\Theta_{N,c,K}||_{H^{s}}\leq C(||\Theta_{N,c,K}||_{L^{2}}+||D^{s}\Theta_{N,c,K}||_{L^{2}})\leq C(||D^{s}\Theta_{N,c,K}||_{L^{2}}+N^{-(2\beta-1)})$$
that

$$\frac{\partial}{\partial t}||D^{s}\Theta_{N,c,K}||_{L^{2}}\leq (CN^{-(\beta-\frac32)}+C||D^{s}\Theta_{N,c,K}||_{L^{2}}+C||D^{s}\Theta_{\lambda,K,N}||_{L^{2}}^2),$$
and using $||D^{s}\Theta_{\lambda,K,N}||_{L^{2}}\leq log(N)N^{-(\beta-\frac32)}$ and integrating we get

$$||D^{s}\Theta_{N,c,K}||_{L^2}\leq \frac{C (e^{Ct}-1)}{N^{\beta-\frac32}}.$$

\end{proof}

\subsection{Strong ill-posedness in supercritical Sobolev spaces}

Now we are ready to prove strong ill-posedness in supercritical Sobolev spaces:

\begin{theorem} (Strong ill-posedness in $H^{\beta}$)
For any $c_{0}>0$, $M>1$, $\beta\in (\frac32, 2)$ and $t_{*}>0$, we can find a $H^{\beta+\frac12}$ function $\theta_{0}(x)$  with $||\theta_{0}(x)||_{H^{\beta}}\leq c_{0}$ such that the unique solution $\theta(x,t)$ in $H^{\beta+\frac12}$ to the SQG equation (\ref{SQG}) with initial conditions $\theta_{0}(x)$ is such that $||\theta(x,t_{*})||_{H^{\beta}}\geq M  c_{0}$.
\end{theorem}

\begin{proof}

First we prove a bound for the pseudo-solution $\bar{\theta}_{N,c,K}$ defined in (\ref{pseudoNcK}). More precisely
$$||f_{2,c,K}(r) \frac{r_{c,K}^{\beta}sin(N\alpha)}{N^{\beta}}||_{L^{2}}\leq \frac{c r_{c,k}^{\beta}}{N^{\beta}},$$

and

$$||f_{2,c,K}(r) \frac{r_{c,K}^{\beta}sin(N\alpha)}{N^{\beta}}||_{H^{2}}\leq \frac{C c r_{c,k}^{\beta-2}}{N^{\beta-2}},$$

which in combination with the interpolation inequality for Sobolev spaces and the bounds for $f_{1,c,K}$ gives us
$$||\bar{\theta}_{N,c,K}(x,0)||_{H^{\beta}}\leq C_{1}c$$

with $C_{1}$ depending only on $\beta$.

Furthermore, at time  $t$ we have that our pseudo-solution fulfils 
$$||\bar{\theta}_{N,c,K}(x,t)-f_{1,c,K}||_{L^{2}}\leq \frac{c r_{c,k}^{\beta}}{N^{\beta}} $$

and we can find the lower bound for the $H^{1}$ norm of $\bar{\theta}_{N,c,K}-f_{1,c,K}$ by using 
\begin{align*}
    &\frac{\partial (\bar{\theta}_{N,c,K}-f_{1,c,K})}{\partial x_{1}}\\
    &=cos(\alpha)\frac{\partial (\bar{\theta}_{N,c,K}-f_{1,c,K})}{\partial r}-\frac{sin(\alpha)}{r}\frac{\partial (\bar{\theta}_{N,c,K}-f_{1,c,K})}{\partial \alpha}\\
\end{align*}
which gives us, after some trigonometric manipulations and using (\ref{velocidadrapida}) that, for $N$ large


$$|| \bar{\theta}_{N,c,K}(x,t)-f_{1,c,K}||_{H^{1}}\geq C\frac{ctKr_{c,K}^{\beta-1}}{N^{\beta-1}}$$

with $C$ a constant.

Furthermore, since $supp(\bar{\theta}_{N,c,K}-f_{1,c,K})\cap supp(f_{1,c,K})=\emptyset$ we have that
\begin{align*}
    &|| \bar{\theta}_{N,c,K}(x,t)||_{H^{1}}\geq || \bar{\theta}_{N,c,K}(x,t)-f_{1,c,K}||_{H^{1}} \\
    &\geq C\frac{ctKr_{c,K}^{\beta-1}}{N^{\beta-1}},\\
\end{align*}

for sufficiently large $N$.
On the other hand the interpolation inequality gives us

$$||\bar{\theta}_{N,c,K}(x,t)||_{H^{1}}\leq  ||\bar{\theta}_{N,c,K}(x,t)||_{H^{\beta}}^{\frac1\beta} ||\bar{\theta}_{N,c,K}(x,t)||_{L^{2}}^{\frac{\beta-1}{\beta}} $$

and using our bounds for $||\bar{\theta}_{N,c,K}(x,t)||_{L^{2}}$ and $||\bar{\theta}_{N,c,K}(x,t)||_{H^{1}}$ we get

$$||\bar{\theta}_{N,c,K}(x,t)||_{H^{\beta}}\geq C_{2}cK^{\beta}t^{\beta}$$

with $C_{2}$ depending only on $\beta$. Therefore, by choosing $c$, $K$  appropriately  we have that, for all $N$ big enough,
$$||\bar{\theta}_{N,c,K}(x,0)||_{H^{\beta}}\leq c_{0}$$

$$||\bar{\theta}_{N,c,K}(x,t^{*})||_{H^{\beta}}\geq 2Mc_{0}.$$

Now, considering the solution $\theta_{N,c,K}$ of (\ref{SQG}) with initial conditions $\bar{\theta}_{N,c,K}(x,0)$, we know that

$$||\theta_{N,c,K}(x,0)||_{H^{\beta}}\leq c_{0},$$

and, using lemma \ref{errorhssobolev}, if $N$ large,
\begin{align*}
    & ||\bar{\theta}_{N,c,K}(x,t^{*})-\theta_{N,c,K}(x,t^{*})||_{H^{\beta}}\\
    &\leq ||\bar{\theta}_{N,c,K}(x,t^{*})-\theta_{N,c,K}(x,t^{*})||_{H^{\beta+\frac12}}\leq \frac{Ct^{*}}{N^{\beta-\frac32}}\\
\end{align*}

and by taking $N$ big enough we can conclude 

$$||\theta_{N,c,K}(x,t^{*})||_{H^{\beta}}\geq ||\bar{\theta}_{N,c,K}(x,t^{*})||_{H^{\beta}}-||\bar{\theta}_{N,c,K}(x,t^{*})-\theta_{N,c,K}(x,t^{*})||_{H^{\beta}}\geq  Mc_{0}.$$

\end{proof}

\subsection{Non existence in supercritical Sobolev spaces}

In this section we prove the following theorem:

\begin{theorem}\label{nonexistencesupercritical}
(Non existence in $H^{\beta}$ in the supercritical case)
For any $t_{0}$, $c_{0}>0$ and $\beta\in(\frac32,2)$ we can find initial conditions $\theta_{0}(x)$,  with $||\theta_{0}(x)||_{H^{\beta}}\leq c_{0}$ such that there exists a solution $\theta(x,t)$ to (\ref{SQG}) with $\theta(x,0)=\theta_{0}(x)$ satisfying $||\theta(x,t)||_{H^{\beta}}=\infty$ for all $t\in(0,t_{0}]$. Furthermore, it is the only  solution with initial conditions $\theta_{0}(x)$ that satisfy $\theta(x,t)\in L^{\infty}_{t}C^{\gamma_{1}}_{x}\cap L^{\infty}_{t}L^{2}_{x}$ ($0<\gamma_{1}<\frac12$) with the property that $||\theta(x,t)||_{H^{\gamma_{2}}}\leq M(t)$ ($1<\gamma_{2}\leq \frac32$) for some function $M(t)$.
\end{theorem}

\begin{remark}
$M(t)$ is not necessarily a bounded function, so this rules out the existence of solutions in $\theta(x,t)\in L^{\infty}_{t}C^{\gamma_{1}}_{x}\cap L^{\infty}_{t}L^{2}_{x}$ such that $\theta(x,t)\in H^{\beta}$ for $t\in(0,t^{*}]$ with any $0<t^{*}\leq t_{0}$.
\end{remark}

\begin{proof}
Let's first note some of the properties that the pseudo-solutions $\bar{\theta}_{N,c,K}$ (for some fixed $\beta$) have:

\begin{itemize}
    \item $\bar{\theta}_{N,c,K}(x,t)$ is in $C^{\infty}$ for all $t\in[0,t_{0}]$, with  $||\bar{\theta}_{N,c,K}(x,t)||_{C^{k}}\leq CcN^{k-\beta}$, $||\bar{\theta}_{N,c,K}(x,t)||_{H^{k}}\leq CcN^{k-\beta}$ for any natural $k\geq 2$ , with the constant $C$ depending on $k$, $K$ and $t_{0}$. Also, for $\beta>s\geq 0$ we have $||\bar{\theta}_{N,c,K}(x,t)||_{H^{s}}\leq C_{1}cN^{s-\beta}+C_{2}c$ with $C_{1}$ depending on $K$, $s$ and $t_{0}$ and $C_{2}$ a constant.
    \item For $N$ large we have the lower bound $||\bar{\theta}_{N,c,K}(x,t)||_{H^{\beta}}\geq Cct^{\beta}K^{\beta}$ with $C$ a constant.
    
    \item $\bar{\theta}_{N,c,K}(x,t)$ is supported in the ball of radius $M$ centered at zero $B_{M}(0)$ for some $M$ independent of the values of the parameters.
\end{itemize}

Furthermore, we have the following result.

\begin{lemma}\label{evolucionvext}
Let $\tilde{\theta}_{N,c,K}$ with $\tilde{\theta}_{N,c,K}(x,0)=\bar{\theta}_{N,c,K}(x,0)$ and satisfying the equation

$$\frac{\partial \tilde{\theta}_{N,c,K}}{\partial t} + (v_{1}(\tilde{\theta}_{N,c,K})+v^{N,c,K}_{1,ext})\frac{\partial \tilde{\theta}_{N,c,K}}{\partial x_{1}} + (v_{2}(\tilde{\theta}_{N,c,K})+v^{N,c,K}_{2,ext})\frac{\partial\tilde{\theta}_{N,c,K}}{\partial x_{2}}= 0$$

with $$\frac{\partial v^{N,c,K}_{2,ext}}{\partial x_{2}}=-\frac{\partial v^{N,c,K}_{1,ext}}{\partial x_{1}}$$

and 
$$||v^{N,c,K}_{i,ext}||_{C^{3}}\leq CN^{-3} $$

with $C$ depending on $c$ and $K$.

Then for any $T>0$ we have that  if $N$ is big enough, then for $t\in[0,T]$ there exists a unique $\tilde{\theta}_{N,c,K}(x,t)\in H^{\beta+\frac12}$  and

$$||\bar{\theta}_{N,c,K}(x,t)-\tilde{\theta}_{N,c,K}(x,t)||_{L^{2}}\leq CtN^{-(2\beta-1)}$$

$$||\bar{\theta}_{N,c,K}(x,t)-\tilde{\theta}_{N,c,K}(x,t)||_{H^{\beta+\frac12}}\leq CtN^{-(\beta-\frac32)}$$

with $C$ depending on $c$, $K$ and $T$.
\end{lemma}

The local well posedness is straight forward since $v^{N,c,K}_{i,ext}$ for $i=1,2$ are $C^{3}$. As for the error bounds, they are obtained in the same way as in lemmas \ref{errorl2sobolev} and \ref{errorhssobolev}, i.e. studying the evolution equation for $\bar{\theta}_{N,c,K}(x,t)-\tilde{\theta}_{N,c,K}(x,t)$  now with  new terms depending on $v^{N,c,K}_{i,ext} \frac{\partial \tilde{\theta}_{N,c,K}(x,t)}{\partial x_{i}}$. These terms, however, are easily bounded by writing 

$$\tilde{\theta}_{N,c,K}(x,t)=(\tilde{\theta}_{N,c,K}(x,t)-\bar{\theta}_{N,c,K}(x,t))+\bar{\theta}_{N,c,K}(x,t)$$

and using our bounds for $v^{N,c,K}_{i,ext}$ and $\bar{\theta}_{N,c,K}(x,t)$.

This new lemma tells us that our pseudo-solutions as in (\ref{pseudoNcK}) stay close to other pseudo-solutions  that have the same initial conditions and an error term in the velocity if that term is small enough. Now, to obtain the initial conditions that will produce instantaneous loss of regularity, we consider

$$\theta(x,0)=\sum_{j=1}^{\infty}T_{R_{j}}(\bar{\theta}_{N_{j},c_{j},K_{j}}(x,0)),$$

with $T_{R}(f(x_{1},x_{2}))=f(x_{1}+R,x_{2})$, and $R_{j}$ yet to be fixed.

We will refer to the solution of (\ref{SQG}) with this initial conditions and $H^{\frac32}$ regularity  (if it exists) as $\theta(x,t)$, keeping in mind that it depends on the values for $R_{j}$, $N_{j}$, $c_{j}$, $K_{j}$, with $j\in \mathds{N}$.

We start by fixing $c_{j}$ and $K_{j}$ with the following properties:

1)
\begin{equation}\label{normabetapequeña}
    ||\bar{\theta}_{N_{j},c_{j},K_{j}}(x,0)||_{H^{\beta}}\leq c_{0}2^{-j},||\bar{\theta}_{N_{j},c_{j},K_{j}}(x,0)||_{L^{1}}\leq c_{0}2^{-j}.
\end{equation}

2) If $N_{j}$ large enough then
\begin{equation}\label{normabetagrande}
    ||\bar{\theta}_{N_{j},c_{j},K_{j}}(x,t)||_{H^{\beta}}\geq tc_{0}2^{j}
\end{equation}
 and
$$||\bar{\theta}_{N_{j},c_{j},K_{j}}(x,t)||_{H^{\frac32}}\leq c_{0}2^{-j}$$
for $t\in[0,t_{0}].$

This gives us a bound for the velocity generated by $\sum_{j=1}^{\infty}T_{R_{j}}(\bar{\theta}_{N_{j},c_{j},K_{j}}(x,t))$, which we will call $v_{max}$.

As for $R_{j}$, we will consider $R_{j}=R_{j-1}+D_{j}+D_{j-1}$, $R_{0}=0$, and we will take $D_{j}=j^{4}N_{j}^{4}+2M+8v_{max}t_{0}$.





Now, we say that a sequence $(w_{j}(x,t))_{j\in\mathds{N}}$ is in the space $W_{(N_{j})_{j\in \mathds{N}}, C_{0}}$ if it satisfies the following four conditions:

1) $w_{j}(x,t)\in H^{\beta+\frac12}$ for $t\in[0,t_{0}],$.

2) $w_{j}(x,t)$ satisfy
\begin{align}\label{wjecuacion}
    &\frac{\partial w_{j}(x,t)}{\partial t}=\\
    &-(v_{1}(w_{j})+v^{j}_{1,ext}(x,t))\frac{\partial w_{j}}{\partial x_{1}} - (v_{2}(w_{j})+v^{j}_{2,ext}(x,t))\frac{\partial w_{j}(x,t)}{\partial x_{2}}\nonumber
\end{align}

$$w_{j}(x,0)=T_{R_{j}}(\bar{\theta}_{N_{j},c_{j},K_{j}}(x,0))$$
with $v^{j}_{i,ext}$ fulfilling
$$||v^{j}_{i,ext}||_{C^{3}}\leq \frac{C_{0}}{j^{4}N_{j}^{3}}$$
and 
$$\frac{\partial v^{j}_{1,ext}}{\partial x_{1}}=-\frac{\partial v^{j}_{2,ext}}{\partial x_{2}}$$

3) $\sum_{j=1}^{\infty}||v_{i}(w_{j}(x,t))||_{L^{\infty}}\leq 2v_{max}$ for $i=1,2.$ 

4) $\sum_{j=1}^{\infty} ||w_{j}(x,t)||_{L^{2}}<\infty$ for $t\in [0,t_{0}]$.

The space $W_{(N_{j})_{j\in \mathds{N}}, C_{0}}$ is a complete metric space with the norm $$||(w_{j})_{j\in\mathds{N}}||_{W}:=sup_{j\in\mathds{N},i=1,2}\text{ess-sup}_{t\in[0,t_{0}]}j^2N_{j}^3||v^{j}_{i,ext}(x,t)||_{C^{3}}.$$ Note that the fourth condition is actually automatically satisfied if  the other three are satisfied, but we include it to emphasize that the norm is well defined. Note also that lemma \ref{evolucionvext} tells us that if $(N_j)_{j\in N}$ are big enough, the condition $$||v^{j}_{i,ext}||_{C^{3}}\leq \frac{C_{0}}{j^{4}N_{j}^{3}}$$
implies that there is a solution to (\ref{wjecuacion}) for $t\in[0,t_{0}]$ with

$$||T_{R_{j}}(\bar{\theta}_{N_{j},c_{j},K_{j}}(x,t))-w_{j}(x,t)||_{L^{2}}\leq CtN_{j}^{-(2\beta-1)}$$
\begin{equation}\label{cotahbeta+12}
    ||T_{R_{j}}(\bar{\theta}_{N_{j},c_{j},K_{j}}(x,t))-w_{j}(x,t)||_{H^{\beta+\frac12}}\leq CtN_{j}^{-(\beta-\frac32)},
\end{equation}
so that the solution satisfies condition 1 and, by taking $(N_{j})_{j\in \mathds{N}}$ large,  also condition 3.

Now, given an element $(w_{j})_{j\in \mathds{N}}$ in $W_{(N_{j})_{j\in \mathds{N}}, C_{0}}$ we define,

$$v^{j_{0}}_{1,ext}((w_{j})_{j\in \mathds{N}})=\frac{\partial}{\partial x_{2}}\bigg(\Lambda^{-1}[(\sum_{j=1}^{\infty}w_{j})-w_{j_{0}}]T_{R_{j_{0}}}\phi(x)\bigg),$$

$$v^{j_{0}}_{2,ext}((w_{j})_{j\in \mathds{N}})=-\frac{\partial}{\partial x_{1}}\bigg((\Lambda^{-1}[(\sum_{j=1}^{\infty}w_{j})-w_{j_{0}}]T_{R_{j_{0}}}\phi(x)\bigg),$$
where $\phi(x)$ is a smooth $C^{\infty}$ function with $\phi(x)=1$ if $x\in B_{4v_{max}+M}(0)$ and $\phi(x)=0$ if $|x|\geq 8vmax+M$.

Note that $||v^{j_{0}}_{i,ext}((w_{j})_{j\in \mathds{N}})||_{C^{3}}\leq \frac{Cc}{j_{0}^{4}N_{j_{0}}^{4}}$, and thus $||v^{j_{0}}_{i,ext}((w_{j})_{j\in \mathds{N}})||_{C^{3}}\leq \frac{C_{0}}{j_{0}^{4}N_{j_{0}}^{3}}$ if $N_{j_{0}}$ is large. Furthermore, if $x\in B_{4vmax+M}(-R_{j_{0}},0)$, then
\begin{equation}\label{vextfacil}
    v^{j_{0}}_{i,ext}((w_{j})_{j\in \mathds{N}})=v_{i} ((\sum_{j=1}^{\infty}w_{j})-w_{j_{0}})
\end{equation}

and, since, $supp(\theta_{j_{0}}(x,t))\subset B_{4vmax+M}(-R_{j_{0}},0)$, we could actually use (\ref{vextfacil}) as our definition of $v^{j_{0}}_{i,ext}((w_{j})_{j\in \mathds{N}})$ without changing anything.

This allows us to define the operator $G$ over a sequence $w$ in the space $ W_{(N_{j})_{j\in\mathds{N}},C_{0}}$ as

$$G(w)=(G^{j}(w))_{j\in \mathds{N}},$$

with $G^{j}(w)(x,t)$ the only $H^{\beta+\frac12}$ function for $t\in[0,t_{0}]$ satisfying 

\begin{align*}
    &\frac{\partial G^{j}(w)}{\partial t} + (v_{1}(G^{j}(w))+v^{j}_{1,ext}(w))\frac{\partial G^{j}(w)}{\partial x_{1}} \\
    &+ (v_{2}(G^{j}(w))+v^{j}_{2,ext}(w))\frac{\partial G^{j}(w)}{\partial x_{2}}= 0,\\
\end{align*}

$$G^{j}(w)(x,0)=T_{R_{j}}(\bar{\theta}_{N_{j},c_{j},K_{j}}(x,0)).$$

The operator $G$ maps (for $(N_{j})_{j\in\mathds{N}}$ large) $ W_{(N_{j})_{j\in\mathds{N}},C_{0}}$ to $ W_{(N_{j})_{j\in\mathds{N}},C_{0}}$ and actually, if we can find a point $w\in W_{(N_{j})_{j\in\mathds{N}},C_{0}}$ such that $G(w)=w$, then 
$$\theta(x,t)=\sum_{j=1}^{\infty}w_{j}(x,t)$$
is a solution to (\ref{SQG}) with initial conditions
$$\theta(x,0)=\sum_{j=1}^{\infty}T_{R_{j}}(\bar{\theta}_{N_{j},c_{j},K_{j}}(x,0)).$$

If we now consider two sequences $w^{1}=(w^{1}_{j})_{j\in\mathds{N}},w^{2}=(w^{2}_{j})_{j\in\mathds{N}}$ $\in W_{(N_{j})_{n\in\mathds{N}},C_{0}}$ and we define 
$||w^{1}-w^{2}||_{L^{2}}= sup_{t\in[0,t_{0}]}\sum_{j=1}^{\infty}||w^{1}_{j}-w^{2}_{j}||_{L^{2}}$
we can compute $||G(w^{1})-G(w^{2})||_{L^{2}}$, by defining $\tilde{w}_{j}=G^{j}(w^{1})-G^{j}(w^{2})$, since it fulfills the evolution equation

\begin{align*}
    \frac{\partial\tilde{w}_{j}}{\partial t}
    &= -\frac{\partial (G^{j}w^{1})}{\partial x_{1}}v_{1}(\tilde{w}_{j})-\frac{\partial\tilde{w}_{j}}{\partial x_{1}}v_{1}(G^{j}(w^{2}))\\
    &-\frac{\partial (G^{j}w^{1})}{\partial x_{2}}v_{2}(\tilde{w}_{j})-\frac{\partial\tilde{w}_{j}}{\partial x_{2}}v_{2}(G^{j}(w^{2}))\\
    &-\frac{\partial (G^{j}w^{1})}{\partial x_{1}}v^{j}_{1,ext}(w^{1}-w^{2})-\frac{\partial\tilde{w}_{j}}{\partial x_{1}}v^{j}_{1,ext}(w^{2})\\
    &-\frac{\partial (G^{j}w^{1})}{\partial x_{2}}v^{j}_{2,ext}(w^{1}-w^{2})-\frac{\partial\tilde{w}_{j}}{\partial x_{2}}v^{j}_{2,ext}(w^{2}).\\
\end{align*}

This gives us a bound for the evolution of the $L^{2}$ norm of $\tilde{w}_{j}$

\begin{align*}
    &\frac{\partial ||\tilde{w}_{j}||_{L^{2}}}{\partial t}\leq C ||G^{j}(w^1)||_{C^{1}}||\tilde{w}_{j}||_{L^{2}}+2 ||G^{j}(w^1)||_{C^{1}}\frac{||w^{1}-w^{2}||_{W}}{j^4 N_{j}^{3}}.\\
\end{align*}

But for $N_{j}$ large we can bound $||G^{j}(w^{1})||_{C^{1}}$ by some constant $\bar{C}_{j}$ using (\ref{cotahbeta+12}) , and thus we obtain, for $t\in[0,t_{0}]$

\begin{align*}
    &||\tilde{w}_{j}(x,t)||_{L^{2}}\leq C\tilde{C}_{j}(e^{Ct_{0}}-1) \frac{||w^{1}-w^{2}||_{W}}{j^4 N_{j}^{3}}\\
\end{align*}
and for $N_{j}$ large
\begin{align*}
    &||\tilde{w}_{j}(x,t)||_{L^{2}}\leq \epsilon \frac{||w^{1}-w^{2}||_{W}}{j^{4}}\\
\end{align*}
with $\epsilon$ as small as we want. Adding over all $j$ we obtain, for $t\in[0,t_{0}]$
    
    \begin{align*}
    &||G(w^{1})(x,t)-G(w^{2})(x,t)||_{L^{2}}\leq C\epsilon ||w^{1}-w^{2}||_{W}.\\
\end{align*}
But we have that

$$||v^{j}_{i,ext}(G(w^1))-v^{j}_{i,ext}(G(w^2))||_{C^{3}}j^4 N_{j}^{3}\leq \frac{C}{N_{j}}\sum_{j=1}^{\infty}||G(w^{1})-G(w^{2})||_{L^{2}}$$
and thus for $(N_{j})_{j\in\mathds{N}}$ big enough 
$$||G(w^{1})-G(w^{2})||_{W}\leq \epsilon ||w^{1}-w^{2}||_{W}$$
with $\epsilon$ arbitrarily small and in particular the map $G$ is a contraction.
Furthermore  the set $W_{(N_{j})_{n\in\mathds{N}},C_{0}}$ is not empty, since it includes at least the point where $v^{j}_{i,ext}=0$ when $(N_{j})_{j\in\mathds{N}}$ is large, and therefore, using the Banach point fixed theorem there exists $G(w)=w\in W_{(N_{j})_{n\in\mathds{N}},C_{0}}$. But as we pointed out earlier that implies that $w$ is a solution to (\ref{SQG}) with initial conditions
    
    $$\theta(x,0)=\sum_{j=1}^{\infty}T_{R_{j}}(\bar{\theta}_{N_{j},c_{j},K_{j}}(x,0)).$$
    
    Properties  (\ref{normabetapequeña}),(\ref{normabetagrande}) and (\ref{cotahbeta+12})  finish the proof that a solution with the desired properties of theorem \ref{nonexistencesupercritical} exists.

    For uniqueness in the space mentioned we call $\theta_{1}(x,t)$ the solution we constructed above and assume the existence of another solution $\theta_{2}(x,t)\in L^{\infty}_{t}C^{\gamma_{1}}_{x}\cap L^{\infty}_{t}L^{2}_{x}$ ($0<\gamma_{1}<\frac12$) with the property that $||\theta(x,t)||_{H^{\gamma_{2}}}\leq M(t)$ ($1<\gamma_{2}\leq\frac32$) for some function $M(t)$. In particular (since it is in $ L^{\infty}_{t}C^{\gamma_{1}}_{x}$), there exists a certain constant $v_{2,max}$ such that $||v_{i}(\theta_{2})||_{L^{\infty}}\leq v_{2,max}$. We start by studying the uniqueness for $t\in[0,\text{min}(t^{*},t_{0})]$ with $t^{*}v_{2,max}=4t_{0}v_{max}$. In particular, we have that $supp(\theta_{2}(x,t))\subset \cup_{j\in\mathds{N}} T_{R_{j}} (B_{t_{0}v_{max}+M}(0))$. We define 
    $$\theta^{j}_{1}(x,t)=1_{B_{4t_{0}v_{max}+M}(-R_{j},0)}\theta_{1}(x,t)$$
    $$\theta^{j}_{2}(x,t)=1_{B_{4t_{0}v_{max}+M}(-R_{j},0)}\theta_{2}(x,t).$$
    
    If we define $\Theta^{j}:=\theta^{j}_{2}-\theta^{j}_{1}$, $\Theta:=\theta_{2}-\theta_{1}$, we get

\begin{align*}
    \frac{\partial\Theta^{j}}{\partial t}
    &= -\frac{\partial \theta^{j}_{1}}{\partial x_{1}}v_{1}(\Theta^{j})-\frac{\partial\Theta^{j}}{\partial x_{1}}v_{1}(\Theta^{j})-\frac{\partial \theta^{j}_{1}}{\partial x_{2}}v_{2}(\Theta^{j})-\frac{\partial\Theta^{j}}{\partial x_{2}}v_{2}(\Theta^{j})\\
    &-\frac{\partial \Theta^{j}}{\partial x_{1}}v_{1}(\theta^{j}_{1})-\frac{\partial \Theta^{j}}{\partial x_{2}}v_{2}(\theta^{j}_{1})-\frac{\partial \theta^{j}_{1}}{\partial x_{1}}v_{1}(\Theta-\Theta^{j})-\frac{\partial\Theta^{j}}{\partial x_{1}}v_{1}(\Theta-\Theta^{j})\\
    &-\frac{\partial \theta^{j}_{1}}{\partial x_{2}}v_{2}(\Theta-\Theta^{j})-\frac{\partial\Theta^{j}}{\partial x_{2}}v_{2}(\Theta-\Theta^{j})\\
    &-\frac{\partial\Theta^{j}}{\partial x_{1}}v_{1}(\theta_{1}-\theta_{1}^{j})-\frac{\partial\theta^{j}}{\partial x_{2}}v_{2}(\theta_{1}-\theta_{1}^{j})\\
\end{align*}
which give us

\begin{align*}
    \frac{\partial ||\Theta_{j}||_{L^{2}}}{\partial t}\leq C||\theta^{j}_{1}||_{C^{1}}||\Theta_{j}||_{L^{2}}+C||\theta^{j}_{1}||_{C^{1}}\frac{||\Theta||_{L^{2}}}{j^{4}N_{j}^4}
\end{align*}
and by taking $N_{j}$ big we get

\begin{align*}
     ||\Theta_{j}||_{L^{2}}\leq \frac{\epsilon||\Theta||_{L^{2}}}{j^{4}}
\end{align*}
and adding over all $j$ and taking $\epsilon$ small

\begin{align*}
     ||\Theta||_{L^{2}}\leq \frac{||\Theta||_{L^{2}}}{2}
\end{align*}

and thus $||\Theta||_{L^{2}}$ for $t\in[0,t^{*}]$. Iterating the argument allows us to prove $||\Theta||_{L^{2}}=0$ for $t\in[0,t_{0}]$.

\end{proof}

\section{Strong ill-posedness in the critical Sobolev space  $H^{2}$}
For this section, we will consider solutions of (\ref{SQG}) that are in layers around zero, each one closer to the origin, so that the exterior layers effect over the inner layers will give us (in the limit) an evolution system of the form

$$\frac{\partial \bar{\theta}}{\partial t}+(v_{1}(\bar{\theta})+K(t)x_{1})\frac{\partial \bar{\theta}}{\partial x_{1}}+ (v_{2}(\bar{\theta})-K(t)x_{2})\frac{\partial \bar{\theta}}{\partial x_{2}}=0$$
$$v_{1}=-\frac{\partial}{\partial x_{2}}\Lambda^{-1} \bar{\theta}=-\mathcal{R}_{2}\bar{\theta}$$
$$v_{2}=\frac{\partial}{\partial x_{1}}\Lambda^{-1} \bar{\theta}=\mathcal{R}_{1}\bar{\theta}$$
$$\bar{\theta}(x,0)=\theta_{0}(x).$$

But first we need to obtain an expression for $\frac{\partial v_{i}(\theta)(0)}{\partial x_{j}}$  ($i,j=1,2$) for $\theta$ with support far away from $0$. We consider first $i=1$. We have

$$v_{1}(\theta)=\frac{\Gamma (3/2)}{\pi^{3/2}}P.V. \int_{\mathds{R}^2} \frac{(-x_{2}+y_{2})\theta(y)}{|x-y|^{3}}dy_{1}dy_{2}. $$

For $\theta$ with support far away from $x=0$ we can just differentiate under the integral sign and when we evaluate at $x=0$ this yields

$$\frac{\partial v_{1}(\theta)}{\partial x_{1}}(x=0)=\frac{\Gamma (3/2)}{\pi^{3/2}}P.V. \int_{\mathds{R}^2} -3y_{1}\frac{y_{2}\theta(y)}{|y|^{5}}dy_{1}dy_{2}, $$

$$\frac{\partial v_{1}(\theta)}{\partial x_{2}}(x=0)=\frac{\Gamma (3/2)}{\pi^{3/2}}P.V. \int_{\mathds{R}^2} (\frac{-3y^2_{2}\theta(y)}{|y|^{5}}+\frac{\theta(y)}{|y|^{3}})dy_{1}dy_{2}. $$

We will consider $\theta(x_{1},x_{2})$ satisfying $\theta(-x_{1},x_{2})=-\theta(x_{1},x_{2})$, $\theta(x_{1},-x_{2})=-\theta(x_{1},x_{2})$, so

$$\frac{\partial v_{1}(\theta)}{\partial x_{1}}(x=0)=\frac{4\Gamma (3/2)}{\pi^{3/2}}P.V. \int_{\mathds{R}_{+}^2} -3y_{1}\frac{y_{2}\theta(y)}{|y|^{5}}dy_{1}dy_{2}, $$

$$\frac{\partial v_{1}(\theta)}{\partial x_{2}}(x=0)=0.$$

If we take a look at the expression for $\frac{\partial v_{1}(\theta)}{\partial x_{1}}$ in polar coordinates and combining all the constant into a certain $C_{0}>0$ we obtain

$$\frac{\partial v_{1}(\theta)}{\partial x_{1}}(x=0)=-C_{0} P.V. \int_{\mathds{R}_{+}\times [0,\pi/2]} \frac{sin(2\alpha')\theta(r',\alpha')}{(r')^{2}}dr'd\alpha'.$$

The expressions for $v_{2}$ are obtained the same way and in fact we have

$$\frac{\partial v_{2}(\theta)}{\partial x_{1}}(x=0)=0,$$
$$\frac{\partial v_{2}(\theta)}{\partial x_{2}}(x=0)=C_{0} P.V. \int_{\mathds{R}_{+}\times [0,\pi/2]} \frac{sin(2\alpha')\theta(r',\alpha')}{(r')^{2}}dr'd\alpha'.$$

Analogously, the second derivatives of $v_{i}$ all vanish.

We will be interested in studying the evolution of initial conditions of the form

$$\sum_{j=1}^{J}\frac{f(b^{-j}r)b^{j}sin(2\alpha)}{j}$$
for $f(r)$ a positive $C^{\infty}$ function with compact support and $\frac12>b>0$. More precisely, we would like to study the behaviour of the unique $H^{4}$ solution with said initial conditions when $b$ tends to zero. One could think that  we can just check the evolution of each of the terms $\frac{f(b^{-j}r)b^{j}sin(2\alpha)}{j}$ and then add them together, hoping that the interaction between them gets small as $b\rightarrow 0$. However this is not true, and we get an interaction depending on $\frac{\partial v_{i}}{\partial x_{i}}$. To get specific results, we fix some positive radial function $f$ in $C^{\infty}$ with $supp(f)\subset \{r\in [1/2,3/2]\}$ and $||f(r)sin(2\alpha)||_{H^{4}}=1$. We define $\theta_{c,J,b}$ as the unique $H^{4}$ solution of

\begin{equation*}
\frac{\partial \theta_{c,J,b}}{\partial t} + v_{1}(\theta_{c,J,b})\frac{\partial \theta_{c,J,b}}{\partial x_{1}} + v_{2}(\theta_{c,J,b})\frac{\partial \theta_{c,J,b}}{\partial x_{2}}= 0
\end{equation*}
with
$$v_{1}(\theta_{c,J,b})=-\frac{\partial}{\partial x_{2}}\Lambda \theta_{c,J,b}=-\mathcal{R}_{2}\theta_{c,J,b}$$
$$v_{2}(\theta_{c,J,b})=\frac{\partial}{\partial x_{1}}\Lambda \theta_{c,J,b}=\mathcal{R}_{1}\theta_{c,J,b}$$

\begin{equation}\label{coniniclayer}
  \theta_{c,J,b}(x,0)=\sum_{j=1}^{J}c\frac{f(b^{-j}r)b^{j}sin(2\alpha)}{j },\quad \frac12>b>0.  
\end{equation}

Note that the odd symmetry is preserved in time.

A few comments need to be made regarding the properties of the transformation $h(r,\alpha)\rightarrow \frac{h(\lambda r,\alpha)}{\lambda}$ (or equivalently $h(x)\rightarrow \frac{h(\lambda x)}{\lambda}$). We have that

\begin{itemize}
    \item If $\lambda>1$, then $||\frac{h(\lambda r,\alpha)}{\lambda}||_{H^{2}}\leq ||h(r,\alpha)||_{H^{2}}$.
    \item If $h(r,\alpha,t)$ is a solution to (\ref{SQG}) with initial conditions $h(r,\alpha,0)$, then $\frac{h(\lambda r,\alpha,t)}{\lambda}$ is a solution to (\ref{SQG}) with initial conditions $\frac{h(\lambda r,\alpha,0)}{\lambda}$.
    \item For $i=1,2$, $j=1,2$ we have $v_{i}(\frac{h(\lambda \cdot , \cdot)}{\lambda})(\frac{r}{\lambda},\alpha)=\frac{1}{\lambda}v_{i}(h(\cdot,\cdot ))(r,\alpha)$,$\quad$ $ \frac{\partial v_{i}(\frac{h(\lambda \cdot ,\cdot)}{\lambda})}{\partial x_{j}}(\frac{r}{\lambda},\alpha)=\frac{\partial v_{i}(h(\cdot,\cdot ))}{\partial x_{j}}(r,\alpha)$.
\end{itemize}

The initial conditions in (\ref{coniniclayer}) fulfil that, taking $c$ small and $J$ big, they have an arbitrarily small $H^{2}$ norm and an arbitrarily big value of $|\frac{\partial v_{1}(\theta_{c,J,b})}{\partial x_{1}}(0,t=0)|$. If  $|\frac{\partial v_{1}(\theta_{c,J,b})}{\partial x_{1}}(0,t)|$ remained big for a long enough time and $\theta$ remained sufficiently regular during that time, we could then use a small perturbation around $x=0$ to obtain a big growth in some $H^{s}$ norm. 

The main problem here is that we cannot assure existence for sufficiently long times using just the a priori bounds, so we need some extra machinery to be able to work with these solutions. For that we consider $\tilde{C}$ the constant fulfilling that, for any $H^{4}$ solution of SQG (\ref{SQG}) we have
\begin{equation}\label{ctildecuad}
\frac{\partial ||\theta(x,t)||_{H^{4}}}{\partial t}\leq \tilde{C}||\theta(x,t)||^{2}_{H^{4}}
\end{equation}
and, fixed constants $t_{0},K>0$, we define $t^{crit}_{t_{0},K,c,J,b}$ as the biggest time  fulfilling that, for all times $t$ satisfying $t^{crit}_{t_{0},K,c,J,b}\geq t\geq 0$ we have

\begin{itemize}
    \item $t\leq t_{0}$.
    \item If $x\in [\frac12b^{n},\frac32b^{n}] $ for $1\leq n\leq J$, then $\phi_{c,J,b}(x,t)\in [b^{n+\frac18},b^{n-\frac18}]$, with $\phi_{c,J,b}(x,t)$ the flow given by $$\frac{d\phi_{c,J,b}(x,t)}{dt}=v(\theta_{c,J,b}(x,t)).$$
    \item $||b^{-j}\theta_{c,J,b}(b^{j}x,t)1_{[b^{\frac18},b^{-\frac18}]}(r)||_{H^{4}}\leq \frac{1}{t_{0}\tilde{C}}$ for $1\leq j\leq J$.
    \item $\int_{0}^{t}|\frac{\partial v_{1}(\theta_{c,J,b})}{\partial x_{1}}(0,s)|ds\leq K$.
\end{itemize}

Let us make a few remark on these conditions. First, due to the odd symmetry of the solution and the initial conditions, $\frac{\partial v_{1}(\theta_{c,J,b})}{\partial x_{1}}$ is always negative and thus 
$$\int_{0}^{t}|\frac{\partial v_{1}(\theta_{c,J,b})}{\partial x_{1}}(0,s)|ds$$
is a monotone function with respect to $t$. Note also that we can check that the norm 
$$||b^{-j}\theta_{c,J,b}(b^{j}x,t)1_{[b^{\frac18},b^{-\frac18}]}(r)||_{H^{4}}$$
is continuous in time by checking the evolution equation for it and using that $\theta_{c,J,b}$ exists locally in time.
Also, depending on the choice of parameters $t^{crit}_{t_{0},K,c,J,b}$ could not exist ( the second and third condition may no bet satisfied for $t=0$), so we will only consider $c<\frac{1}{\tilde{Ct_{0}}}$ and $b<2^{-8}$ to avoid that. Finally, if we only consider the typical a priori bounds, the second and third conditions could make $t^{crit}_{t_{0},K,c,J,b}$ tend to zero as we make $b$ small, which would be a problem for our purposes. However, we have the following lemma.

\begin{lemma}\label{bpequeño}
Fixed $t_{0},K,c$ and $J$ fulfilling $c<\frac{e^{-6K}}{\tilde{C}t_{0}}$ and $K>\text{max}(1,t_{0})$, we have that, if $b$ is small enough, then the unique $H^{4}$ solution $\theta_{c,J,b}$ with initial conditions as in (\ref{coniniclayer}) satisfies

$$||b^{j}\theta_{c,J,b}(b^{-j}x,t)1_{[b^{\frac18},b^{-\frac18}]}(r)||_{H^{4}}< \frac{1}{t_{0}\tilde{C}}$$ 
for $1\leq j\leq J$, $t\in[0,t^{crit}_{t_{0},K,c,J,b}]$ and if $x\in [b^{n}\frac12,b^{n}\frac32] $ then $\phi_{c,J,b}(x,t)\in (b^{n+\frac18},b^{n-\frac18})$ if $0\leq t\leq t^{crit}_{t_{0},K,c,J,b}$.
\end{lemma}
\begin{proof}
Before we get into the proof, we need to define

$$k_{n}(t):=|\frac{\partial v_{1}(\theta_{c,J,b}1_{(b^{n+\frac18},\infty)}(r))}{\partial x_{1}}(0,t)|,$$
$$K_{n}(t):=\int_{0}^{t}k_{n}(s)ds.$$


We will study the evolution of $\theta_{j}:=\theta_{c,J,b}1_{[b^{j+\frac18},b^{j-\frac18}]}(r)$ (these functions obviously depend on $c,J$ and $b$, but we will omit this dependence to obtain a more compact notation). These functions satisfy the evolution equation

$$\frac{\partial \theta_{j}}{\partial t} + v_{1}(\theta_{j})\frac{\partial \theta_{j}}{\partial x_{1}} + v_{2}(\theta_{j})\frac{\partial \theta_{j}}{\partial x_{2}}+v_{2}(\theta_{c,J,b}-\theta_{j})\frac{\partial \theta_{j}}{\partial x_{2}}+ v_{1}(\theta_{c,J,b}-\theta_{j})\frac{\partial \theta_{j}}{\partial x_{1}}= 0. $$

Furthermore, we have that $\theta^{'}_{j}(x,t)=b^{-j}\theta_{j}(b^{j}x,t)$ fulfils the evolution equation
\begin{equation}\label{evolucionthetaj}
    \frac{\partial \theta^{'}_{j}}{\partial t} + v_{1}(\theta^{'}_{j})\frac{\partial \theta^{'}_{j}}{\partial x_{1}} + v_{2}(\theta^{'}_{j})\frac{\partial \theta^{'}_{j}}{\partial x_{2}}+v_{2}(\theta^{',j}_{c,J,b}-\theta^{'}_{j})\frac{\partial \theta^{'}_{j}}{\partial x_{2}}+ v_{1}(\theta^{',j}_{c,J,b}-\theta^{'}_{j})\frac{\partial \theta^{'}_{j}}{\partial x_{1}}= 0,
\end{equation}
with $\theta^{',j}_{c,J,b}(x,t):=b^{-j}\theta_{c,J,b}(b^{j}x,t).$

We want to obtain suitable bounds for the terms depending on $\theta^{',j}_{c,J,b}-\theta^{'}_{j}$. To do this we decompose $\theta^{',j}_{c,J,b}-\theta^{'}_{j}$ as

$$\theta^{',j}_{c,J,b}-\theta^{'}_{j}=\theta^{'}_{+,j}+\theta^{'}_{-,j}$$
with $\theta^{'}_{+,j}=(\theta^{',j}_{c,J,b}-\theta^{'}_{j})1_{[1,\infty]}(r)$ and $\theta^{'}_{-,j}=(\theta^{',j}_{c,J,b}-\theta^{'}_{j})1_{[0,1]}(r)$.

But $\theta^{'}_{-,j}$ satisfies that $||\theta^{'}_{-,j}||_{L^{1}}\leq Cb^{3}$, $d(supp(\theta^{'}_{j}),supp(\theta^{'}_{-,j}))\geq \frac{b^{\frac18}}{2}$, which gives us, if we define $v^{-,j}_{i}(x):=v_{i}(\theta^{'}_{-,j})(x)$

$$||v^{-,j}_{i}(x)1_{supp(\theta^{'}_{j})}||_{C^{4}}\leq Cb^{3-\frac68}.$$

For the term depending on $\theta^{'}_{+,j}$, we use that, for $k\geq 1$

$$||\theta^{',j}_{c,J,b}1_{[b^{-k+\frac18},b^{-k-\frac18}]}||_{L^{1}}\leq Cb^{-3k}$$
$$d(supp(\theta^{',j}_{c,J,b}1_{[b^{-k+\frac18},b^{-k-\frac18}]}),supp(\theta^{'}_{j}))\geq \frac{b^{-k+\frac18}}{2}$$

which gives us, after adding the contributions for all the $k$

$$|\frac{\partial^{2} v_{i}(\theta^{'}_{+,j})}{\partial^{2-j} x_{1}\partial^{j} x_{2}}(x)|\leq Cb^{\frac12}.$$

Therefore, using a second order Taylor expansion for the velocity we obtain that, for $|x|\leq b^{-\frac18}$ 

$$v_{1}(\theta^{'}_{+,j})=k_{j-1}(t)x_{1}+ v^{+,j,error}_{1}(x),$$

with $||v^{+,j,error}_{1}(x)1_{[b^{\frac18},b^{-\frac18}]}(r)||_{L^{\infty}}\leq Cb^{\frac14}$. Furthermore by computing the derivatives of $v_{1}(\theta^{'}_{+,j})$ we actually obtain $||v^{+,j,error}_{1}(x)1_{[b^{\frac18},b^{\frac18}]}(r)||_{C^{4}}\leq Cb^{\frac14}.$

Analogously, we have

$$v_{2}(\theta^{'}_{+,j})=-k_{j-1}(t)x_{2}+ v^{+,j,error}_{2}(x),$$

with $||v^{+,j,error}_{2}(x)1_{[b^{\frac18},b^{-\frac18}]}(r)||_{C^{4}}\leq Cb^{\frac14}.$

Writing $v^{error}_{i}:=v^{+,j,error}_{i}(x)+v^{-,j}_{i}(x)$, we get that (\ref{evolucionthetaj}) is equivalent to

$$\frac{\partial \theta^{'}_{j}}{\partial t} +( v_{1}(\theta^{'}_{j})+v^{error}_{1}+k_{j-1}x_{1})\frac{\partial \theta^{'}_{j}}{\partial x_{1}} + (v_{2}(\theta^{'}_{j})+v^{error}_{2}-k_{j-1}x_{2})\frac{\partial \theta^{'}_{j}}{\partial x_{2}}= 0, $$
with $||v^{error}_{i}||_{C^{4}}\leq Cb^{\frac14}$
To obtain the evolution of the $H^{4}$ norm, we note that, with our definition of $H^{4}$ norm

$$\frac{\partial || \theta^{'}_{j}||_{H^{4}}}{\partial t}=\sum_{i=0}^{4}\sum_{j=0}^{i}\frac{\partial ||\frac{\partial^{i}\theta^{'}_{j}}{\partial^{j}x_{1}\partial^{i-j} x_{2}}||_{L^{2}}}{\partial t}$$

and

\begin{align*}
    &\frac{\partial ||\frac{\partial^{i}\theta^{'}_{j}}{\partial^{j}x_{1}\partial^{i-j} x_{2}}||^{2}_{L^{2}} }{\partial t}\\
    &=2(\frac{\partial^{i}\theta^{'}_{j}}{\partial^{j}x_{1}\partial^{i-j} x_{2}},( v_{1}(\theta^{'}_{j})+v^{error}_{1}+k_{j-1}x_{1})\frac{\partial \theta^{'}_{j}}{\partial x_{1}} + (v_{2}(\theta^{'}_{j})+v^{error}_{2}-k_{j-1}x_{2})\frac{\partial \theta^{'}_{j}}{\partial x_{2}})_{L^{2}}.
\end{align*}

But using $||v^{error}_{i}||_{C^{4}}\leq Cb^{\frac14}$ and incompressibility we get, for $i=0,1,...,4$, $j=0,...,i$

$$|(\frac{\partial^{i}\theta^{'}_{j}}{\partial^{j}x_{1}\partial^{i-j} x_{2}},\frac{\partial^{i}(v^{error}_{1}\frac{\partial \theta^{'}_{j}}{\partial x_{1}})}{\partial^{j}x_{1}\partial^{i-j} x_{2}}+\frac{\partial^{i}(v^{error}_{2}\frac{\partial \theta^{'}_{j}}{\partial x_{2}})}{\partial^{j}x_{1}\partial^{i-j} x_{2}})_{L^{2}}|\leq Cb^{\frac14}||\theta^{'}_{j}||^{2}_{H^{4}}$$

and 

$$|(\frac{\partial^{i}\theta^{'}_{j}}{\partial^{j}x_{1}\partial^{i-j} x_{2}},\frac{\partial^{i}(k_{j-1}x_{1}\frac{\partial \theta^{'}_{j}}{\partial x_{1}})}{\partial^{j}x_{1}\partial^{i-j} x_{2}}-\frac{\partial^{i}(k_{j-1}x_{2}\frac{\partial \theta^{'}_{j}}{\partial x_{2}})}{\partial^{j}x_{1}\partial^{i-j} x_{2}})_{L^{2}}|\leq ik_{j-1}||\frac{\partial^{i}\theta^{'}_{j}}{\partial^{j}x_{1}\partial^{i-j} x_{2}}||^{2}_{L^{2}}$$

which gives us, by adding all the terms and including the contribution from the terms depending on $ v_{1}(\theta^{'}_{j})\frac{\partial \theta^{'}_{j}}{\partial x_{1}}$ and $ v_{2}(\theta^{'}_{j})\frac{\partial \theta^{'}_{j}}{\partial x_{2}}$

\begin{align}\label{ecuacionnormaj}
    &\frac{\partial ||\theta^{'}_{j}||_{H^{4}}}{\partial t}=\frac{\partial ||b^{j}\theta_{c,J,b}(b^{-j}x,t)1_{[b^{\frac18},b^{-\frac18}]}(r)||_{H^{4}}}{\partial t}\\
    &\leq (4k_{j-1}+Cb^{\frac14}) ||b^{j}\theta_{c,J,b}(b^{-j}x,t)1_{[b^{\frac18},b^{-\frac18}]}(r)||_{H^{4}}\nonumber\\
    &+\tilde{C}||b^{j}\theta_{c,J,b}(b^{-j}x,t)1_{[b^{\frac18},b^{-\frac18}]}(r)||^2_{H^{4}},\nonumber\\\nonumber
\end{align}
with $\tilde{C}$ given by (\ref{ctildecuad}).

Using that, by hypothesis 
$$||b^{j}\theta_{c,J,b}(b^{-j}x,t)1_{[b^{\frac18},b^{-\frac18}]}(r)||_{H^{4}}\leq \frac{1}{t_{0}\tilde{C}},$$ $$||b^{j}\theta_{c,J,b}(b^{-j}x,0)1_{[b^{\frac18},b^{-\frac18}]}(r)||_{H^{4}}\leq c$$ 
and integrating  (\ref{ecuacionnormaj}) we get

$$||b^{j}\theta_{c,J,b}(b^{-j}x,t)1_{[b^{\frac18},b^{-\frac18}]}(r)||_{H^{4}}\leq ce^{4K_{j-1}(t)+(\frac{1}{t_{0}}+Cb^{\frac14})t}$$

and using $K_{j-1}(t)\leq K$, and taking $b$ small enough 

$$||b^{j}\theta_{c,J,b}(b^{-j}x,t)1_{[b^{\frac18},b^{-\frac18}]}(r)||_{H^{4}}< ce^{6K}<\frac{1}{\tilde{C}t_{0}},$$

which gives us the first property we wanted.

As for the bounds for $\phi_{c,J,b}(x,t)$, we again work in the equivalent problem with $\theta'_{c,J,b}$ and note that we just proved that  
$$|v_{i}(\theta^{',j}_{c,J,b})(x)1_{[b^{\frac18},b^{-\frac18}]}|\leq (k_{J}(t)+Cb^{\frac14})|x|+|v_{i}(\theta^{'}_{j})|(x),$$
and since $|v_{i}(\theta^{'}_{j})|\leq min(C,C|x|)$ (by using our bounds in $H^{4}$ plus $v_{i}(\theta^{'}_{j})(x=0)=0$), integrating in time we have that, for $b$ small, the particles under that flow starting in $[\frac12,\frac32]$ will stay in $(e^{-C},e^{C})\subset(b^{\frac18},b^{-\frac18})$, with $C$ depending on $K$ and $t_{0}$ and we are done by returning to the original problem.
\end{proof}

Note that last lemma tells us that for $b$ small enough, at $t=t^{crit}_{t_{0},K,c,J,b}$, either $t=t_{0}$ or $\int_{0}^{t^{crit}_{t_{0},K,c,J,b}}|\frac{\partial v_{1}(\theta_{c,J,b})}{\partial x_{1}}(0,s)|ds=K$. Our next goal is to prove that, if the right conditions are met, we will actually have $\int_{0}^{t^{crit}_{t_{0},K,c,J,b}}|\frac{\partial v_{1}(\theta_{c,J,b})}{\partial x_{1}}(0,s)|ds=K$.

\begin{lemma}\label{estiramientogrande}
For fixed $t_{0},K$ and $c$ fulfilling $c<\frac{e^{-6K}}{\tilde{C}t_{0}}$ and $K>\text{max}(1,t_{0})$, we can find $J$ and $b$ such that at time $t=t^{crit}_{t_{0},K,c,J,b}$ we have that

$\int_{0}^{t^{crit}_{t_{0},K,c,J,b}}|\frac{\partial v_{1}(\theta_{c,J,b})}{\partial x_{1}}(0,s)|ds= K$.
\end{lemma}

\begin{proof}
We start by studying the trajectories of particles with $|x|\in[b^{J+\frac18},b^{-\frac18}]$.




In the proof of lemma \ref{bpequeño}  we obtained that, for $|x|\in[b^{\frac18},b^{-\frac18}]$,
\begin{equation}
    v_{1}(\theta^{',j}_{c,J,b})=v(\theta^{'}_{j})+v^{error}_{1}(x)+k_{j-1}(t)x_{1}
\end{equation}
$$v_{2}(\theta^{',j}_{c,J,b})=v(\theta^{'}_{j})+v^{error}_{2}(x)-k_{j-1}(t)x_{2}$$
(let us remember that here $\theta^{'}_{j}$ actually depends on $c,J$ and $b$ but we omit it), with $||v^{error}_{i}(x)||_{C^{4}}\leq C_{1}b^{\frac14}$ for $i=1,2$, with $C_{1}$ depending on $c,j$ and $J$, and $||v(\theta^{'}_{j})||_{C^{1}}\leq C_{2}$ with $C_{2}$ depending on $t_{0}$. By returning to the original problem, we get that, for $|x|\in[b^{j+\frac18},b^{j-\frac18}]$

\begin{equation}\label{ecuacionesj}
    v_{1}(\theta_{c,J,b})=v(\theta_{j})+v^{error,j}_{1}(x)+k_{j-1}(t)x_{1}
\end{equation}
$$v_{2}(\theta_{c,J,b})=v(\theta_{j})+v^{error,j}_{2}(x)-k_{j-1}(t)x_{2}$$

with $||v^{error}||_{C^{1}}\leq Cb^{\frac14}$ and $||v(\theta_{j})||_{C^{1}}\leq C_{2}$ with $C_{2}$ depending on $t_{0}$.

We are interested in studying the $\phi$ associated to this problem in polar coordinates for particles starting in $(r,\alpha)\in([\frac12,\frac32],[0,2\pi])$. We study separately the evolution of the radial coordinate and of the angular coordinate for simplicity.

For the radial coordinate, if we call $\phi^{j}_{r}(r_{0},\alpha_{0},t)$ the flow associated to (\ref{ecuacionesj}) that gives us the radial coordinate of the particle that was initially in $(r_{0},\alpha_{0})$, we have that,

$$\frac{\phi^{j}_{r}(r_{0},\alpha_{0},t)}{r_{0}}\leq e^{\int^{t}_{0} k_{j-1}(s)ds+C_{1}b^{\frac14}t+C_{2}t}\leq e^{K+C_{1}b^{\frac14}t+C_{2}t}.$$

As for the change in the angular coordinate, we are interested in finding bounds for how fast a particle can approach the lines $\alpha=i\frac{\pi}{2}$, $i=0,1,2,3$. All four cases are equivalent, so we will consider $i=0$. We have that

$$v_{\alpha}(r,0,t)=0$$
and, since for $i=1,2$ $||\frac{\partial v_{\alpha}}{\partial x_{i}}||_{C^{1}}\leq C(|k_{j-1}|+C_{1}b+C_{2}) $ (with $C$ a universal positive constant) we get, defining $\phi^{j}_{\alpha}$ similarly as we did with $\phi^{j}_{r}(r_{0},\alpha_{0},t)$,

$$\frac{\phi^{j}_{\alpha}(r_{0},\alpha_{0},t)}{\alpha_{0}}\geq e^{-C(\int^{t}_{0} k_{j-1}(s)ds+C_{1}b^{\frac14}t+C_{2}t)}\geq e^{-C(K+C_{1}b^{\frac14}t+C_{2}t)}.$$

Now we are ready to obtain bounds for
$$\int_{0}^{t^{crit}_{t_{0},K,c,J,b}}|\frac{\partial v_{1}(\theta_{c,J,b})}{\partial x_{1}}(0,s)|ds.$$

Since the transformation

$$\theta_{c,J,b}(x)\rightarrow \frac{\theta_{c,J,b}(\lambda x)}{\lambda}$$
does not change the value of $\frac{\partial v_{1}(\theta_{c,J,b})}{\partial x_{1}}(0,s)$ and by linearity, we have that, for $s=0$ we can compute

$$\frac{\partial v_{1}(\theta_{c,J,b})}{\partial x_{1}}(x=0,t=0)=\sum_{j=1}^{J}\frac{c}{j}\frac{\partial v_{1}(f(r)sin(2\alpha))}{\partial x_{1}}(x=0)=C(\sum_{j=1}^{J}\frac{c}{j})\geq Ccln(J).$$

For times $t>0$, writing for the flow map $\phi_{c,J,b}(x,t)=(\phi_{1,c,J,b}(x,t),\phi_{2,c,J,b}(x,t))$

\begin{align*}
    &|\frac{\partial v_{1}(\theta_{c,J,b}(r,\alpha,t))}{\partial x_{1}}|=C\int_{\mathds{R}_{+}^2} y_{1}\frac{y_{2}\theta_{c,J,b}(y,t)}{|y|^{5}}dy_{1}dy_{2}\\
    &=C\int_{\mathds{R}_{+}^2} y_{1}\frac{y_{2}\theta_{c,J,b}(\phi_{c,J,b}^{-1}(y,t),0)}{|y|^{5}}dy_{1}dy_{2}\\
    &=C\int_{\mathds{R}_{+}^2} \phi_{1,c,J,b}(\tilde{y},t)\frac{\phi_{2,c,J,b}(\tilde{y},t)\theta_{c,J,b}(\tilde{y},0)}{|\phi_{c,J,b}(\tilde{y},t)|^{5}}d\tilde{y}_{1}d\tilde{y}_{2}\\
    &=C\int_{\mathds{R}_{+}^2}\phi_{1,c,J,b}(\tilde{y},t)\frac{\phi_{2,c,J,b}(\tilde{y},t)}{|\phi_{c,J,b}(\tilde{y},t)|^{5}}\frac{|\tilde{y}|^{5}}{\tilde{y}_{1}\tilde{y}_{2}} \frac{\tilde{y}_{1}\tilde{y}_{2}\theta_{c,J,b}(\tilde{y},0)}{|\tilde{y}|^{5}}d\tilde{y}_{1}d\tilde{y}_{2}\\
\end{align*}
with  $C$ a constant, but (passing to polar coordinates to obtain the bound more easily)

\begin{align*}
    &\phi_{1,c,J,b}(x,t)\frac{\phi_{2,c,J,b}(x,t)}{|\phi_{c,J,b}(x,t)|^{5}}|\frac{|x|^{5}}{x_{1}x_{2}}\\
    &=\frac{sin(2\phi^{\alpha}_{c,J,b}(r,\alpha))}{sin(2\alpha)}\frac{r^{3}}{\phi^{r}_{c,J,b}(r,\alpha)}\geq e^{-C(K+c_{1}b^{\frac14}t+C_{2}t)}\\
\end{align*}
for some $C$, and thus

\begin{align*}
    &|\frac{\partial v_{1}(\theta_{c,J,b}(r,\alpha,t))}{\partial x_{1}}|\\
    &\geq Ce^{-C(K+c_{1}b^{\frac14}t+C_{2}t)}\int_{\mathds{R}_{+}^2}\ \frac{\tilde{y}_{1}\tilde{y}_{2}\theta_{c,J,b}(\tilde{y},0)}{|\tilde{y}|^{5}}d\tilde{y}_{1}d\tilde{y}_{2}\\
\end{align*}
and integrating in time
$$\int_{0}^{t^{crit}_{t_{0},K,c,J,b}}|\frac{\partial v_{1}(\theta_{c,J,b})}{\partial x_{1}}(0,s)|ds\geq t^{crit}_{t_{0},K,c,J,b}Ccln(J)e^{-C(K+C_{1}b^{\frac14}t_{0}+C_{2}t_{0})} $$

To finish our prove, we just fix some $K,$ $t_{0}$ and $c$ fulfilling our hypothesis, we take $J$ big enough so that 

$$t_{0}Ccln(J)e^{-C(K+C_{2}t_{0})}>K+1$$

and then take $b$ small enough so that using lemma \ref{bpequeño} either $t_{0}=t^{crit}_{t_{0},K,c,J,b}$ or

$$\int_{0}^{t^{crit}_{t_{0},K,c,J,b}}|\frac{\partial v_{1}(\theta_{c,J,b})}{\partial x_{1}}(0,s)|ds=K$$

and such that 

$$t_{0}Ccln(J)e^{-C(K+C_{1}b^{\frac14}t_{0}+C_{2}t_{0})}>K.$$

The result then follows by contradiction, since if we assume $t_{0}=t^{crit}_{t_{0},K,c,J,b}$ we obtain

$$\int_{0}^{t_{0}}|\frac{\partial v_{1}(\theta_{c,J,b})}{\partial x_{1}}(0,s)|ds\geq t_{0}Ccln(J)e^{-C(K+C_{1}b^{\frac14}t_{0}+C_{2}t_{0})}>K. $$

\end{proof}

\begin{corollary}\label{corolariocondinic}
There are initial conditions $\theta^{initial}_{K,t_{0},\tilde{c}}\in H^{4}$  with $||\theta^{initial}_{K,t_{0},\tilde{c}}||_{H^{2}}\leq \tilde{c}$ such that there exists $0<t^{crit}_{K,t_{0},\tilde{c}}\leq t_{0}$ and a solution $\theta_{K,t_{0},\tilde{c}}(x,t)$ to (\ref{SQG}) with $\theta^{initial}_{K,t_{0},\tilde{c}}$ as initial conditions fulfilling
$$ \int_{0}^{t^{crit}_{K,t_{0},\tilde{c}}}\frac{\partial v_{1}(\theta_{K,t_{0},\tilde{c}})}{\partial x_{1}}(0,s)ds=-K,$$
$$||\theta_{K,t_{0},\tilde{c}}(x,t)||_{H^{4}}\leq M_{K,t_{0},\tilde{c}}.$$

Furthermore we have $supp(\theta^{initial}_{K,t_{0},\tilde{c}})\subset\{r\in(a_{1},\frac32)\}$, $supp(\theta_{K,t_{0},\tilde{c}}(x,t))\subset\{r\in(a_{1},a_{2})\}$ with $a_{1}$, $a_{2}$ depending on $K$, $t_{0}$ and $\tilde{c}$.

\end{corollary}

\begin{proof}
The initial conditions and solution are the ones obtained in lemma \ref{estiramientogrande}, we only need to note that $||\theta_{c,J,b}||_{H^{2}}= c(\sum_{j=1}^{J}\frac{1}{j^2})^{\frac12}\leq Cc$, and thus we need to take $Cc\leq\tilde{c}$ and then apply lemma \ref{estiramientogrande}. As for the condition regarding the support, we just need to use that since the solution remains in $H^{4}$ the velocity is $C^1$ and that the velocity at $(x_{1},x_{2})=(0,0)$ is zero and thus particles can only approach the origin exponentially fast.
\end{proof}

\begin{theorem}\label{teoremacrit}
For any $c_{0}>0$, $M>1$ and $t_{*}>0$, we can find a $H^{2+\frac14}$ function $\theta_{0}(x)$  with $||\theta_{0}(x)||_{H^{2}}\leq c_{0}$ such that the only solution $\theta(x,t)\in H^{2+\frac14}$ to the SQG equation (\ref{SQG}) with initial conditions $\theta_{0}(x)$ is such that there exists $t\leq t^{*}$ with $||\theta(x,t)||_{H^{2}}\geq M  c_{0}$.
\end{theorem}

\begin{proof}
We consider the pseudo-solution
\begin{align}\label{bartheta}
    &\bar{\theta}_{M,t^{*},c_{0},N}=\theta_{K=4M,t_{0}=t^{*},\tilde{c}=\frac{c_{0}}{2}}(x,t)\\
    &+\frac{c_{0}}{4}g_{1}(e^{G(t)}N^{\frac12}x_{1})g_{2}(e^{-G(t)}N^{\frac12}x_{2})\frac{sin(e^{G(t)}Nx_{1})}{N^{\frac32}}\nonumber \\ \nonumber
\end{align}

with $\theta_{K,t_{0},\tilde{c}}$ given by corollary \ref{corolariocondinic} with $\tilde{c}=\frac{c_{0}}{2}$, $t_{0}=t^{*}$ and $K=4M$,
$$G(t)=-\int_{0}^{t}\frac{\partial v_{1}(\theta_{K,t_{0},\tilde{c}})}{\partial x_{1}}(0,s)ds$$
and $g_{1}(x_{1}),g_{2}(x_{2})$  $C^{\infty}$ functions with support in $[-1,1]$ and $||g_{i}||_{L^{2}}=1$.
We will define
 $$f^{1}_{M,t^{*},c_{0}}(x,t):=\theta_{K=4M,t_{0}=t^{*},\tilde{c}=\frac{c_{0}}{2}}(x,t)$$
 $$f^{2}_{c_{0},N}(x,t):=\frac{c_{0}}{4}g_{1}(e^{G(t)}N^{\frac12}x_{1})g_{2}(e^{-G(t)}N^{\frac12}x_{2})\frac{sin(e^{G(t)}Nx_{1})}{N^{\frac32}}$$
 
 for a more compact notation.

These pseudo-solutions have the following properties:
\begin{itemize}
    \item For $N$ large, $||\bar{\theta}_{M,t^{*},c_{0},N}(t=0)||_{H^{2}}\leq c_{0}$.
    \item There exists a $t_{crit}\leq t^{*}$ (given by corollary \ref{corolariocondinic}) such that, for $N$ large, we have
    
    $$||\bar{\theta}_{M,t^{*},c_{0},N}(t=t_{crit})||_{H^{2}}\geq \frac{c_{0}}{16} e^{16M}>c_{0}e^{M}$$
    where we used that, since $g_{1},g_{2}\in C^{1}$ and have compact support, for $\lambda>0$
    $$\text{lim}_{N\rightarrow\infty}||N^{\frac12}g_{1}(\lambda N^{\frac12}x_{1})g_{2}(\lambda N^{\frac12}x_{2})sin(\lambda Nx_{1})||_{L^{2}}=\frac{1}{\sqrt{2}}||g(x_{1})||_{L^{2}}.$$
\end{itemize}

Furthermore they fulfil the evolution equation
\begin{align*}
    &\frac{\bar{\theta}_{M,t^{*},c_{0},N} }{\partial t} + v_{1}(f^{1}_{M,t^{*},c_{0}})\frac{\partial f^{1}_{M,t^{*},c_{0}}}{\partial x_{1}} + v_{2}(f^{1}_{M,t^{*},c_{0}})\frac{\partial f^{1}_{M,t^{*},c_{0}}}{\partial x_{2}}\\
    &x_{1}\frac{\partial v_{1}(f^{1}_{M,t^{*},c_{0}})}{\partial x_{1}}\frac{\partial f^{2}_{c_{0},N}}{\partial x_{1}}+ x_{2}\frac{\partial v_{2}(f^{1}_{M,t^{*},c_{0}})}{\partial x_{2}}\frac{\partial f^{2}_{c_{0},N}}{\partial x_{2}}= 0\\
\end{align*}

and thus it is a pseudo-solution with source term

\begin{align*}
    &F_{M,t^{*},c_{0},N}(x,t)=F^{1}_{M,t^{*},c_{0},N}(x,t)+F^{2}_{M,t^{*},c_{0},N}(x,t)+F^{3}_{M,t^{*},c_{0},N}(x,t),\\
\end{align*}

$$F^{1}_{M,t^{*},c_{0},N}(x,t):=-(v_{1}(f^{2}_{c_{0},N})\frac{\partial f^{2}_{c_{0},N}}{\partial x_{1}} + v_{2}(f^{2}_{c_{0},N})\frac{\partial f^{2}_{c_{0},N}}{\partial x_{2}})$$

$$F^{2}_{M,t^{*},c_{0},N}(x,t):=-(v_{1}(f^{2}_{c_{0},N})\frac{\partial f^{1}_{M,t^{*},c_{0}}}{\partial x_{1}} + v_{2}(f^{1}_{M,t^{*},c_{0}})\frac{\partial f^{2}_{c_{0},N}}{\partial x_{2}})$$
\begin{align*}
    &F^{3}_{M,t^{*},c_{0},N}(x,t):=(x_{1}\frac{\partial v_{1}(f^{1}_{M,t^{*},c_{0}})}{\partial x_{1}}-v_{1}(f^{1}_{M,t^{*},c_{0}}))\frac{\partial f^{2}_{c_{0},N}}{\partial x_{1}}\\
    &+ (x_{2}\frac{\partial v_{2}(f^{1}_{M,t^{*},c_{0}})}{\partial x_{2}}-v_{2}(f^{1}_{M,t^{*},c_{0}}))\frac{\partial f^{2}_{c_{0},N}}{\partial x_{2}}.\\
\end{align*}

As usual we want to find bounds for the source term for $t\in[0,t_{crit}]$. For $F^{1}_{M,t^{*},c_{0},N}(x,t)$ it is easy to obtain that

$$||F^{1}_{M,t^{*},c_{0},N}(x,t)||_{L^{2}}\leq CN^{-\frac52},\quad ||F^{1}_{M,t^{*},c_{0},N}(x,t)||_{H^{3}}\leq CN^{\frac12}$$

with $C$ depending on $M$ and $c_{0}$.

For $F^{2}_{M,t^{*},c_{0},N}(x,t)$, using that $||f^{2}_{c_{0},N}||_{L^{1}}\leq CN^{-\frac52}$ and that the support of $f^{1}_{M,t^{*},c_{0}}$ lies away from $0$, we get

$$||F^{2}_{M,t^{*},c_{0},N}(x,t)||_{L^{2}}\leq CN^{-\frac{5}{2}},\quad ||F^{2}_{M,t^{*},c_{0},N}(x,t)||_{H^{3}}\leq CN^{-\frac{5}{2}}$$

with $C$ depending on $M$, $t^{*}$ and $c_{0}$.

Finally, for $F^{3}_{M,t^{*},c_{0},N}(x,t)$, using that, for $i=1,2$ 
$$x_{i}\frac{\partial v_{i}(f^{1}_{M,t^{*},c_{0}})}{\partial x_{1}}-v_{i}(f^{1}_{M,t^{*},c_{0}})$$
vanishes to second order around $0$, that the third derivatives of $v_{i}(f^{1}_{M,t^{*},c_{0}})$ are bounded around $0$, and $supp(f^{2}_{c_{0},N})\subset[-N^{-\frac12},N^{-\frac12}]\times[-N^{-\frac12},N^{-\frac12}]$, we get

$$||F^{3}_{M,t^{*},c_{0},N}||_{L^{2}}\leq CN^{-\frac72},\quad ||F^{3}_{M,t^{*},c_{0},N}||_{H^{3}}\leq CN^{-\frac12},$$

with $C$ depending on $M$, $t^{*}$ and $c_{0}$.

With all this combined and using the interpolation inequality, we get

$$||F_{M,t^{*},c_{0},N}||_{L^{2}}\leq CN^{-\frac52},\quad ||F_{M,t^{*},c_{0},N}||_{H^{2+\frac14}}\leq CN^{-\frac14}.$$

This allows us to obtain, in a similar way as in lemmas \ref{errorl2}, \ref{errorhs}, \ref{errorl2sobolev} and \ref{errorhssobolev} that, if $\theta_{M,t^{*},c_{0},N}(x,t)$ is the solution to (\ref{SQG}) with $\theta_{M,t^{*},c_{0},N}(x,0)=\bar{\theta}_{M,t^{*},c_{0},N}(x,0)$ then

$$||\theta_{M,t^{*},c_{0},N}(x,t)-\bar{\theta}_{M,t^{*},c_{0},N}(x,t)||_{H^{2+\frac14}}\leq CtN^{-\frac14}$$

and this combined with the properties of $\bar{\theta}_{M,t^{*},c_{0},N}(x,t)$ finishes the proof.

\end{proof}

\begin{theorem}\label{teoremacritnon}
For any $c_{0}>0$ there exist initial conditions $\theta(x,0)$ with $||\theta(x,0)||_{H^{2}}$ $\leq c_{0}$ such that any solution $\theta(x,t)$ to (\ref{SQG}) satisfies

$$\text{ess-sup}_{t\in[0,\epsilon]}||\theta(x,t)||_{H^{2}}=\infty$$
for any $\epsilon>0$.
\end{theorem}

\begin{proof}

We start by fixing some arbitrary $c_{0}>0$ and defining

$$\bar{\theta}_{n,R,N}(x,t):=T_{R}(\bar{\theta}_{M=4^{n},t^{*}=2^{-n},c_{0}=2^{-n},N}),$$

with $\bar{\theta}_{M,t^{*},c_{0},N}$ as in (\ref{bartheta}) and $T_{R}(f(x_{1},x_{2}))=f(x_{1}+R,x_{2})$. We will also refer to the first  time when
$$||\bar{\theta}_{n,R,N}(x,t)||_{H^{2}}\geq 2^{n}$$
(which we already know exists and is smaller than $2^{-n}$) as $t_{crit,n}.$

We will study the initial conditions

\begin{equation}\label{condicionesinicialesfinales}
    \theta(x,0)=\sum_{n=1}^{\infty}\bar{\theta}_{n,R_{n},N_{n}}(x,0),
\end{equation}

which fulfil  $||\theta(x,0)||_{H^{2}}\leq c_{0}$ if each $N_{n}$ is big enough, and we will prove by contradiction that if we choose appropriately $(R_{n})_{n\in\mathds{N}}$ and $(N_{n})_{n\in \mathds{N}}$ there cannot exists a solution $\theta(x,t)$ with this initial conditions that satisfies 
\begin{equation}\label{pepsilon}
    \text{ess-sup}_{t\in[0,\epsilon]}||\theta(x,t)||_{H^{2}}\leq P
\end{equation}

for some $\epsilon$, $P$. Note also that $\bar{\theta}_{n,R_{n},N_{n}}(x,0)$ is supported in $B_{\frac32}(-R_{n},0)$.  We can assume that our $L^{2}$ norm is conserved, since this will be true if equation (\ref{pepsilon}) holds (for the time intervals that we will consider).  We will assume without loss of generality that $\epsilon\leq 1$, and we define  $v_{max}$ as the maximum velocity that a function  $f$ with $||f||_{H^{2}}\leq 1$, $||f||_{L^{2}}\leq ||\theta(x,0)||_{L^{2}}$ can produce. With this in mind, we write

$$R_{n}=D_{n}+D_{n+1}+4v_{max}2^{n-1}+R_{n-1}+3$$

with $D_{n}=N^{2}_{n}$ and we will prove that, if $N_{n}$ is big enough, then any solution to (\ref{SQG}) with initial conditions (\ref{condicionesinicialesfinales}) will satisfy

\begin{equation}
    \text{ess-sup}_{t\in[0,2^{-n}]}||\theta(x,t)||_{H^{2}}\geq 2^{n-1}.
\end{equation}

Note that with this definition of $R_{n}$, we have, for any $i\neq n$ that

$$d(\text{supp}(T_{R_{n}}(\bar{\theta}_{n,R_{n},N_{n}}(x,0))),\text{supp}(T_{R_{i}}(\bar{\theta}_{i,R_{i},N_{i}}(x,0))))\geq 4v_{max}2^{n-1}+D_{n}$$

For this, we focus on the evolution of 
$$\theta_{n}(x,t):=1_{B_{D_{n}+2v_{max}2^{n-1}+\frac32}(-R_{n},0)}\theta(x,t)$$

and we will assume that 

\begin{equation}
    \text{ess-sup}_{t\in[0,2^{-n}]}||\theta(x,t)||_{H^{2}} < 2^{n-1}.
\end{equation}

Then if $t\in[0,2^{-n}]$, we have that $\theta_{n}(x,t)$ will fulfil the evolution equation

$$\frac{\partial \theta_{n}}{\partial t} + (v_{1}(\theta_{n})+v_{1}(\theta-\theta_{n}))\frac{\partial \theta_{n}}{\partial x_{1}} + (v_{2}(\theta_{n})+v_{2}(\theta-\theta_{n}))\frac{\partial\theta_{n}}{\partial x_{2}}= 0$$

and $||v_{i}(\theta-\theta_{n}) 1_{B_{v_{max}2^{n}}(R_{n})}||_{C^{3}}\leq CN_{n}^{-4}$ since $d(\text{supp}(\theta-\theta_{n}),\text{supp}(\theta_{n}))\geq N^{2}_{n}$

But then we can argue as in lemmas \ref{errorl2sobolev}, \ref{errorhssobolev} and \ref{evolucionvext} to show that, for $t\in[0,t_{crit,n}]$, if $N_{n}$ is large, we have that 

$$||\theta_{n}(x,t)-T_{R_{n}}(\bar{\theta}_{n,R_{n},N_{n}}(x,t))||_{H^{2+\frac14}}\leq CN_{n}^{-\frac14}.$$

Since for some $t_{crit,n}\in[0,2^{-n}]$ we have that

$$||T_{R_{n}}(\bar{\theta}_{n,R_{n},N_{n}}(x,t_{crit,n}))||_{H^{2}}\geq 2^{n}, $$
and the $H^{2}$ norm of $T_{R_{n}}(\bar{\theta}_{n,R_{n},N_{n}}(x,t))$ is continuous in time, we arrive to a contradiction by taking $N_{n}$ big enough and repeating this argument for each $n\in \mathds{N}$.

\end{proof}

\begin{remark}\label{remarkfinal}
The proof can be adapted to work in the critical spaces $W^{1+\frac2p,p}$ for $p\in(1,\infty]$. For this, note that it is easy to obtain a version of corollary \ref{corolariocondinic} but with small $W^{1+\frac2p,p}$, since the function

$$\sum_{j=1}^{J}c\frac{f(b^{-j}r)b^{j}sin(2\alpha)}{j }$$

has a $W^{1+\frac2p,p}$ norm as small as we want by taking $c$ small. As for the perturbation, we need to consider

$$\lambda g_{1}(N^{b}x_{1})g_{2}(N^{b}x_{2})\frac{sin(Nx_{1})}{N^{1+a}},$$

with $a=a(p),b=b(p)\geq 0$ values that keep the norm $W^{1+\frac2p,p}$ bounded (but not tending to zero) as $N\rightarrow \infty$ (for example, in $W^{1,\infty}$ we consider $a=0$) and $\lambda>0$. Taking $b=\frac12$ and arguing as in theorems \ref{teoremacrit} and \ref{teoremacritnon} allows us to obtain ill-posedness for a wide range of $p$, but we need to include some refinements to obtain the result for all $p\in(1,\infty]$. Namely, 
approximations for the velocity similar to those obtained in lemma \ref{aproxfinal} are needed and we have to include one extra time dependent term in the pseudo-solution.
\end{remark}

\section*{Acknowledgements}
This work is supported in part by the Spanish Ministry of Science
and Innovation, through the “Severo Ochoa Programme for Centres of Excellence in R$\&$D
(CEX2019-000904-S)” and MTM2017-89976-P. DC and LMZ were partially supported by
the ERC Advanced Grant 788250.

\bibliographystyle{abbrv}

\begin{thebibliography}{93}

\bibitem{Bourgaincm} J. Bourgain and D. Li. Strong illposedness of the incompressible Euler equation in integer $C^m$ spaces. Geom. Funct. Anal.
25 (2015), no. 1, 1–86.

\bibitem{Bourgainsobolev} J. Bourgain and D. Li. Strong ill-posedness of the incompressible Euler equation in borderline
sobolev spaces. Inventiones Mathematicae, 201(1), 2014, 97-157.

\bibitem{Bressan} Bressan, A., Nguyen, K.T.: Global existence of weak solutions for the Burgers–Hilbert equation. SIAM J. Math. Anal. 46(4), 2884–2904, 2014.

\bibitem{Buckmaster} T. Buckmaster, S. Shkoller, and V. Vicol. Nonuniqueness of weak solutions to the SQG equation. Comm. Pure
Appl. Math. Volume 72, Issue 9  (2019) 1809-1874.

\bibitem{growthc1delta} A. Castro, D. Cordoba, and F. Gancedo, Singularity formations for a
surface wave model, Nonlinearity, 23 (2010), pp. 2835–2847.


\bibitem{Globalsmooth} A. Castro, D. Cordoba, and J. Gomez-Serrano. Global smooth solutions for the inviscid SQG equation. Memoirs
of the AMS, (2020) Volume 266, Number 1292.

\bibitem{ccz} A. Castro, D. Cordoba and F. Zheng. The lifespan of classical solutions for the inviscid Surface Quasi-Geostrophic equation. arXiv: 2007.04692.

\bibitem{ChaeWu} D. Chae and J. Wu. Logarithmically regularized inviscid models in borderline Sobolev spaces. J. Math. Phys., 53(11):115601, 15, 2012.

\bibitem{cls} P. Constantin, M.-C. Lai, R. Sharma, Y.-H. Tseng, and J. Wu. New numerical results for the
surface quasi-geostrophic equation. J. Sci. Comput., 50(1):1–28, 2012.


\bibitem{Majda} P. Constantin, A. Majda and E. Tabak. Formation of strong fronts in the 2D quasi-geostrophic thermal
active scalar. Nonlinearity 7, (1994), 1495–1533.

\bibitem{const1} P. Constantin and H. Q. Nguyen. Local and global strong solutions for SQG in bounded domains. Phys. D 376/377 (2018), 195–203.

\bibitem{const2} P. Constantin and H. Q. Nguyen. Global weak solutions for SQG in bounded domains. Comm. Pure Appl. Math. 71 (2018), no. 11, 2323–2333.

\bibitem{c} D. Cordoba. Nonexistence of simple hyperbolic blow-up for the quasi-geostrophic equation. Ann.
of Math. (2), 148(3):1135–1152, 1998.

\bibitem{cf} D. Cordoba and C. Fefferman. Growth of solutions for QG and 2D Euler equations. J. Amer.
Math. Soc., 15(3):665–670, 2002.


\bibitem{Elgindi} T.M. Elgindi and N. Masmoudi. $L^{\infty}$ ill-posedness for a class of equations arising in hydrodynamics. Arch. Ration. Mech. 235(3)  1979-2025, (2020).

\bibitem{Elgindisobolev} T.M. Elgindi and I.J. Jeong. ll-posedness for the Incompressible Euler Equations in Critical Sobolev Spaces. Annals of PDE 3(1), 19pp, (2017).

\bibitem{SmallscaleSQG} S. He and A. Kiselev. Small-scale creation for solutions of the SQG equation. Duke Math. J. 170 (5) 1027 - 1041, (2021).

\bibitem{Held} I. Held, R. Pierrehumbert, S. Garner and K. Swanson. Surface quasi-geostrophic dynamics, J. Fluid Mech.,
282, (1995), 1–20.

\bibitem{Injee} I.J. Jeong, J. Kim, Strong illposedness for SQG in critical Sobolev spaces. arXiv: 2107.07739.

\bibitem{Jolly} M.S. Jolly, A. Kumar and V.R. Martinez. On local well-posedness of logarithmic inviscid regularizations
of generalized SQG equations in borderline Sobolev spaces. Communications on
pure and applied analysis. Volume 21, Number 1, January 2022.

\bibitem{Kukavica} I. Kukavica, V. Vicol, and F. Wang. On the ill-posedness of active scalar equations with odd singular kernels. In New trends in differential equations, control theory and optimization, pages 185–200. World Sci. Publ., Hackensack, NJ, 2016. 

\bibitem{Nazarov} A. Kiselev and F. Nazarov. A simple energy pump for the periodic 2D surface quasi-geostrophic equation,
175–179, Abel Symp., 7, Springer, Heidelberg, 2012.

\bibitem{Kwon} H. Kwon,
Strong ill-posedness of logarithmically regularized 2D Euler equations in the borderline Sobolev space,
Journal of Functional Analysis, Volume 280, Issue 7, 108822, 2021.

\bibitem{Dongli} D. Li. On Kato–Ponce and fractional Leibniz. Rev. Mat. Iberoam. 35 (2019), 23-100.

\bibitem{Marchand} F. Marchand. Existence and regularity of weak solutions to the quasi-geostrophic equations in the spaces
$L^p$ or $\dot{H}^{-1/2}$. Comm. Math. Phys., 277(1):45–67, 2008.

\bibitem{oy} K. Ohkitani and M. Yamada. Inviscid and inviscid-limit behavior of a surface quasigeostrophic flow. Phys. Fluids, 9(4):876–882, 1997.

\bibitem{Ped} J. Pedlosky. 1979 Geophysical Fluid Dynamics. Springer.

\bibitem{Resnic} S. G. Resnick. Dynamical problems in non-linear advective partial differential equations. PhD thesis,
University of Chicago, Department of Mathematics, 1995.

\bibitem{scott} R.K. Scott. A scenario for finite-time singularity in the quasigeostrophic model. Journal of
Fluid Mechanics, 687:492–502, 2011.

\bibitem{Wu} J.Wu. Solutions of the 2D quasi-geostrophic equation in
Hölder spaces. Nonlinear Analysis 62 (2005) 579 – 594.











\end{thebibliography}

\end{document}